\documentclass[11pt]{article}

\usepackage{amsmath,amsthm}
\usepackage{amssymb,bbm}
\usepackage{mathabx}
 \usepackage{IEEEtrantools}
\usepackage{mathrsfs}
\usepackage{color}
\usepackage{tikz-cd}
\usetikzlibrary{arrows,decorations.markings}
 \usetikzlibrary{fit,matrix}
\usepackage{multirow}

\usepackage{authblk}

\numberwithin{equation}{section}

\newcommand{\mychi}{\raisebox{.4ex}[1ex][0.5ex]{$\chi$}}

\counterwithin*{equation}{section}

\newtheorem{dfn}{Definition}[section]
\newtheorem{thm}[dfn]{Theorem}

\newtheorem{lmma}[dfn]{Lemma}

\newtheorem{ppsn}[dfn]{Proposition}

\newtheorem{crlre}[dfn]{Corollary}

\newtheorem{xmpl}[dfn]{Example}

\newtheorem{rmrk}[dfn]{Remark}

\newcommand{\bdfn}{\begin{dfn}}
\newcommand{\bthm}{\begin{thm}}
\newcommand{\blmma}{\begin{lmma}}
\newcommand{\bppsn}{\begin{ppsn}}
\newcommand{\bcrlre}{\begin{crlre}}
\newcommand{\bxmpl}{\begin{xmpl}}
\newcommand{\brmrk}{\begin{rmrk}}

\newcommand{\edfn}{\end{dfn}}
\newcommand{\ethm}{\end{thm}}
\newcommand{\elmma}{\end{lmma}}
\newcommand{\eppsn}{\end{ppsn}}
\newcommand{\ecrlre}{\end{crlre}}
\newcommand{\exmpl}{\end{xmpl}}
\newcommand{\ermrk}{\end{rmrk}}
\newcommand{\exrcs}{\end{xrcs}}

\begin{document}

\author[1]{\sc Manabendra Giri\thanks{MG thanks the Indian Statistical Institute 
for financial support in the form of a PhD fellowship.}
\thanks{manabgiri18r@isid.ac.in, manabendra991@gmail.com}}
\author[1]{\sc Arup Kumar Pal\thanks{arup@isid.ac.in, arupkpal@gmail.com}}
\affil[1]{Indian Statistical Institute, Delhi, INDIA}
\title{Irreducible representations of the crystallization of the 
quantized function algebras $C(SU_{q}(n+1))$}

\maketitle

% \dedication{Dedicated to the memory of  K. R. Parthasarathy (1936--2023)}

\begin{center}
\textit{Dedicated to the memory of  K. R. Parthasarathy (1936--2023)}
\end{center}
%------
% Insert your abstract.
%------
\begin{abstract}
Crystallization of the $C^*$-algebras $C(SU_{q}(n+1))$ was introduced  by
Giri \& Pal in \cite{GirPal-2022tv} as a $C^*$-algebra $C(SU_{0}(n+1))$ given by a
finite set of generators and relations. Here we study 
representations of the $C^*$-algebra $C(SU_{0}(n+1))$ and prove a factorization
theorem for its irreducible representations. This leads to a complete
classification of all irreducible representations of this $C^*$-algebra. As an
important consequence, we prove that all the irreducible representations
of $C(SU_{0}(n+1))$ arise exactly as $q\to 0+$ limits of irreducible 
representations of $C(SU_{q}(n+1))$. We also present a few other important
corollaries of the classification theorem.
\end{abstract}
{\bf AMS Subject Classification No.:}
20G42, %Quantum groups (quantized function algebras) and their representations
46L67, %Quantum groups (operator algebraic aspects)
58B32, %Geometry of quantum groups
\\

%------
% INSERT THE BODY OF THE PAPER HERE (except
% acknowledgments, funding info and bibliography)
%------

%--------------------------------------------------------------------
\tableofcontents
%%--------------------------------------------------------------------

%%--------------------------------------------------------------------
 \section{Introduction}
%%--------------------------------------------------------------------
The theory of crystal bases in quantum group theory was developed in the early
1990's independently by  Kashiwara (\cite{Kas-1991aa}, \cite
{Kas-1993aa}), and Lusztig (\cite{Lus-1990aa}, \cite{Lus-1990ab}). In
particular, Kashiwara (\cite{Kas-1993aa}) had produced crystal bases for the
coordinate ring $\mathcal{O}(G_{q})$ of regular functions for the $q$-deformation 
of a complex simple Lie group $G$ viewed as a module over the corresponding quantized
universal enveloping algebra $U_{q}(\mathfrak{g})$ of the Lie algebra
$\mathfrak{g}$ of the group $G$. Loosely speaking, this can be thought of 
as a $q\to 0$ limit of a basis of the module $\mathcal{O}(G_{q})$ in the
algebraic set up. In the present paper, we 
will be concerned with the notion of crystallization of the algebra 
(or *-algebra) structure of $\mathcal{O}(G_{q})$. 
The $q\to 0$ limits of coordinate function algebras on different quantum 
spaces have appeared in the literature, see for example \cite{Wor-1987aa},
\cite{HonSzy-2002aa}, \cite{HonSzy-2003aa}  etc.
But these were defined in a rather ad hoc manner.
A more systematic approach was taken up in \cite{GirPal-2022tv}
where the notion of crystallization for the $C^*$-algebras $C(K_{q})$, 
$K=SU(n+1)$, $n\geq 2$ was introduced by the present authors.
This was followed by the work
of Giselsson \cite{Gis-2023aa} who looked at the $SU(3)$ case and proved 
that $C(SU_{0}(3))$ is the $C^{*}$-algebra of a graph of rank 2, and
Matassa \& Yuncken \cite{MatYun-2023aa} who introduced the notion
of crystallization in a wider set up, including in particular all connected simply 
connected compact Lie groups $K$ and proved that the crystallized 
algebras $C(K_{0})$ are higher rank graph $C^{*}$-algebras, in the process 
setting up a link between the theory of crystal bases and crystallization 
of coordinate function algebras.

The two definitions of crystallization given in \cite{GirPal-2022tv}
and \cite{MatYun-2023aa} are slightly different.
In simple terms, the difference between the two definitions are as follows.
Matassa \& Yuncken (\cite{MatYun-2023aa}) take a specific 
faithful representation $\pi_{q}$ which act on the same Hilbert space 
$\mathcal{H}$ for different values of $q$, and define the $C^{*}$-subalgebra 
of the bounded operators on $\mathcal{H}$ formed by the $q\to 0+$ limits of 
the operators $\pi_{q}(a)$ for certain elements $a\in \mathcal{O}(K_{q})$
to be the crystallization of $C(K_{q})$.
In \cite{GirPal-2022tv} on the other hand, we observe that for every 
irreducible representation of $C(K_{q})$, the canonical generating 
elements, scaled appropriately, obey the same set of relations and 
use those relations to define the crystallized algebra.
In this article, we will study irreducible representations of the 
crystallization of the $C^*$-algebras $C(SU_{q}(n+1))$  as defined in 
\cite{GirPal-2022tv}. In particular, we will classify all irreducible representations of the crystallized algebra.
As one of the consequences, we realize  our crystallized
algebra as a $C^*$-subalgebra of operators on the Hilbert space
on which the Soibelman representation (see Section~3 for the definition)
of the $C^*$-algebras $C(SU_{q}(n+1))$ reside. 

One immediate outcome of the way the crystal limit of 
$C(SU_{q}(n+1))$ is defined in \cite{GirPal-2022tv}
is that for any irreducible representation $\pi_{q}$ of $C(SU_{q}(n+1))$,
taking limits of the scaled generators, one gets a representation
of the crystallized algebra $C(SU_{0}(n+1))$. In the present paper, we prove that
the representations that arise in this way are irreducible, inequivalent, and 
they give all the irreducible representations of $C(SU_{0}(n+1))$.
Recall (Chapter~3, \cite{KorSoi-1998ab}) that in the $q\neq 0$ case,
Soibelman develops the representation theory of $C(K_{q})$ in terms of
highest weight modules in classifying all the irreducible representations of
$C(K_{q})$. It is not clear to us at this point if his techniques can be 
modified and used for classifying the irreducible representations of 
$C(SU_{0}(n+1))$. The proofs presented here use elementary tools involving 
Hilbert space operators.

The classification result that we prove here was proved for the $n=2$ case
in \cite{GirPal-2022tv}. However, the technique of proof used here is different.
For the $n=2$ case, an irreducible representation was classified by extracting the
invariant quantities needed to describe it (these correspond to the different
cases in the proof of the main theorem in \cite{GirPal-2022tv}) and then directly
setting up a unitary equivalence for the two representations. However, for higher
values of $n$, it becomes difficult to make a similar proof work. The proof we 
give here is a recursive one. We prove a factorization theorem for irreducible
representations, and by repeated application of it, we are then able to use 
the classification statement for $n-1$ case, which then gives us the result.

The main motivation behind looking at the crystallizations of the 
$C^{*}$-algebras $C(K_{q})$ and their representations come from
the local index computations of Connes for $SU_{q}(2)$ (\cite{Con-2004ab})
and from a subsequent result in Chakraborty \& Pal  (\cite{ChaPal-2022ey}) 
where they prove an approximate equivalence involving the GNS representation 
of the Haar state of $SU_{q}(2)$ and the Soibelman 
representation of $C(SU_{q}(2))$ on the Hilbert space
$\ell^{2}(\mathbb{Z})\otimes\ell^{2}(\mathbb{N})$.
The main observation in \cite{ChaPal-2022ey} was that this approximate 
equivalence is what is going on structurally behind the computations 
of Connes (\cite{Con-2004ab}) on the equivariant spectral triple 
for $SU_{q}(2)$. To obtain the unitary that gives the approximate 
equivalence, it was observed that `at $q=0$', the above two 
representations give two representations of a $q=0$ version the $C^{*}$-algebra  
$C(SU_{q}(2))$. One gets this $C^*$-algebra if one replaces $q$ by 0 
in the defining relations of $C(SU_{q}(2)$ (see \cite{Wor-1987aa}).
Up to unitary equivalence, the GNS representation at $q=0$ is then an
amplification of the other representation at $q=0$, and this same
unitary then gives an approximate equivalence for nonzero values of $q$.
By Soibelmen's classification result, for all connected simply 
connected compact Lie groups $K$, the $C^{*}$-algebra $C(K_{q})$ admits a faithful
representation on a nice Hilbert space, just like that of $C(SU_{q}(2))$ 
mentioned above (this representation has been referred to as the 
Soibelmen representation in \cite{MatYun-2023aa} -- we will adopt the 
same terminology here).
We  expect that study of the crystallizations of $C(K_{q})$ for other 
compact Lie groups $K$ and their representations will shed more light 
on the relation between GNS representations of the Haar states for these 
quantum groups and their Soibelman representations, thereby leading to a 
better understanding of the equivariant spectral triples for these quantum 
groups given by Neshveyev \& Tuset in \cite{NesTus-2010ab}.
From an application standpoint, the origin of crystallization of 
quantized function algebras can thus be traced back to the computations 
by Connes on equivariant spectral
triples for $SU_{q}(2)$ (and similar computations on equivariant spectral triples
for other homogeneous spaces in \cite{SuiDabLan-2005aa}, \cite{PalSun-2010aa}
etc.\ following Connes' ideas) and the more recent paper \cite{ChaPal-2022ey}
where the role of the algebra at $q=0$ was further clarified. 

It is important to note that for the crystallizations of $C(K_{q})$ to be
useful in the way indicated above, the crystal limits must have the
property that representations of $C(K_{q})$ give rise to representations
of the crystallized algebra $C(K_{0})$ by sending the generators of
$C(K_{0})$ to limits of the (scaled) generators for $C(K_{q})$. Our main 
result in this paper says that in the type $A_{n}$ case that we have
considered, this is the case for irreducible representations of $C(K_{q})$.

Let us now give a brief outline of the content of this paper. In Section~2, we
recall the crystallzation of the $C^{*}$-algebras $C(SU_{q}(n+1))$. 
One consequence of the way the crystallized algebra is defined is
that every irreducible representation of $C(SU_{q}(n+1))$ gives rise to a
representation of the crystallized algebra $C(SU_{0}(n+1))$ in the limit $q\to 0+$.
In Section~3, we prove that each of these representations is irreducible, and
representations arising from limits of two inequivalent irreducible
representations are inequivalent. Sections~4--7 are devoted to proving that
any irreducible representation of $C(SU_{0}(n+1))$ must occur in this way. In
Section~4, we prove that the generating elements $z_{i,j}^{}$ are all partial
isometries. We also prove a few relations that play a very useful role in the
subsequent sections. In Section~5, we study the image $\pi(C(SU_{0}(n+1)))$ and its
closure for a representation $\pi$ of $C(SU_{0}(n+1))$ on a complex separable
Hilbert space. In Section~6, we further specialize to an irreducible
representation $\pi$, which helps us draw stronger conclusions on the
operators we define in Section~5. In particular, we are able to decompose the
operators $\pi(z_{i,j}^{})$ is a certain way (Theorem~\ref
{th:factorization-2}) that helps us prove a factorization theorem for
irreducible representations (Theorem~\ref{th:factorization-1}) in Section~7.
This factorization theorem, in turn, helps us give a recursive proof of the
classification theorem for all irreducible representations of $C(SU_{0}(n+1))$. In
particular, the classification result tells us that any irreducible
representation must occur as a $q\to 0+$ limit of an irreducible
representation for $C(SU_{q}(n+1))$. In the final section, we present this and 
a few other important corollaries of the classification theorem.

It will be good to mention a couple of very important questions 
related to our study that we do not touch upon in this paper.
The first is about the two notions of crystallization given in 
\cite{GirPal-2022tv} and \cite{MatYun-2023aa} mentioned earlier:
whether they result in the same $C^*$-algebra. We make a few remarks on this 
in the last section (see Remark~\ref{rm:comparison}). 
 The other is the possible $q$-invariance of the $C^{*}$-algebras
 arising as quantized function algebras of a Lie group or its homogeneous space.
 This has been investigated in various cases in the literature,
 see for example  \cite{Wor-1987aa}, \cite{HonSzy-2002aa}, 
 \cite{HonSzy-2003aa}, \cite{Nag-1999aa} and \cite{She-1997ac}. 
 For the case of $SU(n+1)$, Giselsson proved in \cite{Gis-2020aa} that 
 the $C^*$-algebras  $C(SU_{q}(n+1)$ are isomorphic for $q\in (0,1)$.
 So the question that remains to be settled in this case is whether
 the crystallized algebra $C(SU_{0}(n+1)$ is isomorphic to the ones 
 for $q\in(0,1)$. The results in the last section, where the crystallized 
 algebra is realized as a $C^*$-algebra of operators on the same space
 as the Soibelmen representation of all the $C(SU_{q}(n+1)$'s,
 might be useful in settling this question.

A word about the notations used for the quantized function algebras. 
For a Lie group $K$, We will denote the algebra of regular functions 
on the $q$-deformation $K_{q}$ by $\mathcal{O}(K_{q})$ and the $C^{*}$-algebra 
of continuous functions on it by $C(K_{q})$. For the group $SU(n+1)$, however,
we deviate from this and use the more widely used notations 
$\mathcal{O}(SU_{q}(n+1))$ and $C(SU_{q}(n+1))$ respectively.
%%--------------------------------------------------------------------

%%--------------------------------------------------------------------
\section{The crystallized algebras $C(SU_{0}(n+1))$ and $\mathcal{O}(SU_{0}(n+1))$}
%%--------------------------------------------------------------------
Let us start by recalling the following theorem from \cite{GirPal-2022tv}.
\bthm[\cite{GirPal-2022tv}] %\texttt{th:an0-defrel}]
    %%---------------------------
    \label{th:an0-defrel}
    %%---------------------------
There is a universal $C^{*}$-algebra generated by elements 
$z_{i,j}^{}$, $1\leq i,j\leq n+1$ satisfying the following relations:
{\allowdisplaybreaks
\begin{IEEEeqnarray}{rCl}
	 z_{i,j}^{}z_{i,l}^{} &=& 0 \qquad\qquad \qquad\quad \; 
	                     \text{if }j<l,\label{eq:an0-defrel-1}\\
	 z_{i,j}^{}z_{k,j}^{} &=& 0 
	 \qquad\qquad \qquad\quad \; \text{if }i<k,\label{eq:an0-defrel-2}\\
      z_{i,l}^{}z_{k,j}^{} - z_{k,j}^{}z_{i,l}^{} &=& 0 
      \qquad\qquad \qquad\quad \;  \text{if }i<k\text{ and }j<l.
      \label{eq:an0-defrel-3}\\
     z_{i,l}^{}z_{k,j}^{}  &=&0
         \qquad\qquad \quad\quad  \;
            \begin{cases}
		     \text{if } i<k,\; j<l\\
	          \text{and } \max\{i,j\}\geq \min\{k,l\}.
			\end{cases}\label{eq:an0-defrel-4}\\
z_{i,j}^{}z_{k,l}^{} - z_{k,l}^{}z_{i,j}^{}&=& z_{i,l}^{}z_{k,j}^{}
    \qquad\quad \quad  
            \begin{cases}
		     \text{if } i<k,\; j<l\\
	          \text{and } \max\{i,j\}+1 = \min\{k,l\},
			\end{cases}\label{eq:an0-defrel-5}\\
    &&\nonumber\\
  z_{i,j}^{}z_{k,l}^{} - z_{k,l}^{}z_{i,j}^{}  &=& 0
    \qquad\qquad \quad\quad \; 
            \begin{cases}
		     \text{if } i<k,\; j<l\\
	          \text{and } \max\{i,j\}+1 < \min\{k,l\},
			\end{cases}\label{eq:an0-defrel-6}\\[1ex]
    z_{1,1}^{}z_{2,2}^{}\ldots z_{n+1,n+1}^{} &=& 1,\label{eq:an0-defrel-7}\\[1ex]
    z_{i,j}^{}z_{r,s}^{*} - z_{r,s}^{*}z_{i,j}^{} 
        &=& 0,\qquad\qquad \qquad\quad \;  \text{if }i\neq r\text{ and }j\neq s,
        \label{eq:an0-defrel-8}
%%--------------------------------
\end{IEEEeqnarray}
}
\begin{IEEEeqnarray}{rCl}
	z_{r,s}^{*}
	  &=& \begin{cases}
	 (z_{1,1}^{}\dots z_{s-1,s-1}^{})(z_{s,s+1}^{}z_{s+1,s+2}^{}
    \dots z_{r-1,r}^{})(z_{r+1,r+1}^{}\dots z_{n+1,n+1}^{})&\text{if } r>s,\\
	     (z_{1,1}^{}\dots z_{r-1,r-1}^{})(z_{r+1,r}^{}z_{r+2,r+1}^{}
        \dots z_{s,s-1}^{})(z_{s+1,s+1}^{}\dots z_{n+1,n+1}^{}) &\text{if } r<s,\\
	    (z_{1,1}^{}\dots z_{s-1,s-1}^{})(z_{s+1,s+1}^{}z_{s+2,s+2}^{}
       \dots  z_{n+1,n+1}^{})  &\text{if } r=s.
		\end{cases}
    %%---------------------------
    \label{eq:an0-defrel-9}
    %%---------------------------
\end{IEEEeqnarray}
\ethm
%%---------------------------------------------------------
%%---------------------------------------------------------
The $C^{*}$-algebra in the above theorem is defined to be 
the crystallization of the  family of $C^{*}$-algebras 
$C(SU_{q}(n+1))$ in \cite{GirPal-2022tv}. This will be denoted by
$C(SU_{0}(n+1))$.
The *-algebra given by the above set of relations is similarly defined as 
the crystallization of the quantized algebra of regular functions
$\mathcal{O}(SU_{q}(n+1))$. We will denote this *-algebra by 
$\mathcal{O}(SU_{0}(n+1))$.

%%---------------------------------------------------------
%%---------------------------------------------------------
The following is an easy consequence of the relations satisfied by the 
generators of the algebra $\mathcal{O}(SU_{0}(n+1))$.
%%---------------------------------------------------------
\bppsn %[\texttt{pr:an0-xtrrel-1}]
    %%---------------------------
    \label{pr:an0-xtrrel-1}
    %%---------------------------
Let $1\leq i\leq j \leq n+1$. Then
\begin{IEEEeqnarray}{rCl}
 (z_{i,i}^{}z_{i+1,i+1}^{}\ldots z_{j,j}^{})^{*}
    &=& (z_{1,1}^{}z_{2,2}^{}\ldots z_{i-1,i-1}^{})
           (z_{j+1,j+1}^{}z_{j+2,j+2}^{}\ldots z_{n+1,n+1}^{}).
    %%---------------------------
    \label{eq:xtrrel-10}
    %%---------------------------
\end{IEEEeqnarray}
\eppsn
%%---------------------------------------------------------
\begin{proof}
Note that for $i=j$, it follows from (\ref{eq:an0-defrel-9}).
We will assume that the equality holds for $i\leq j<n+1$ and 
prove that it holds for $i\leq j+1$ also. Using (\ref{eq:an0-defrel-6})
and (\ref{eq:an0-defrel-9}), we get
\begin{IEEEeqnarray*}{rCl}
 \IEEEeqnarraymulticol{3}{l}{
    (z_{i,i}^{}z_{i+1,i+1}^{}\ldots z_{j+1,j+1}^{})^{*}}\\
   \qquad\qquad &=& 
    z_{j+1,j+1}^{*}\left(z_{i,i}^{}z_{i+1,i+1}^{}\ldots z_{j,j}^{}\right)^{*}\\
    &=& (z_{1,1}^{}z_{2,2}^{}\ldots z_{j,j}^{})
           (z_{j+2,j+2}^{}z_{j+3,j+3}^{}\ldots z_{n+1,n+1}^{})\\
      &&\qquad     (z_{1,1}^{}z_{2,2}^{}\ldots z_{i-1,i-1}^{})
           (z_{j+1,j+1}^{}z_{j+2,j+2}^{}\ldots z_{n+1,n+1}^{})\\
    &=& (z_{j+2,j+2}^{}z_{j+3,j+3}^{}\ldots z_{n+1,n+1}^{})
            (z_{1,1}^{}z_{2,2}^{}\ldots z_{j,j}^{}) \\
      &&\qquad   (z_{j+1,j+1}^{}z_{j+2,j+2}^{}\ldots z_{n+1,n+1}^{})
                   (z_{1,1}^{}z_{2,2}^{}\ldots z_{i-1,i-1}^{})\\
    &=&  (z_{j+2,j+2}^{}z_{j+3,j+3}^{}\ldots z_{n+1,n+1}^{})
                   (z_{1,1}^{}z_{2,2}^{}\ldots z_{i-1,i-1}^{})\\
    &=& (z_{1,1}^{}z_{2,2}^{}\ldots z_{i-1,i-1}^{})
           (z_{j+2,j+2}^{}z_{j+3,j+3}^{}\ldots z_{n+1,n+1}^{}).
\end{IEEEeqnarray*}
This completes the proof.
\end{proof}
%%--------------------------------------------------------------------

%%--------------------------------------------------------------------
\section{$q\to 0+$ limits of representations}
%%--------------------------------------------------------------------
In this section we will be concerned with representations of 
 $C(SU_{0}(n+1))$) that arise as $q\to 0+$ limits of the irreducible 
 representations of the $C^{*}$-algebra $C(SU_{q}(n+1))$. 
We will start by recalling briefly the 
irreducible representations  $C(SU_{q}(n+1))$ from \cite{KorSoi-1998ab} where
irreducible representations of $C(K_{q})$ for all connected simply 
connected compact Lie groups $K$ were classified.
It should be mentioned in this context that following Soibelman's ideas,
the irreducible representations of the quantized function algebras
for more general classes of spaces were later studied and classified
by Dijkhuizen \& Stokman \cite{StoDij-1999aa} 
and Neshveyev \& Tuset \cite{NesTus-2012ab}. The results in these
papers will have an important role to play in
studying crystallizations of quantized function algebras in
those cases, exactly like Soibelman's results play in the present case. 

The Weyl group for $SU_q(n+1)$ is isomorphic to the permutation group 
$\mathfrak{S}_{n+1}$ on $n+1$ symbols. Denote by
$s_i$ the transposition $(i,i+1)$. Then $\{s_{1},s_{2},\ldots,s_{n}\}$ form a
set of generators for $\mathfrak{S}_{n+1}$. 
For $a,b\in\mathbb{N}$ with $a\leq b$, let us denote by $s_{[a,b]}$ the product
$s_{b}s_{b-1}\ldots s_{a}$ in $\mathfrak{S}_{n+1}$.
Let $1\leq k\leq n$,
$1\leq b_{k}<b_{k-1}<\cdots <b_{1}\leq n$, and $1\leq a_{i}\leq b_{i}$ for 
$1\leq i\leq k$.
It follows from the strong exchange condition and the deletion condition in the characterization 
of Coxeter system (see~\cite{BjoBre-2005aa}, \cite{Gar-1997aa}) that 
$\omega=s_{[a_{k},b_{k}]}s_{[a_{k-1},b_{k-1}]}\ldots  s_{[a_{1},b_{1}]}$
is a reduced word in $\mathfrak{S}_{n+1}$. It is well-known that if $k=n$,
$a_{i}=1$  and $b_{i}=i$ for $1\leq i\leq n$, then one gets the longest element
$\omega_{0}$.

Let $S$  be the left shift operator and $N$ be the number 
operator on $\ell^{2}(\mathbb{N})$:
\[
Se_n=\begin{cases}
       e_{n-1}& \text{if } n\geq 1,\cr
       0 & \text{if }n=0.
	  \end{cases}, \qquad Ne_{k}=ke_{k},\quad k\in\mathbb{N}.
\]
Let us denote the generators of the $C^{*}$-algebra $C(SU_{q}(n+1))$ by
$u_{i,j}^{}(q)$. Denote by $\psi_{s_r}^{(q)}$ the following 
representation of $C(SU_{q}(n+1))$ on $\ell^{2}(\mathbb{N})$:
%%---------------------------------------------------------
\begin{equation}
\psi_{s_r}^{(q)}(u_{i,j}^{}(q))=\begin{cases}
                   S\sqrt{I-q^{2N}} &\text{if }i=j=r,\cr
                  \sqrt{I-q^{2N}}S^{*} & \text{if } i=j=r+1,\cr
                   -q^{N+1} &  \text{if } i=r, j=r+1,\cr
                    q^{N}   &  \text{if }i=r+1, j=r,\cr
                    \delta_{i,j}I & \text{otherwise}.
					  \end{cases}
	%%---------------------------
	\label{eq:irr-rep-1}
	%%---------------------------
\end{equation}
%%---------------------------------------------------------
For a reduced word $\omega=s_{i_1}s_{i_2}\ldots
s_{i_m}\in\mathfrak{S}_{n+1}$, define $\psi_{\omega}^{(q)}$ to be
$\psi_{s_{i_1}}^{(q)}\ast\psi_{s_{i_2}}^{(q)}\ast\ldots\ast\psi_{s_{i_m}}^{
(q)}$. Here, for two representations $\phi$ and $\psi$, $\phi\ast\psi$
denotes the representation $(\phi\otimes\psi)\Delta_{q}$, where 
$\Delta_{q}$ is the comultiplication map on $C(SU_{q}(n+1))$.

Next, let $\lambda\equiv(\lambda_{1},\ldots,\lambda_{n})\in (S^1)^{n}$.
Define
%%---------------------------------------------------------
\begin{equation}
	\mychi_{\lambda}(u_{i,j}^{}(q))=
	\begin{cases}
            \lambda_{i}\delta_{i,j} &  \text{if }i=1,\cr
            \bar{\lambda}_{n} \delta_{i,j} &  \text{if }i=n+1,\cr
            \bar{\lambda}_{i-1}\lambda_{i}\delta_{i,j} &  \text{otherwise}.
	\end{cases}
	%%---------------------------
	\label{eq:irr-rep-2}
	%%---------------------------
\end{equation}
%%---------------------------------------------------------
A well-known result of Soibelman (Theorem 6.2.7, page~121, \cite{KorSoi-1998ab})
says that for any
reduced word $\omega\in\mathfrak{S}_{n+1}$ and a tuple $\lambda\equiv
(\lambda_{1},\ldots,\lambda_{n})\in (S^1)^{n}$, the representation
\begin{equation}
  \psi_{\lambda,\omega}^{(q)}:=\mychi_{\lambda}\ast \psi_{\omega}^{(q)}
    %%---------------------------
       \label{eq:irr-rep-3}
    %%---------------------------
\end{equation}
is an irreducible representation of $C(SU_{q}(n+1))$, and these give all the
irreducible representations of $C(SU_{q}(n+1))$.
If one takes the longest word $\omega_{0}$ in $\mathfrak{S}_{n+1}$ and takes
the direct integral of the representations $\psi_{\lambda,\omega_{0}}^{(q)}$
with respect to $\lambda$, one gets a faithful representation $\psi^{(q)}$ of
$C(SU_{q}(n+1))$ acting on 
$\ell^{2}(\mathbb{Z})^{\otimes n}\otimes 
\ell^{2}(\mathbb{N})^{\otimes\frac{n(n+1)}{2}}$
(or equivalently on 
$L^{2}(S^{1})^{\otimes n}\otimes \ell^{2}(\mathbb{N})^{\otimes\frac{n(n+1)}{2}}$).
Following \cite{MatYun-2023aa}, this will be called the \textbf{Soibelman representation}.

Let us define
\begin{IEEEeqnarray}{rCl}
	Z_{i,j}^{} &=& \lim_{q\to 0+}(-q)^{\min\{i-j,0\}}
	 \psi_{\lambda,\omega}^{(q)}(u_{i,j}^{}(q)),
	%%---------------------------
	   \label{eq:irred-q0-1}
	%%---------------------------
\end{IEEEeqnarray}
where $\psi_{\lambda,\omega}^{(q)}$ is the representation of $C(SU_{q}(n+1))$ defined
above. It has been shown in~\cite{GirPal-2022tv} 
(Propositions~2.1, 2.3, 2.4, 2.6 and~2.7, \cite{GirPal-2022tv}) that the 
above limits exist and the operators $Z_{i,j}^{}$ obey the relations
(\ref{eq:an0-defrel-1}--\ref{eq:an0-defrel-9}). Therefore
\begin{IEEEeqnarray}{rCl}
\psi_{\lambda,\omega}^{}(z_{i,j}^{})=
 \lim_{q\to 0+}(-q)^{\min\{i-j,0\}}
 \psi_{\lambda,\omega}^{(q)}(u_{i,j}^{}(q)),\qquad
  i,j\in\{1,2,\ldots,n+1\}
	%%---------------------------
	\label{eq:irred-q0-2}
	%%---------------------------
\end{IEEEeqnarray}
defines a representation of $C(SU_{0}(n+1))$ on the Hilbert space
$\ell^{2}(\mathbb{N}^{\ell(\omega)})$, where $\ell(\omega)$ denotes the length
of the element $\omega$. 
In particular, if $\omega$ is the identity element, its reduced word is
the empty word and the representation $\psi_{\lambda,\omega}^{}$ in
this case is an one dimensional representation given by
%%---------------------------------------------------------
\begin{equation}
   \psi_{\lambda,id}^{}(z_{i,j}^{})=
   \begin{cases}
            \lambda_{i}\delta_{i,j} &  \text{if }i=1,\cr
            \bar{\lambda}_{n} \delta_{i,j} &  \text{if }i=n+1,\cr
            \bar{\lambda}_{i-1}\lambda_{i}\delta_{i,j} &  \text{otherwise}.
   \end{cases}
   %%---------------------------
   \label{eq:irred-q0-2aa}
   %%---------------------------
\end{equation}
%%---------------------------------------------------------
If $\lambda=(1,1,\ldots,1)$, we will denote 
$\psi_{\lambda,\omega}^{}$ by just $\psi_{\omega}^{}$.
The remaining part of this section is devoted to the proof of
irreducibility of $\psi_{\lambda,\omega}^{}$ and their 
inequivalence for different pairs $(\lambda,\omega)$.

Let us start with the observation that the representation $\psi_{s_{r}}^{}$
is given by
%%---------------------------------------------------------
\begin{equation}
\psi_{s_r}^{}(z_{i,j}^{})=\begin{cases}
                   S &\text{if }i=j=r,\cr
                  S^{*} & \text{if } i=j=r+1,\cr
                   P_{0} &  \text{if } i\neq j \text{ and } i,j\in\{r, r+1\},\cr
                    \delta_{i,j}I & \text{otherwise},
                      \end{cases}
    %%---------------------------
    \label{eq:irred-q0-2a}
    %%---------------------------
\end{equation}
where $P_{0}$ denotes the rank one projection $|e_{0}\rangle\langle e_{0}|=I-S^{*}S$.
%%---------------------------------------------------------
The following lemma is  
an important observation that will be used in the proof
of the main theorem in the last section. The proof is straightforward
and hence omitted.
%%---------------------------------------------------------
\blmma %[\texttt{lm:irred-q0-0}]
    %%---------------------------
    \label{lm:irred-q0-0}
    %%---------------------------
Let 
$\omega=s_{[a_{k},b_{k}]}s_{[a_{k-1},b_{k-1}]}\ldots s_{[a_{1},b_{1}]}$ and let
$\omega^{\prime}= \omega s_{a_{1}-1}$. Then
%%---------------------------------------------------------
\begin{IEEEeqnarray}{rCl}
   \psi_{\lambda,\omega'}(z_{i,j}^{}) &=& 
      \sum_{m=\min\{i,j\}}^{\max\{i,j\}}
        \psi_{\lambda,\omega}(z_{i,m}^{})\otimes \psi_{s_{a_{1}-1}}(z_{m,j}^{}).
    %%---------------------------
    \label{eq:irred-q0-3}
    %%---------------------------
\end{IEEEeqnarray}
%%---------------------------------------------------------
\elmma
%%---------------------------------------------------------

%%---------------------------------------------------------
%%---------------------------------------------------------
  \paragraph{Diagrams for the representations $\psi_{\lambda,\omega}$.}
  We will next introduce certain diagrams for the representations
$\psi_{\lambda,\omega}$ similar to those introduced in \cite{ChaPal-2003ac} for
the case of $SU_{q}(n+1)$ for nonzero $q$. These diagrams will be very
helpful in understanding the images of the generatoring elements
$z_{i,j}^{}$ of the $C^*$-algebra $C(SU_{q}(n+1))$. In particular, they will
play a very important role in proving the irreducibility of the representations
$\psi_{\lambda,\omega}$ in the present section and later on while characterising
all irreducible representations, these diagrams will once again he helpful 
by providing us the necessary ideas and intuition.

Let us start with the diagrams for the 
representations $\psi_{s_{i}}$ and $\psi_{s_{i}s_{j}}$.

%%---------------------------------------
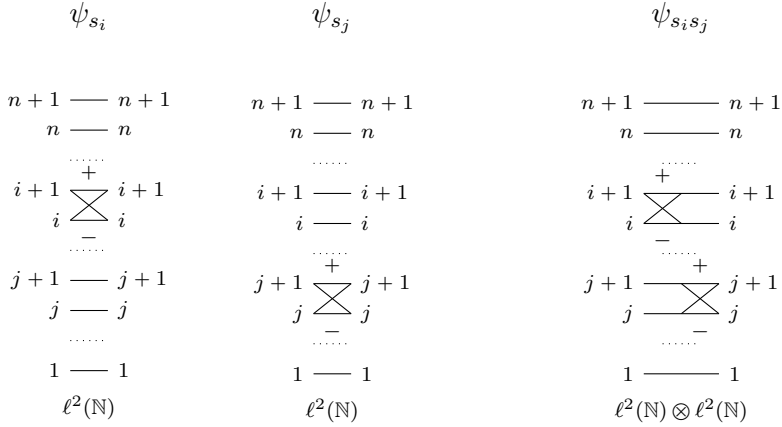
\begin{tabular}{p{80pt}p{80pt}p{160pt}}
%%---------------------------------------
% \input{An0-irred-diagram-2a.tex}
%%---------------------------------------
\begin{center}
\begin{tikzpicture}[xscale=.5, yscale=.4]
	\draw[] (2.5,11) node[above]{$\psi_{s_{i}}$};
	\draw[] (0,0) +(2,9)  -- +(3,9);
	\draw[] (0,0) +(2,8)  -- +(3,8);
	\draw[dotted] (0,0) +(2,7)  -- +(3,7);
	\draw[] (0,0) +(2,6)  -- +(3,6);
	\draw[] (0,0) +(2,5)  -- +(3,5);
	\draw[dotted] (0,0) +(2,4)  -- +(3,4);
	\draw[] (0,0) +(2,3)  -- +(3,3);
	\draw[] (0,0) +(2,2)  -- +(3,2);
	\draw[dotted] (0,0) +(2,1)  -- +(3,1);
	\draw[] (0,0) +(2,0)  -- +(3,0);
	\draw[] (2,5) +(0,0)  -- +(1,1);
	\draw[] (2,6) +(0,0)  -- +(1,-1);
	\draw[] (2.5,6) node[above]{\scriptsize $+$};
	\draw[] (2.5,5) node[below]{\scriptsize $-$};
	\draw[] (2.5,-.5) node[below]{\scriptsize $\ell^{2}(\mathbb{N})$};
	\draw[] (2,0) node[left]{\scriptsize 1};
	\draw[] (2,2) node[left]{\scriptsize $j$};
	\draw[] (2,3) node[left]{\scriptsize $j+1$};
	\draw[] (2,5) node[left]{\scriptsize $i$};
	\draw[] (2,6) node[left]{\scriptsize $i+1$};
	\draw[] (2,8) node[left]{\scriptsize $n$};
	\draw[] (2,9) node[left]{\scriptsize $n+1$};
	\draw[] (3,0) node[right]{\scriptsize 1};
	\draw[] (3,2) node[right]{\scriptsize $j$};
	\draw[] (3,3) node[right]{\scriptsize $j+1$};
	\draw[] (3,5) node[right]{\scriptsize $i$};
	\draw[] (3,6) node[right]{\scriptsize $i+1$};
	\draw[] (3,8) node[right]{\scriptsize $n$};
	\draw[] (3,9) node[right]{\scriptsize $n+1$};
\end{tikzpicture}
\end{center}
%%---------------------------------------
&
%%---------------------------------------
% \input{An0-irred-diagram-2b.tex}
%%---------------------------------------
\begin{center}
\begin{tikzpicture}[xscale=.5, yscale=.4]
	\draw[] (2.5,11) node[above]{$\psi_{s_{j}}$};
	\draw[] (0,0) +(2,9)  -- +(3,9);
	\draw[] (0,0) +(2,8)  -- +(3,8);
	\draw[dotted] (0,0) +(2,7)  -- +(3,7);
	\draw[] (0,0) +(2,6)  -- +(3,6);
	\draw[] (0,0) +(2,5)  -- +(3,5);
	\draw[dotted] (0,0) +(2,4)  -- +(3,4);
	\draw[] (0,0) +(2,3)  -- +(3,3);
	\draw[] (0,0) +(2,2)  -- +(3,2);
	\draw[dotted] (0,0) +(2,1)  -- +(3,1);
	\draw[] (0,0) +(2,0)  -- +(3,0);
	\draw[] (2,2) +(0,0)  -- +(1,1);
	\draw[] (2,3) +(0,0)  -- +(1,-1);
	\draw[] (2.5,3) node[above]{\scriptsize $+$};
	\draw[] (2.5,2) node[below]{\scriptsize $-$};
	\draw[] (2.5,-.5) node[below]{\scriptsize $\ell^{2}(\mathbb{N})$};
	\draw[] (2,0) node[left]{\scriptsize 1};
	\draw[] (2,2) node[left]{\scriptsize $j$};
	\draw[] (2,3) node[left]{\scriptsize $j+1$};
	\draw[] (2,5) node[left]{\scriptsize $i$};
	\draw[] (2,6) node[left]{\scriptsize $i+1$};
	\draw[] (2,8) node[left]{\scriptsize $n$};
	\draw[] (2,9) node[left]{\scriptsize $n+1$};
	\draw[] (3,0) node[right]{\scriptsize 1};
	\draw[] (3,2) node[right]{\scriptsize $j$};
	\draw[] (3,3) node[right]{\scriptsize $j+1$};
	\draw[] (3,5) node[right]{\scriptsize $i$};
	\draw[] (3,6) node[right]{\scriptsize $i+1$};
	\draw[] (3,8) node[right]{\scriptsize $n$};
	\draw[] (3,9) node[right]{\scriptsize $n+1$};
\end{tikzpicture}
\end{center}
%%---------------------------------------
&
%%---------------------------------------
% \input{An0-irred-diagram-2c.tex}
%%---------------------------------------
\begin{center}
\begin{tikzpicture}[xscale=.5, yscale=.4]
	\draw[] (3,11) node[above]{$\psi_{s_{i}s_{j}}$};
	\draw[] (0,0) +(2,9)  -- +(3,9) -- +(4,9);
	\draw[] (0,0) +(2,8)  -- +(3,8) -- +(4,8);
	\draw[dotted] (0,0) +(2.5,7)  -- +(3.5,7);
	\draw[] (0,0) +(2,6)  -- +(3,6)  -- +(4,6);
	\draw[] (0,0) +(2,5)  -- +(3,5)  -- +(4,5);
	\draw[dotted] (0,0) +(2.5,4)  -- +(3.5,4);
	\draw[] (0,0) +(2,3)  -- +(3,3)  -- +(4,3);
	\draw[] (0,0) +(2,2)  -- +(3,2)  -- +(4,2);
	\draw[dotted] (0,0) +(2.5,1)  -- +(3.5,1);
	\draw[] (0,0) +(2,0)  -- +(3,0)  -- +(4,0);
	\draw[] (2,5) +(0,0)  -- +(1,1);
	\draw[] (2,6) +(0,0)  -- +(1,-1);
	\draw[] (3,2) +(0,0)  -- +(1,1);
	\draw[] (3,3) +(0,0)  -- +(1,-1);
	\draw[] (2.5,6) node[above]{\scriptsize $+$};
	\draw[] (2.5,5) node[below]{\scriptsize $-$};
	\draw[] (3.5,3) node[above]{\scriptsize $+$};
	\draw[] (3.5,2) node[below]{\scriptsize $-$};
	\draw[] (3,-.5) node[below]{\scriptsize $\ell^{2}(\mathbb{N})\otimes\ell^{2}(\mathbb{N})$};
	\draw[] (2,0) node[left]{\scriptsize 1};
	\draw[] (2,2) node[left]{\scriptsize $j$};
	\draw[] (2,3) node[left]{\scriptsize $j+1$};
	\draw[] (2,5) node[left]{\scriptsize $i$};
	\draw[] (2,6) node[left]{\scriptsize $i+1$};
	\draw[] (2,8) node[left]{\scriptsize $n$};
	\draw[] (2,9) node[left]{\scriptsize $n+1$};
	\draw[] (4,0) node[right]{\scriptsize 1};
	\draw[] (4,2) node[right]{\scriptsize $j$};
	\draw[] (4,3) node[right]{\scriptsize $j+1$};
	\draw[] (4,5) node[right]{\scriptsize $i$};
	\draw[] (4,6) node[right]{\scriptsize $i+1$};
	\draw[] (4,8) node[right]{\scriptsize $n$};
	\draw[] (4,9) node[right]{\scriptsize $n+1$};
\end{tikzpicture}
\end{center}
%%---------------------------------------
\end{tabular}
%%---------------------------------------

In the leftmost diagram above, each path from a node $k$ on the 
left to a node $l$ on the right stands for an 
operator on $\ell^{2}(\mathbb{N})$. A horizontal
unlabelled line stands for the identity operator, 
a horizontal line labelled with a `$+$' sign 
       stands for the right shift $S^*$ and 
one labelled with a `$-$' sign stands for 
the left shift $S$. Diagonal
lines going upward or downward represent the projection $P_{0}$.
The image $\psi_{s_i}(z_{m,l}^{})$ is the operator represented by 
the path from $m$ on the left
to $l$ on the  right,
and is zero if there is no such path.
Similarly the diagram in the middle stands for the
representation $\psi_{s_{j}}$.

The diagram on the right above is for the representation $\psi_{s_{i}s_{j}}$,
where one simply puts the two diagrams representing $\psi_{s_i}$ and $\psi_{s_j}$
adjacent to each other, with the one for $\psi_{s_{j}}$ to the right of that for
$\psi_{s_i}$. 
The image $\psi_{s_{i}s_{j}}(z_{m,l}^{})$ of $z_{m,l}^{}$ is
an operator on $\ell^{2}(\mathbb{N})\otimes \ell^{2}(\mathbb{N})$ 
given by \textit{the
sum of the paths from the node $m$ on the extreme left to the node $l$ on the 
extreme right that do not contain both upward and downward diagonals}. 
It would be zero if there is no such path.

Next, we come to the description of $\psi_\omega$.
As we have already remarked, we will take a reduced word for $\omega$
of the form
$\omega=s_{[a_{k},b_{k}]}s_{[a_{k-1},b_{k-1}]}\ldots  s_{[a_{1},b_{1}]}$
where
$1\leq k\leq n$,
$1\leq b_{k}<b_{k-1}<\cdots <b_{1}\leq n$, and $1\leq a_{i}\leq b_{i}$ for 
$1\leq i\leq k$.
To get the diagram for $\psi_\omega$, we simply put
the diagrams for the $\psi_{s_{r}}$'s appearing in this word
adjacent to each other in the order they appear in the word from left to right. 
Thus for example, if 
 $\omega=s_{[2,4]}s_{[6,7]}s_{[1,8]}s_{[5,9]}s_{[3,10]}$, 
 then the following is the diagram for  $\psi_\omega$:
%%---------------------------------------
 % \input{An0-irred-diagram-3a.tex}
 %%---------------------------------------
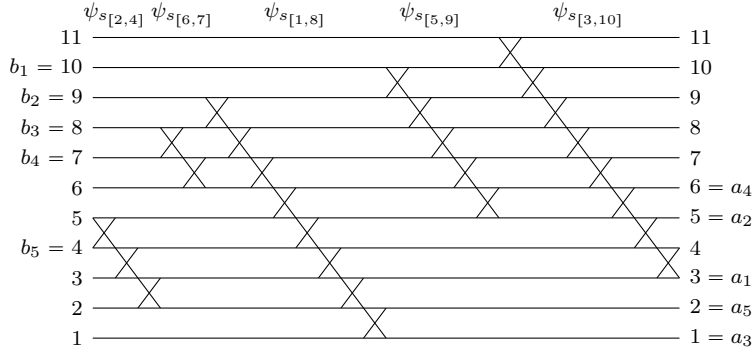
 \begin{center}
\begin{tikzpicture}[xscale=.3, yscale=.4]
	\draw[] (24,10) node[above]{\scriptsize $\psi_{s_{[3,10]}}$};
	\draw[] (17,10) node[above]{\scriptsize $\psi_{s_{[5,9]}}$};
	\draw[] (11,10) node[above]{\scriptsize $\psi_{s_{[1,8]}}$};
	\draw[] (6,10) node[above]{\scriptsize $\psi_{s_{[6,7]}}$};
	\draw[] (3,10) node[above]{\scriptsize $\psi_{s_{[2,4]}}$};
	\draw[] (0,0) +(2,10)  -- +(28,10);
	\draw[] (0,0) +(2,9)  -- +(28,9);
	\draw[] (0,0) +(2,8)  -- +(28,8);
	\draw[] (0,0) +(2,7)  -- +(28,7);
	\draw[] (0,0) +(2,6)  -- +(28,6);
	\draw[] (0,0) +(2,5)  -- +(28,5);
	\draw[] (0,0) +(2,4)  -- +(28,4);
	\draw[] (0,0) +(2,3)  -- +(28,3);
	\draw[] (0,0) +(2,2)  -- +(28,2);
	\draw[] (0,0) +(2,1)  -- +(28,1);
	\draw[] (0,0) +(2,0)  -- +(28,0);
	\draw[] (20,10) +(0,0)  -- +(8,-8);
	\draw[] (20,9) +(0,0)  -- +(1,1);
	\draw[] (21,8) +(0,0)  -- +(1,1);
	\draw[] (22,7) +(0,0)  -- +(1,1);
	\draw[] (23,6) +(0,0)  -- +(1,1);
	\draw[] (24,5) +(0,0)  -- +(1,1);
	\draw[] (25,4) +(0,0)  -- +(1,1);
	\draw[] (26,3) +(0,0)  -- +(1,1);
	\draw[] (27,2) +(0,0)  -- +(1,1);
	\draw[] (15,9) +(0,0)  -- +(5,-5);
	\draw[] (15,8) +(0,0)  -- +(1,1);
	\draw[] (16,7) +(0,0)  -- +(1,1);
	\draw[] (17,6) +(0,0)  -- +(1,1);
	\draw[] (18,5) +(0,0)  -- +(1,1);
	\draw[] (19,4) +(0,0)  -- +(1,1);
	\draw[] (7,8) +(0,0)  -- +(8,-8);
	\draw[] (7,7) +(0,0)  -- +(1,1);
	\draw[] (8,6) +(0,0)  -- +(1,1);
	\draw[] (9,5) +(0,0)  -- +(1,1);
	\draw[] (10,4) +(0,0)  -- +(1,1);
	\draw[] (11,3) +(0,0)  -- +(1,1);
	\draw[] (12,2) +(0,0)  -- +(1,1);
	\draw[] (13,1) +(0,0)  -- +(1,1);
	\draw[] (14,0) +(0,0)  -- +(1,1);
	\draw[] (5,7) +(0,0)  -- +(2,-2);
	\draw[] (5,6) +(0,0)  -- +(1,1);
	\draw[] (6,5) +(0,0)  -- +(1,1);
	\draw[] (2,4) +(0,0)  -- +(3,-3);
	\draw[] (2,3) +(0,0)  -- +(1,1);
	\draw[] (3,2) +(0,0)  -- +(1,1);
	\draw[] (4,1) +(0,0)  -- +(1,1);
	\draw[] (2,0) node[left]{\scriptsize 1};
	\draw[] (2,1) node[left]{\scriptsize 2};
	\draw[] (2,2) node[left]{\scriptsize 3};
	\draw[] (2,3) node[left]{\scriptsize $b_{5}=4$};
	\draw[] (2,4) node[left]{\scriptsize 5};
	\draw[] (2,5) node[left]{\scriptsize 6};
	\draw[] (2,6) node[left]{\scriptsize $b_{4}=7$};
	\draw[] (2,7) node[left]{\scriptsize $b_{3}=8$};
	\draw[] (2,8) node[left]{\scriptsize $b_{2}=9$};
	\draw[] (2,9) node[left]{\scriptsize $b_{1}=10$};
	\draw[] (2,10) node[left]{\scriptsize 11};
	\draw[] (28,0) node[right]{\scriptsize $1=a_{3}$};
	\draw[] (28,1) node[right]{\scriptsize $2=a_{5}$};
	\draw[] (28,2) node[right]{\scriptsize $3=a_{1}$};
	\draw[] (28,3) node[right]{\scriptsize 4};
	\draw[] (28,4) node[right]{\scriptsize $5=a_{2}$};
	\draw[] (28,5) node[right]{\scriptsize $6=a_{4}$};
	\draw[] (28,6) node[right]{\scriptsize 7};
	\draw[] (28,7) node[right]{\scriptsize 8};
	\draw[] (28,8) node[right]{\scriptsize 9};
	\draw[] (28,9) node[right]{\scriptsize 10};
	\draw[] (28,10) node[right]{\scriptsize 11};
\end{tikzpicture}
\end{center}
 %%---------------------------------------
 Finally, for $\lambda=(\lambda_{1},\ldots,\lambda_{n})\in S^{1}$,
 the diagram for $\psi_{\lambda,\omega}$ is the same as
 that for $\psi_{\omega}$, except that \textit{the paths starting at $m$
  on the left get multiplied by a scalar} which is
 $\overline{\lambda}_{m-1}\lambda_{m}$ if $1\leq m\leq n+1$ (with the convention that
 $\lambda_{0}=\lambda_{n+1}=1$).
%%---------------------------------------
\brmrk
  %%---------------------------------------
  \label{rm:diag-1}
  %%---------------------------------------
 In view of the above, we will not need to keep track
 of the $\lambda_{i}$'s except for the last theorem in this section 
 (Theorem~\ref{th:irred-q0-4}) where we prove that $\psi_{\lambda,\omega}$'s are
 inequivalent for different pairs $(\lambda,\omega)$. Accordingly, 
 all the operator equalities we write
 before Theorem~\ref{th:irred-q0-4} in this section will be valid only modulo
 multiplication by a complex number of modulus one.
 \ermrk
%%---------------------------------------
 In the diagram for $\psi_{\lambda,\omega}$, we will refer to the 
 part corresponding to $\psi_{s_{[a_{j},b_{j}]}}$ 
 as the \textbf{$s_{[a_{j},b_{j}]}$-section} and the
 smallest rectangle containing  the
 diagonal lines in the $s_{[a_{j},b_{j}]}$-section as the
 \textbf{$s_{[a_{j},b_{j}]}$-rectangle}.
 For example, the different $s_{[a_{j},b_{j}]}$-rectangles have been 
 shown in red colour and the $s_{[5,9]}$-section has been
 shown in blue colour in the next diagram.
 Observe that if $i>b_{j}+1$, then the horizontal line starting 
 at node $i$ on the left passes above the $s_{[a_{m},b_{m}]}$-rectangles
 for $m\geq j$.
 For an operator $R$ given by a path, we will denote by $R_{j}$ the operator
 given by the $s_{[a_{j},b_{j}]}$-section of that path. In such a case one has
%%---------------------------------------------------------
\begin{IEEEeqnarray*}{rCl}
R &=& R_{k}\otimes R_{k-1}\otimes\ldots\otimes R_{1},
\end{IEEEeqnarray*}
%%---------------------------------------------------------
where $R_{j}$ acts on $\ell^{2}(\mathbb{N})^{\otimes (b_{j}-a_{j}+1)}$ for all $j$.
%%---------------------------------------
 % \input{An0-irred-diagram-3b.tex}
%%---------------------------------------
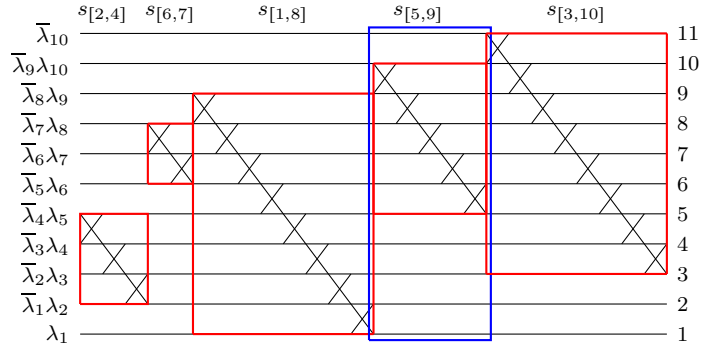
\begin{center}
\begin{tikzpicture}[xscale=.3, yscale=.4]
	\draw[] (24,10) node[above]{\scriptsize $s_{[3,10]}$};
	\draw[] (17,10) node[above]{\scriptsize $s_{[5,9]}$};
	\draw[] (11,10) node[above]{\scriptsize $s_{[1,8]}$};
	\draw[] (6,10) node[above]{\scriptsize $s_{[6,7]}$};
	\draw[] (3,10) node[above]{\scriptsize $s_{[2,4]}$};
	\draw[] (0,0) +(2,10)  -- +(28,10);
	\draw[] (0,0) +(2,9)  -- +(28,9);
	\draw[] (0,0) +(2,8)  -- +(28,8);
	\draw[] (0,0) +(2,7)  -- +(28,7);
	\draw[] (0,0) +(2,6)  -- +(28,6);
	\draw[] (0,0) +(2,5)  -- +(28,5);
	\draw[] (0,0) +(2,4)  -- +(28,4);
	\draw[] (0,0) +(2,3)  -- +(28,3);
	\draw[] (0,0) +(2,2)  -- +(28,2);
	\draw[] (0,0) +(2,1)  -- +(28,1);
	\draw[] (0,0) +(2,0)  -- +(28,0);
	\draw[] (20,10) +(0,0)  -- +(8,-8);
	\draw[] (20,9) +(0,0)  -- +(1,1);
	\draw[] (21,8) +(0,0)  -- +(1,1);
	\draw[] (22,7) +(0,0)  -- +(1,1);
	\draw[] (23,6) +(0,0)  -- +(1,1);
	\draw[] (24,5) +(0,0)  -- +(1,1);
	\draw[] (25,4) +(0,0)  -- +(1,1);
	\draw[] (26,3) +(0,0)  -- +(1,1);
	\draw[] (27,2) +(0,0)  -- +(1,1);
	\draw[] (15,9) +(0,0)  -- +(5,-5);
	\draw[] (15,8) +(0,0)  -- +(1,1);
	\draw[] (16,7) +(0,0)  -- +(1,1);
	\draw[] (17,6) +(0,0)  -- +(1,1);
	\draw[] (18,5) +(0,0)  -- +(1,1);
	\draw[] (19,4) +(0,0)  -- +(1,1);
	\draw[] (7,8) +(0,0)  -- +(8,-8);
	\draw[] (7,7) +(0,0)  -- +(1,1);
	\draw[] (8,6) +(0,0)  -- +(1,1);
	\draw[] (9,5) +(0,0)  -- +(1,1);
	\draw[] (10,4) +(0,0)  -- +(1,1);
	\draw[] (11,3) +(0,0)  -- +(1,1);
	\draw[] (12,2) +(0,0)  -- +(1,1);
	\draw[] (13,1) +(0,0)  -- +(1,1);
	\draw[] (14,0) +(0,0)  -- +(1,1);
	\draw[] (5,7) +(0,0)  -- +(2,-2);
	\draw[] (5,6) +(0,0)  -- +(1,1);
	\draw[] (6,5) +(0,0)  -- +(1,1);
	\draw[] (2,4) +(0,0)  -- +(3,-3);
	\draw[] (2,3) +(0,0)  -- +(1,1);
	\draw[] (3,2) +(0,0)  -- +(1,1);
	\draw[] (4,1) +(0,0)  -- +(1,1);
	\draw[red,thick] (2,1) +(0,0)  -- +(0,3) -- +(3,3) -- +(3,0) -- +(0,0);
	\draw[red,thick] (5,5) +(0,0)  -- +(0,2);
	\draw[red,thick] (5,5) +(0,2)  -- +(2,2);
	\draw[red,thick] (5,5) +(2,2)  -- +(2,0);
	\draw[red,thick] (5,5) +(2,0)  -- +(0,0);
	\draw[red,thick] (7,0) +(0,0)  -- +(0,8);
	\draw[red,thick] (7,0) +(0,8)  -- +(8,8);
	\draw[red,thick] (7,0) +(8,8)  -- +(8,0);
	\draw[red,thick] (7,0) +(8,0)  -- +(0,0);
	\draw[red,thick] (15,4) +(0,0)  -- +(0,5);
	\draw[red,thick] (15,4) +(0,5)  -- +(5,5);
	\draw[red,thick] (15,4) +(5,5)  -- +(5,0);
	\draw[red,thick] (15,4) +(5,0)  -- +(0,0);
	\draw[blue,thick] (14.8,-.2) +(0,0)  -- +(0,10.4)  -- +(5.4,10.4) -- +(5.4,0) -- +(0,0);
	\draw[red,thick] (20,2) +(0,0)  -- +(0,8);
	\draw[red,thick] (20,2) +(0,8)  -- +(8,8);
	\draw[red,thick] (20,2) +(8,8)  -- +(8,0);
	\draw[red,thick] (20,2) +(8,0)  -- +(0,0);
	\draw[] (2,0) node[left]{\scriptsize $\lambda_{1}$};
	\draw[] (2,1) node[left]{\scriptsize $\overline{\lambda}_{1}\lambda_{2}$};
	\draw[] (2,2) node[left]{\scriptsize $\overline{\lambda}_{2}\lambda_{3}$};
	\draw[] (2,3) node[left]{\scriptsize $\overline{\lambda}_{3}\lambda_{4}$};
	\draw[] (2,4) node[left]{\scriptsize $\overline{\lambda}_{4}\lambda_{5}$};
	\draw[] (2,5) node[left]{\scriptsize $\overline{\lambda}_{5}\lambda_{6}$};
	\draw[] (2,6) node[left]{\scriptsize $\overline{\lambda}_{6}\lambda_{7}$};
	\draw[] (2,7) node[left]{\scriptsize $\overline{\lambda}_{7}\lambda_{8}$};
	\draw[] (2,8) node[left]{\scriptsize $\overline{\lambda}_{8}\lambda_{9}$};
	\draw[] (2,9) node[left]{\scriptsize $\overline{\lambda}_{9}\lambda_{10}$};
	\draw[] (2,10) node[left]{\scriptsize $\overline{\lambda}_{10}$};
	\draw[] (28,0) node[right]{\scriptsize 1};
	\draw[] (28,1) node[right]{\scriptsize 2};
	\draw[] (28,2) node[right]{\scriptsize 3};
	\draw[] (28,3) node[right]{\scriptsize 4};
	\draw[] (28,4) node[right]{\scriptsize 5};
	\draw[] (28,5) node[right]{\scriptsize 6};
	\draw[] (28,6) node[right]{\scriptsize 7};
	\draw[] (28,7) node[right]{\scriptsize 8};
	\draw[] (28,8) node[right]{\scriptsize 9};
	\draw[] (28,9) node[right]{\scriptsize 10};
	\draw[] (28,10) node[right]{\scriptsize 11};
\end{tikzpicture}
\end{center}
%%---------------------------------------
The diagram for $\psi_{\lambda,\omega}$ and the terminology introduced 
above will be helpful in understanding the operators in the image 
$\psi_{\lambda,\omega}(C(SU_{0}(n+1)))$ as well as certain quantities 
associated with them that we will define in the remaining part of this section. 
They are also going to be helpful in getting insight into some of the proofs 
in sections 5 and 6,
though we will not explicitly use or refer to the diagrams there.
%%---------------------------------------------------------

Given a reduced word 
$\omega=s_{[a_{k},b_{k}]}s_{[a_{k-1},b_{k-1}]}\ldots s_{[a_{1},b_{1}]}$,
we will next define certain elements in the $C^*$-algebra 
$C(SU_{0}(n+1))$ that will play
a crucial role in the proof of irreducibility of the representation
$\psi_{\lambda,\omega}$.
For $2\leq j\leq k$ and $a_{j}\leq i\leq b_{j}+1$, define two integers 
$n(j,i)$ and $n'(j,i)$  as follows:
%%---------------------------------------------------------
\begin{IEEEeqnarray}{rCl}
   n(j,i) &=& \max\bigl\{m: 1\leq m\leq j-1 \text{ and } a_{m}\leq i\bigr\},
   %%------------------------
   \label{eq:irred-q0-3a}\\
   %%------------------------
   n'(j,i) &=& \max\bigl\{m: 1\leq m\leq j-1 \text{ and } a_{m}+1 \leq i\bigr\},
   %%------------------------
   \label{eq:irred-q0-3b}
   %%------------------------
\end{IEEEeqnarray}
%%---------------------------------------------------------
with the usual convention that maximum of an empty set will be 
taken as $-\infty$. 
Thus $n(j,i)$ is the maximum of all those $m\in\{1,2,\ldots,j-1\}$ for which the 
$i$\raisebox{.4ex}{th} horizontal line  (from the bottom) 
intersects the $s_{[a_{m},b_{m}]}$-rectangle, and $n(j,i)=-\infty$
when the $i$\raisebox{.4ex}{th} horizontal line is 
below the $s_{[a_{m},b_{m}]}$-rectangles for $1\leq m\leq j-1$.

Next, for $1\leq j\leq k$ and for $a_{j}\leq i\leq b_{j}+1$, 
define elements $v_{j,i}^{}(\omega)\equiv v_{j,i}^{}$ as follows:
%%---------------------------------------------------------
\begin{IEEEeqnarray*}{rCl}
    v_{1,i}^{}  &=& z_{b_{1}+1,i}^{} ,\\
    v_{j,i}^{}  &=& 
         \begin{cases}
          z_{b_{j}+1,i}^{} v_{n(j,i),i+1}^{}, & 
                   \text{if $n(j,i)$ is finite}\\
          z_{b_{j}+1,i}^{} & \text{if } n(j,i)=-\infty
          \end{cases}
          \qquad 2\leq j\leq k.
          %%------------------------
          \label{eq:irred-q0-3c}
          %%------------------------
\end{IEEEeqnarray*}
%%---------------------------------------------------------
Let us denote by $Z_{j,i}^{}$  the operator
$\psi_{\lambda,\omega}^{}(z_{j,i}^{})$ for $1\leq j,i\leq n+1$
and by $V_{j,i}^{}$ the operator
$\psi_{\lambda,\omega}^{}(v_{j,i}^{})$ for $1\leq j\leq k$ and 
$a_{j}\leq i\leq b_{j}+1$.
We will next study the operators $V_{j,i}^{}$.
In general, $Z_{b_{j}+1,i}^{}$ is given by a sum of paths from
$b_{j}+1$ to $i$ going downwards. 
What the next proposition says is that the operator $V_{j,i}^{}$ will be
given by the lowest among these paths (appearing in the sum
for $Z_{b_{j}+1,i}^{}$) with some of its $s_{[a_{r},b_{r}]}$-sections
slightly modified for $1\leq r\leq j-1$ so that at most one shift operator
occurs in each such section. 
%%---------------------------------------------------------
\bppsn %[\texttt{pr:irred-q0-1}]
	%%---------------------------
	\label{pr:irred-q0-1}
	%%---------------------------
 Let $\ell_{r}=b_{r}-a_{r}+1$ for $1\leq r \leq k$.
 Then for $a_{j} \leq i\leq b_{j}+1$, one has
%%---------------------------------------------------------
\begin{IEEEeqnarray}{rCl}
    V_{j,i} 
    &=&  T_{k,i}^{(j)}\otimes T_{k-1,i}^{(j)}\otimes\cdots \otimes T_{1,i}^{(j)},
    %%---------------------------
    \label{eq:irred-q0-4}
    %%---------------------------
\end{IEEEeqnarray}
%%---------------------------------------------------------
    where $T_{r,i}^{(j)}$'s are of the following form:
    \begin{IEEEeqnarray}{ll}
      I^{\otimes {\ell_{r}}} & \text{if }  r > j, 
             \IEEEyesnumber\IEEEyessubnumber* \label{eq:irred-q0-5a}\\[.5ex]
       P_{0}^{\otimes {(b_{r}+1-i)}} \otimes S^{*}\otimes I^{\otimes(i-a_{r}-1)} 
                      & \text{if }  r = j\text{ and } i \geq a_{r}+1,
                      \label{eq:irred-q0-5b} \\[.5ex]
       P_{0}^{\otimes(b_{r}-a_{r}+1)} 
                      & \text{if }  r = j\text{ and } i= a_{r},
                      \label{eq:irred-q0-5c} \\[.5ex]
        P_{0}^{\otimes {(b_{r}-a_{r}-1)}}\otimes I\otimes S^{*}
                     & \text{if }  r < j \text{ and } i= a_{r}+1, 
                     \label{eq:irred-q0-5d}\\[.5ex]
      P_{0}^{\otimes {(b_{r}-i)}}\otimes I\otimes S^{*}\otimes I^{\otimes(i-a_{r}-1)}
                     & \text{if }  r < j \text{ and } a_{r}+1 < i < b_{r}, 
                     \label{eq:irred-q0-5e}\\[.5ex]
        \text{tensor product of $P_{0}$'s and $I$'s}\qquad
                     & \text{if }  r < j \text{ and } i \leq a_{r} 
                     \label{eq:irred-q0-5f}.
    %%---------------------------
    % \label{eq:irred-q0-5}
    %%---------------------------
    \end{IEEEeqnarray}
\eppsn
%%---------------------------------------------------------
\begin{proof}
For $j=1$, it is easy to check from the diagram 
for $\psi_{\lambda,\omega}^{}$ that (\ref{eq:irred-q0-4}) holds 
for $a_{1} \leq i\leq b_{1}+1$. To help follow the argument,
let us keep the following diagram before us in which the 
$s_{[a_{m},b_{m}]}$-section and the the 
$s_{[a_{m'},b_{m'}]}$-section have been enclosed by pairs of coloured vertical lines.
%%---------------------------------------
 % \input{An0-irred-diagram-3c.tex}
%%---------------------------------------
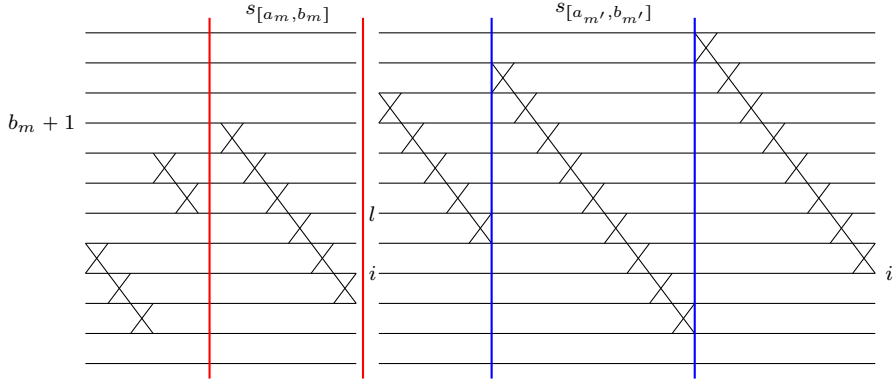
\begin{center}
\begin{tikzpicture}[xscale=.3, yscale=.4]
	\draw[] (9,11) node[above]{\scriptsize $s_{[a_{m},b_{m}]}$};
	\draw[] (23,11) node[above]{\scriptsize $s_{[a_{m'},b_{m'}]}$};
	\draw[] (0,11) +(0,0) -- +(12,0)  +(13,0) -- +(27,0) -- +(35,0);
	\draw[] (0,10) +(0,0) -- +(12,0)  +(13,0) -- +(27,0) -- +(35,0);
	\draw[] (0,9) +(0,0) -- +(12,0)  +(13,0) -- +(27,0) -- +(35,0);
	\draw[] (0,8) +(0,0) -- +(12,0)  +(13,0) -- +(27,0) -- +(35,0);
	\draw[] (0,7) +(0,0) -- +(12,0)  +(13,0) -- +(27,0) -- +(35,0);
	\draw[] (0,6) +(0,0) -- +(12,0)  +(13,0) -- +(27,0) -- +(35,0);
	\draw[] (0,5) +(0,0) -- +(12,0)  +(13,0) -- +(27,0) -- +(35,0);
	\draw[] (0,4) +(0,0) -- +(12,0)  +(13,0) -- +(27,0) -- +(35,0);
	\draw[] (0,3) +(0,0) -- +(12,0)  +(13,0) -- +(27,0) -- +(35,0);
	\draw[] (0,2) +(0,0) -- +(12,0)  +(13,0) -- +(27,0) -- +(35,0);
	\draw[] (0,1) +(0,0) -- +(12,0)  +(13,0) -- +(27,0) -- +(35,0);
	\draw[] (0,0) +(0,0) -- +(12,0)  +(13,0) -- +(27,0) -- +(35,0);
	\draw[] (27,11) +(0,0)  -- +(8,-8);
	\draw[] (27,10) +(0,0)  -- +(1,1);
	\draw[] (28,9) +(0,0)  -- +(1,1);
	\draw[] (29,8) +(0,0)  -- +(1,1);
	\draw[] (30,7) +(0,0)  -- +(1,1);
	\draw[] (31,6) +(0,0)  -- +(1,1);
	\draw[] (32,5) +(0,0)  -- +(1,1);
	\draw[] (33,4) +(0,0)  -- +(1,1);
	\draw[] (34,3) +(0,0)  -- +(1,1);
	\draw[] (18,10) +(0,0)  -- +(9,-9);
	\draw[] (18,9) +(0,0)  -- +(1,1);
	\draw[] (19,8) +(0,0)  -- +(1,1);
	\draw[] (20,7) +(0,0)  -- +(1,1);
	\draw[] (21,6) +(0,0)  -- +(1,1);
	\draw[] (22,5) +(0,0)  -- +(1,1);
	\draw[] (23,4) +(0,0)  -- +(1,1);
	\draw[] (24,3) +(0,0)  -- +(1,1);
	\draw[] (25,2) +(0,0)  -- +(1,1);
	\draw[] (26,1) +(0,0)  -- +(1,1);
	\draw[] (13,9) +(0,0)  -- +(5,-5);
	\draw[] (13,8) +(0,0)  -- +(1,1);
	\draw[] (14,7) +(0,0)  -- +(1,1);
	\draw[] (15,6) +(0,0)  -- +(1,1);
	\draw[] (16,5) +(0,0)  -- +(1,1);
	\draw[] (17,4) +(0,0)  -- +(1,1);
	\draw[] (6,8) +(0,0)  -- +(6,-6);
	\draw[] (6,7) +(0,0)  -- +(1,1);
	\draw[] (7,6) +(0,0)  -- +(1,1);
	\draw[] (8,5) +(0,0)  -- +(1,1);
	\draw[] (9,4) +(0,0)  -- +(1,1);
	\draw[] (10,3) +(0,0)  -- +(1,1);
	\draw[] (11,2) +(0,0)  -- +(1,1);
	\draw[] (3,7) +(0,0)  -- +(2,-2);
	\draw[] (3,6) +(0,0)  -- +(1,1);
	\draw[] (4,5) +(0,0)  -- +(1,1);
	\draw[] (0,4) +(0,0)  -- +(3,-3);
	\draw[] (0,3) +(0,0)  -- +(1,1);
	\draw[] (1,2) +(0,0)  -- +(1,1);
	\draw[] (2,1) +(0,0)  -- +(1,1);
    \draw[red,thick] (5.5,-.5) +(0,0)  -- +(0,12);
    \draw[red,thick] (12.3,-.5) +(0,0)  -- +(0,12);
    \draw[blue,thick] (18,-.5) +(0,0)  -- +(0,12);
    \draw[blue,thick] (27,-.5) +(0,0)  -- +(0,12);
	\draw[] (12.1,5) node[right]{\scriptsize $l$};
	\draw[] (12.1,3) node[right]{\scriptsize $i$};
	\draw[] (0,8) node[left]{\scriptsize $b_{m}+1$};
	\draw[] (35,3) node[right]{\scriptsize $i$};
\end{tikzpicture}
\end{center}
%%---------------------------------------

Let us now assume that (\ref{eq:irred-q0-4}) 
holds for $j\leq m-1$ and $a_{j} \leq i\leq b_{j}+1$.
 We will prove the statement for $j=m$ and $a_{m} \leq i\leq b_{m}+1$.
 Let us look at the case $n(m,i)=-\infty$ first. In this case,
 $V_{j,i}=Z_{b_{m}+1,i}$ and $i<a_{j}$ for all $1\leq j\leq m-1$, i.e.\ the 
 $i$\raisebox{.4ex}{th} horizontal line is 
below the $s_{[a_{j},b_{j}]}$-rectangles for $1\leq j\leq m-1$. 
Therefore $T_{j,i}=I^{\otimes(\ell_{j})}$ for $1\leq j\leq m-1$. 
Since the horizontal line
 starting at $b_{m}+1$ on the left passes above the $s_{[a_{j},b_{j}]}$-rectangles
 for $j>m$, there is exactly one path going from $b_{m}+1$ on the left to $i$ on 
 the right and it does not intersect
 the  $s_{[a_{j},b_{j}]}$-rectangles for  $1\leq j\leq m-1$. Thus
%%---------------------------------------
 \begin{IEEEeqnarray*}{rCl}
    Z_{b_{m}+1,i}^{} &=& I^{\otimes(\sum_{r=m+1}^{k}\ell_{r})}\otimes 
\Bigl(P_{0}^{\otimes(b_{m}+1-i)} \otimes S^{*}\otimes I^{\otimes(i-a_{m}-1)}\Bigr) 
      \otimes       I^{\otimes(\sum_{r=1}^{m-1}\ell_{r})}\\
    &=& T_{k,i}^{(m)}\otimes T_{k-1,i}^{(m)}\otimes\cdots T_{1,i}^{(m)}.
 \end{IEEEeqnarray*}
%%---------------------------------------
  Next assume $m':=n(m,i)$ is finite.
  We have $V_{m,i}=Z_{b_{m}+1,i}^{} V_{m',i+1}^{}$.
  Since $m' < m$, by our hypothesis,
%%---------------------------------------------------------
\begin{IEEEeqnarray*}{rCl}
    V_{m',i+1} 
    &=&  T_{k,i+1}^{(m')}\otimes T_{k-1,i+1}^{(m')}\otimes\cdots T_{1,i+1}^{(m')}.
\end{IEEEeqnarray*}
%%---------------------------------------------------------
Note that since $m' < m$, one has $T_{r,i+1}^{(m')}=I^{\otimes(\ell_{r})}$ 
for $m\leq r\leq k$.
It is also immediate that $Z_{b_{m}+1,i}^{}$ is of the form
$I^{\otimes(\sum_{r=m+1}^{k}\ell_{r})}\otimes R$ for some operator $R$ 
that acts on the last
$\sum_{r=1}^{m}\ell_{r}$ copies of $\ell^{2}(\mathbb{N})$. 
The operator $R$ is given by
$\sum_{l=i}^{b_{m}+1}R(l)$ where $R(l)$ is the sum of operators given 
by the paths that 
go from $b_{m}+1$ to $l$ in the $s_{[a_{m},b_{m}]}$-section and then from
$l$ to $i$ in the $s_{[a_{m-1},b_{m-1}]}\ldots s_{[a_{1},b_{1}]}$-section.

Since the $s_{[a_{r},b_{r}]}$-rectangle for $m'<r<m$ lie above the 
$i$\raisebox{.4ex}{th} horizontal line, we conclude that for $l>i$, the
$s_{[a_{m'},b_{m'}]}$-section of $R(l)$ is of the form 
\[
I^{\otimes(b_{m'}-j)}\otimes S\otimes (\text{tensor product of $P_{0}$, $S^{*}$ and $I$'s}) 
\]
for some $j$ with $i+1\leq j\leq l$, while
$T_{m',i+1}^{(m')}=P_{0}^{\otimes(b_{m'}-i)}\otimes S^{*}
\otimes  I^{\otimes(i-a_{m'})}$.
Therefore 
%%---------------------------------------------------------
\begin{IEEEeqnarray*}{rCl}
R(l)\left(T_{m,i+1}^{(m')}\otimes T_{m-1,i+1}^{(m')}
    \otimes\cdots T_{1,i+1}^{(m')}\right) 
  &=& 0\qquad \text{for }l>i.
\end{IEEEeqnarray*}
%%---------------------------------------------------------
Hence
%%---------------------------------------------------------
\begin{IEEEeqnarray*}{rCl}
R\left(T_{m,i+1}^{(m')}\otimes T_{m-1,i+1}^{(m')}\otimes
      \cdots T_{1,i+1}^{(m')}\right) 
  &=&  R(i)\left(T_{m,i+1}^{(m')}\otimes T_{m-1,i+1}^{(m')}\otimes
      \cdots T_{1,i+1}^{(m')}\right).
\end{IEEEeqnarray*}
%%---------------------------------------------------------
Let us write $T_{r,i}=R(i)_{r}T_{r,i+1}^{(m')}$ for $1\leq r\leq m$.
We will show that $T_{r,i}$ is given by (\ref{eq:irred-q0-5a}--\ref{eq:irred-q0-5f})
for $j=m$.
Note that the path $R(i)$ consists of the diagonal line coming down from
$b_{m}+1$ to $i$ in the $s_{[a_{m},b_{m}]}$-section and its horizonal
continuation in the $s_{[a_{m-1},b_{m-1}]}\ldots s_{[a_{1},b_{1}]}$-section.
Let us now look at the following three cases separately.\\[1ex]
\textbf{Case I}: $1\leq r\leq m-1$ and $a_{r} < i$.\\
Since $a_{r} < i$, we have $r\leq m'$. Therefore in this case, 
$R(i)_{r}=I^{\otimes(b_{r}-i)}\otimes S\otimes S^{*}\otimes I^{\otimes(i-a_{r}-1)}$,
 and
$T_{r,i+1}^{(m')}=P_{0}^{\otimes(b_{r}-i)}
     \otimes S^{*}\otimes I^{\otimes(i-a_{r})}$. Therefore
$T_{r,i}=R(i)_{r}T_{r,i+1}^{(m')}=P_{0}^{\otimes(b_{r}-i)}
   \otimes I\otimes S^{*}\otimes  I^{\otimes(i-a_{r}-1)}$.
Since $r<m$ and $i>a_{r}$, from  (\ref{eq:irred-q0-5e}) it follows that
this is of the form $T_{r,i}^{(m)}$.
\\[1ex]
\textbf{Case II}: $1\leq r\leq m-1$ and $a_{r}= i$.\\
As in the earlier case, we have $r\leq m'$. Hence  
$R(i)_{r}=I^{\otimes(\ell_{r}-1)}\otimes S$.
Now if $r=m'$, then by (\ref{eq:irred-q0-5b}) 
we have
$T_{r,i+1}^{(m')}= P_{0}^{\otimes (\ell_{r}-1)}\otimes S^{*}$
and if $r<m'$, then by (\ref{eq:irred-q0-5d}), we have
$T_{r,i+1}^{(m')}= P_{0}^{\otimes (\ell_{r}-1)}\otimes I\otimes S^{*}$
Thus $T_{r,i}=R(i)_{r}T_{r,i+1}^{(m')} $ is either
$P_{0}^{\otimes (\ell_{r}-1)}\otimes I$ or
$P_{0}^{\otimes (\ell_{r}-2)}\otimes I\otimes I$.
Since $r<m$ and $i=a_{r}$, by (\ref{eq:irred-q0-5f})
$T_{r,i}$ is of the form $T_{r,i}^{(m)}$, as required.\\[1ex]
\textbf{Case III}: $1\leq r\leq m-1$ and $i<a_{r}$.\\
In this case, $R(i)_{r}=I^{\otimes(\ell_{r})}$, so that $T_{r,i}=T_{r,i+1}^{(m')}$.
If $r>m'$, then $T_{r,i+1}^{(m')}= I^{\otimes \ell_{r}}$.
Since  $i+1\leq a_{r}$, if $r=m'$, we must have $i+1=a_{r}$. 
Hence by (\ref{eq:irred-q0-5e}) $T_{r,i+1}^{(m')}$ is a tensor product
of $P_{0}$'s and $I$'s. 
If $r<m'$ and $i+1 \leq a_{r}$, then also by (\ref{eq:irred-q0-5f}), 
$T_{r,i+1}^{(m')}$ is a tensor product
of $P_{0}$'s and $I$'s. 
Since $T_{r,i}=T_{r,i+1}^{(m')}$, in all these subcases, it is of the required form.
\end{proof}
%%---------------------------------------------------------

We now list a few important corollaries of the above proposition.
%%---------------------------------------------------------
\bcrlre %[\texttt{cr:irred-q0-1a}]
    %%---------------------------
    \label{cr:irred-q0-1a}
    %%---------------------------
   For a bounded operator $T$, let $|T|$ denote $(T^{*}T)^{\frac{1}{2}}$. 
   Then one has
%%---------------------------------------------------------
\begin{IEEEeqnarray*}{rCl}
    |V_{j,a_{j}}^{}|\cdot|V_{j-1,a_{j-1}}^{}|\cdots|V_{1,a_{1}}^{}| 
    &=&  T_{k}^{}\otimes T_{k-1}^{}\otimes\cdots \otimes T_{1}^{},
    %%---------------------------
    % \label{eq:irred-q0-1aa}
    %%---------------------------
\end{IEEEeqnarray*}
where
    \begin{IEEEeqnarray*}{rCl}
    T_{r}^{} &=&
       \begin{cases}
      I^{\otimes {\ell_{r}}} & \text{if }  r > j, 
             \IEEEyesnumber\\[.5ex]%\IEEEyessubnumber* \\[.5ex]
       P_{0}^{\otimes {(b_{r}+1-a_{r})}} 
                      & \text{if }  r \leq j.
        \end{cases}
    %%---------------------------
    % \label{eq:irred-q0-1ab}
    %%---------------------------
    \end{IEEEeqnarray*}
%%---------------------------------------------------------
\ecrlre
%%---------------------------------------------------------
%%---------------------------------------------------------
\bcrlre %[\texttt{cr:irred-q0-1b}]
    %%---------------------------
    \label{cr:irred-q0-1b}
    %%---------------------------
 Let  $a_{j} \leq i\leq b_{j}+1$. Assume $n'(j,i)>-\infty$. Then
%%---------------------------------------------------------
\begin{IEEEeqnarray}{rCl}
    V_{n'(j,i),i}^{*}V_{j,i}^{} 
    &=&  T_{k}^{}\otimes T_{k-1}^{}\otimes\cdots \otimes T_{1}^{},
    %%---------------------------
    \label{eq:irred-q0-1ba}
    %%---------------------------
\end{IEEEeqnarray}
%%---------------------------------------------------------
    where 
    \begin{IEEEeqnarray}{rCl}
    T_{r}^{} &=&
       \begin{cases}
      I^{\otimes {(b_{r}+1-a_{r})}} & \text{if }  r > j, 
             \IEEEyesnumber\\[.5ex]%\IEEEyessubnumber* \\[.5ex]
       P_{0}^{\otimes {(b_{r}+1-a_{r})}} 
                      & \text{if }  r = j\text{ and } i = a_{j}, \\[.5ex]
       P_{0}^{\otimes {(b_{r}+1-i)}} \otimes S^{*}\otimes I^{\otimes(i-a_{r}-1)} 
                      & \text{if }  r = j\text{ and } i \geq a_{j}+1, \\[.5ex]
         \text{tensor product of $P_{0}$'s and $I$'s}\qquad
                     & \text{if } r<j.  
        \end{cases}\label{eq:irred-q0-1bc}
    \end{IEEEeqnarray}
\ecrlre
%%---------------------------------------------------------
\begin{proof}
Observe that for $r<j$, the section $T_{r,i}^{(j)}$ is of 
the form (\ref{eq:irred-q0-5d}) or (\ref{eq:irred-q0-5e}) 
(i.e.\ contains the right shift) exactly when $r<j$ and $a_{r}<i$.
Let $n'(j,i)$ be as defined in (\ref{eq:irred-q0-3b}).
Then for the operator $V_{n'(j,i),i}^{}$ also, the corresponding sections
$T_{r,i}^{(s)}$ exactly match $T_{r,i}^{(j)}$ and the remaining ones
are tensor products of $I$ and $P_{0}$'s.
Therefore the result follows.
\end{proof}
%%---------------------------------------------------------
\bcrlre %[\texttt{cr:irred-q0-1c}]
    %%---------------------------
    \label{cr:irred-q0-1c}
    %%---------------------------
 Let  $a_{j} \leq i\leq b_{j}+1$. 
Define
\[
E_{j,i}\equiv E_{j,i}(\omega):=\begin{cases}
      V_{n'(j,i),i}^{*}V_{j,i}^{} 
    |V_{j-1,a_{j-1}}^{}|\cdot|V_{j-2,a_{j-2}}^{}|\cdots|V_{1,a_{1}}^{}| 
              & \text{if } n'(j,i)>-\infty,\\
   V_{j,i}^{} 
    |V_{j-1,a_{j-1}}^{}|\cdot|V_{j-2,a_{j-2}}^{}|\cdots|V_{1,a_{1}}^{}| 
              & \text{if } n'(j,i)= -\infty.
      \end{cases}
\]
 Then
%%---------------------------------------------------------
\begin{IEEEeqnarray}{rCl}
    % \IEEEeqnarraymulticol{3}{l}{
    E_{j,i}(\omega) %:=V_{n'(j,i),i}^{*}V_{j,i}^{} 
    %|V_{j-1,a_{j-1}}^{}|\cdot|V_{j-2,a_{j-2}}^{}|\cdots|V_{1,a_{1}}^{}| 
     &=&  T_{k}^{}\otimes T_{k-1}^{}\otimes\cdots \otimes T_{1}^{},
    %%---------------------------
    \label{eq:irred-q0-1ca}
    %%---------------------------
\end{IEEEeqnarray}
%%---------------------------------------------------------
    where 
    \begin{IEEEeqnarray}{rCl}
    T_{r}^{} &=&
       \begin{cases}
      I^{\otimes {(b_{r}+1-a_{r})}} & \text{if }  r > j, 
             \IEEEyesnumber\\[.5ex]%\IEEEyessubnumber* \\[.5ex]
       P_{0}^{\otimes {(b_{r}+1-a_{r})}} 
                      & \text{if }  r = j\text{ and } i = a_{j}, \\[.5ex]
       P_{0}^{\otimes {(b_{r}+1-i)}} \otimes S^{*}\otimes I^{\otimes(i-a_{r}-1)} 
                      & \text{if }  r = j\text{ and } i \geq a_{j}+1, \\[.5ex]
         P_{0}^{\otimes {(b_{r}+1-a_{r})}}
                     & \text{if } r<j.  
        \end{cases}\label{eq:irred-q0-1cb}
    \end{IEEEeqnarray}
\ecrlre
%%---------------------------------------------------------

%%---------------------------------------------------------
\bppsn %[\texttt{pr:irred-q0-2}]
	%%---------------------------
	\label{pr:irred-q0-2}
	%%---------------------------
 For $1\leq j\leq k$ and $a_{j}+1\leq i\leq b_{j}+1$, let
 \begin{IEEEeqnarray*}{rCl}
    R_{j,i} &=& I^{\otimes(\sum_{r=j+1}^{k}\ell_{r})}\otimes 
\Bigl(P_{0}^{\otimes(b_{j}+1-i)} \otimes S^{*}\otimes I^{\otimes(i-a_{j}-1)}\Bigr) 
      \otimes       I^{\otimes(\sum_{r=1}^{j-1}\ell_{r})},\\
    S_{j,i} &=& I^{\otimes(\sum_{r=j+1}^{k}\ell_{r})}\otimes 
      \Bigl(I^{\otimes(b_{j}+1-i)} \otimes S^{*}\otimes I^{\otimes(i-a_{j}-1)}\Bigr)
       \otimes       I^{\otimes(\sum_{r=1}^{j-1}\ell_{r})}.
 \end{IEEEeqnarray*}
 Then $R_{j,i}$ and  $S_{j,i}$ belong to 
 $\psi_{\lambda,\omega}^{}(C(SU_{0}(n+1)))^{\prime\prime}$
 for $1\leq j\leq k$ and $a_{j}+1\leq i\leq b_{j}+1$.
\eppsn
%%---------------------------------------------------------
\begin{proof}
%%------------------------------------------------------ 
This is easy to see for $j=1$, as 
one has  $R_{1,i}^{}=V_{1,i}^{}$  for $a_{1}+1\leq i\leq b_{1}+1$, 
$S_{1,1+b_{1}}=V_{1,1+b_{1}}$ and
\begin{IEEEeqnarray*}{rCll}
 S_{1,i}^{}    
      &=& 
      \sum_{k_{i},\ldots k_{b_{j}}\in\mathbb{N}}^{} 
        S_{1,1+b_{1}}^{k_{b_{1}}}\ldots S_{1,i+1}^{k_{i}}
        V_{1,i}^{}
        (S_{1,1+b_{1}}^{k_{b_{1}}})^{*}\ldots (S_{1,i+1}^{k_{i}})^{*}
           & \qquad\text{if }a_{1}\leq i \leq b_{1}.
\end{IEEEeqnarray*}
Fix a $j<k$. Assume that 
$R_{r,i},S_{r,i}\in\psi_{\lambda,\omega}^{}(C(SU_{0}(n+1)))^{\prime\prime}$
for $1\leq r\leq j$. 
Take $i$ such that $a_{j+1}+1\leq i\leq b_{j+1}+1$. 
% Let us denote the operators in Corollary~\ref{cr:irred-q0-1c} by $E_{j,i}^{}$.
Let $E_{j,i}^{}$ be as in Corollary~\ref{cr:irred-q0-1c}.
Then
\begin{IEEEeqnarray*}{rCl}
 R_{j+1,i}^{} &=&
   \sum_{r=1}^{j-1}\sum_{k_{r,a_{r}},\ldots k_{r,b_{r}}\in\mathbb{N}}^{}
     S_{r,1+b_{r}}^{k_{r,b_{r}}}\ldots S_{r,1+a_{r}}^{k_{r,a_{r}}}
        E_{j+1,i}^{}
        (S_{r,1+a_{r}}^{k_{r,a_{r}}})^{*}\ldots (S_{r,1+b_{r}}^{k_{r,b_{r}}})^{*}.
\end{IEEEeqnarray*}
Since $E_{j+1,i}^{}$ belongs to the image $\psi_{\lambda,\omega}^{}(C(SU_{0}(n+1)))$,
the operator $R_{j+1,i}^{}$ belongs to 
$\psi_{\lambda,\omega}^{}(C(SU_{0}(n+1)))^{\prime\prime}$.
Now note that
\begin{IEEEeqnarray*}{rCl}
 S_{j+1,i}^{}    
      &=& 
      \sum_{k_{i},\ldots k_{b_{j}}\in\mathbb{N}}^{} 
        R_{j+1,1+b_{j}}^{k_{b_{j}}}\ldots R_{j+1,i+1}^{k_{i}}
        R_{j+1,i}^{}
        (R_{j+1,i+1}^{k_{i}})^{*}\ldots (R_{j+1,1+b_{j}}^{k_{b_{j}}})^{*},
\end{IEEEeqnarray*}
Therefore it follows that  $S_{j+1,i}$ also belongs to the strong operator closure
of the image $\psi_{\lambda,\omega}^{}(C(SU_{0}(n+1)))$.
\end{proof}
%%---------------------------------------------------------

%%---------------------------------------------------------

%%---------------------------------------------------------
\bthm %[\texttt{th:irred-q0-3}]
	%%---------------------------
	\label{th:irred-q0-3}
	%%---------------------------
  Let $\lambda\in (S^{1})^{n}$ and 
  $\omega=s_{[a_{k},b_{k}]}s_{[a_{k-1},b_{k-1}]}\ldots  s_{[a_{1},b_{1}]}$
  be a reduced word in $\mathfrak{S}_{n+1}$.
  Then the representation $\psi_{\lambda,\omega}$ is an irreducible representation.
\ethm
%%---------------------------------------------------------
\begin{proof}
This is an immediate consequence of the previous proposition.
\end{proof}
%%---------------------------------------------------------

%%---------------------------------------------------------
\bthm %[\texttt{th:irred-q0-4}]
	%%---------------------------
	\label{th:irred-q0-4}
	%%---------------------------
Let 
$\lambda,\lambda^{\prime}\in (S^{1})^{n}$
and let
\[
\omega=s_{[a_{k},b_{k}]}s_{[a_{k-1},b_{k-1}]}\ldots  s_{[a_{1},b_{1}]},\qquad
\omega_{}^{\prime}=s_{[a_{k'}^{\prime},b_{k'}^{\prime}]}
  s_{[a_{k'-1}^{\prime},b_{k'-1}^{\prime}]}\ldots  
      s_{[a_{1}^{\prime},b_{1}^{\prime}]}
\]
be two reduced words in  $\mathfrak{S}_{n+1}$.
If $(\lambda,\omega)\neq (\lambda',\omega')$, then the representations 
$\psi_{\lambda,\omega}$ and $\psi_{\lambda',\omega'}$ are inequivalent.
\ethm
%%---------------------------------------------------------

%%---------------------------------------------------------
\begin{proof}
Let $1\leq j\leq n+1$.
Define 
$m_{k+1}(j,\omega) = j$ and for $1\leq i\leq k$, define $m_{i}(j,\omega)$ 
recursively  as follows:
\begin{IEEEeqnarray*}{rCl}
   m_{i}(j,\omega)  &=& 
    \begin{cases}
     m_{i+1}(j,\omega)+1 & \text{if } a_{i}\leq m_{i+1}(j,\omega)\leq b_{i},\\
     m_{i+1}(j,\omega) & \text{if } m_{i+1}(j,\omega) \geq b_{i}+1
                             \text{ or } m_{i+1}(j,\omega) < a_{i}.
          \end{cases}
\end{IEEEeqnarray*}
In the diagram for $\psi_{\lambda,\omega}$, one starts at $j$
and proceeds horizontally and goes one step up whenever one encounters 
a diagonal going upward.  The numbers $m_{r}(j,\omega)$ simply keep 
track of the horizontal line one is on when one completes the  
$s_{[a_{r},b_{r}]}$-section.
Thus $m_{1}(j,\omega)$ is the maximum $k$ such 
that $\psi_{\lambda,\omega}(z_{j,k}^{})\neq 0$.

Let us first assume that $\omega=\omega'$ and $\lambda\neq\lambda'$.
Let 
\[
j:=\max\{1\leq i\leq n: \lambda_{i}^{}\neq\lambda_{i}^{\prime}\}.
\]
If $b_{i+1}< j<b_{i}$ for some $i$, then 
\begin{IEEEeqnarray*}{rCl}
\text{spec}(\psi_{\lambda,\omega}(z_{j+1,m_{1}^{}(j+1,\omega)})) &=&
   \{0,\overline{\lambda_{j}^{}}\lambda_{j+1}^{} \}, \\
\text{spec}(\psi_{\lambda',\omega}(z_{j+1,m_{1}^{}(j+1,\omega)})) &=&
  \{0,\overline{\lambda_{j}^{\prime}}\lambda_{j+1}^{\prime} \},
\end{IEEEeqnarray*}
i.e.\ the spectrums of the same element are different for the two representations.
Therefore $\psi_{\lambda,\omega}$ and $\psi_{\lambda',\omega'}$ are inequivalent.
If $j=b_{i}$ for some $i$, then write $m=n'(i,a_{i})$.
Now using Corollary~\ref{cr:irred-q0-1c} and keeping track of the 
effect of $\lambda$ and $\lambda'$, we get that if $m>-\infty$, then
\begin{IEEEeqnarray*}{rCl}
\text{spec}\bigl(\psi_{\lambda^{\prime},\omega}
            (E_{i,a_{i}}^{}(\omega))\bigr) 
   &=& \Bigl\{0,\overline{\lambda}_{b_{m}}^{\prime}\overline{\lambda}_{j}^{\prime}
                        \lambda_{b_{m}+1}^{\prime}\lambda_{j+1}^{\prime}\Bigr\}\\
   &=& \Bigl\{0,\overline{\lambda}_{b_{m}}^{}\overline{\lambda}_{j}^{\prime}
                        \lambda_{b_{m}+1}^{}\lambda_{j+1}^{}\Bigr\},\\
\text{spec}\bigl(\psi_{\lambda^{},\omega^{}}(E_{i,a_{i}}^{}(\omega))\bigr)
   &=&  \Bigl\{0,\overline{\lambda}_{b_{m}}^{}\overline{\lambda}_{j}^{}
                        \lambda_{b_{m}+1}^{}\lambda_{j+1}^{}\Bigr\},
\end{IEEEeqnarray*}
and if $m=-\infty$, then
\begin{IEEEeqnarray*}{rCl}
\text{spec}\bigl(\psi_{\lambda^{\prime},\omega}
            (E_{i,a_{i}}^{}(\omega))\bigr) 
   &=& \Bigl\{0,\overline{\lambda}_{j}^{\prime}\lambda_{j+1}^{\prime}\Bigr\}\\
   &=& \Bigl\{0,\overline{\lambda}_{j}^{\prime}\lambda_{j+1}^{}\Bigr\},\\
\text{spec}\bigl(\psi_{\lambda^{},\omega^{}}(E_{i,a_{i}}^{}(\omega))\bigr)
   &=&  \Bigl\{0,\overline{\lambda}_{j}^{}\lambda_{j+1}^{}\Bigr\},
\end{IEEEeqnarray*}
Thus $\psi_{\lambda,\omega}$ and $\psi_{\lambda',\omega'}$ are 
inequivalent  in this case also.

Let us next assume that $\omega\neq \omega^{\prime}$. We will split 
this into different cases and
in each case, produce an element of $C(SU_{0}(n+1))$ whose image 
under one representation will be zero while the image under the 
other is nonzero, thus proving that the two representations are not equivalent.\\
\textbf{Case I}:\\ There exists  $j \leq \min\{k,k'\}$ such that
$b_{i}^{}=b_{i}^{\prime}$, $a_{i}^{}=a_{i}^{\prime}$ for $1\leq i\leq j-1$ 
and $b_{j}^{}\neq b_{j}^{\prime}$.\\
Assume without loss in generality that $b_{j}^{}> b_{j}^{\prime}$. Then 
$m_{1}(b_{j},\omega)>m_{1}(b_{j},\omega')$. Therefore
\begin{IEEEeqnarray*}{rCl}
  \psi_{\lambda^{\prime},\omega^{\prime}}(z_{b_{j},m_{1}(b_{j},\omega)}) &=& 0,\\
  \psi_{\lambda,\omega}(z_{b_{j},m_{1}(b_{j},\omega)}) &\neq & 0.
\end{IEEEeqnarray*}
\textbf{Case II}:\\ There exists  $j \leq \min\{k,k'\}$ such that
$b_{i}^{}=b_{i}^{\prime}$, $a_{i}^{}=a_{i}^{\prime}$ for $1\leq i\leq j-1$, 
$b_{j}^{}= b_{j}^{\prime}$ and $a_{j}^{}\neq a_{j}^{\prime}$.\\
Assume without loss in generality that $a_{j}^{}> a_{j}^{\prime}$.
Once again using Corollary~\ref{cr:irred-q0-1c}, we observe that 
modulo multiplication by complex numbers of modulus one,
$\psi_{\lambda,\omega}(E_{j,a_{j}}(\omega))$ is a projection, while 
$\psi_{\lambda',\omega'}(E_{j,a_{j}}(\omega))$ is a tensor product of shifts 
and projections. Therefore $\psi_{\lambda,\omega}$ and $\psi_{\lambda',\omega'}$
are inequivalent.\\
\textbf{Case III}:\\ $b_{i}^{}=b_{i}^{\prime}$, $a_{i}^{}=a_{i}^{\prime}$ 
for $1\leq i\leq k'$ and $k'< k$.\\
In this case $m_{1}(b_{k},\omega)> m_{1}(b_{k},\omega')$. Therefore
\begin{IEEEeqnarray*}{rCl}
  \psi_{\lambda^{\prime},\omega^{\prime}}(z_{b_{k}, m_{1}(b_{k},\omega)}^{}) &=& 0,\\
  \psi_{\lambda,\omega}(z_{b_{k}, m_{1}(b_{k},\omega)}^{}) &\neq & 0.
\end{IEEEeqnarray*}
\textbf{Case IV}:\\ $b_{i}^{}=b_{i}^{\prime}$, $a_{i}^{}=a_{i}^{\prime}$ 
for $1\leq i\leq k$ and $k< k'$.\\
This case is identical to Case III with $k$ and $k'$ interchanged.
\end{proof}
%%---------------------------------------------------------

%%--------------------------------------------------------------------
% \input{an0-irr-sec4.tex}
%%--------------------------------------------------------------------
%%--------------------------------------------------------------------
 \section{Projections in $C(SU_{0}(n+1))$}
%%--------------------------------------------------------------------
In this section we will restrict ourselves to the {$*$}-subalgebra
$\mathcal{O}(SU_{0}(n+1))$ of $C(SU_{0}(n+1))$ generated by the $z_{i,j}^{}$'s 
and derive certain relations that will in particular show that the 
generators $z_{i,j}^{}$ are all partial isometries.
%%--------------------------------------------------------------------
\blmma %[\texttt{lm:proj-1}]
	%%---------------------------
	\label{lm:proj-1}
	%%---------------------------
    Let us denote by $p_{i,j}^{}$ the element $z_{i,j}^{*}z_{i,j}^{}$
    and by $q_{i,j}^{}$ the element $z_{i,j}^{}z_{i,j}^{*}$.
  For $1\leq i \leq j\leq n+1$, one has the following:
  \begin{IEEEeqnarray*}{rClCl}
      \sum_{k=1}^{i}q_{k,j}^{} &=& \sum_{k=j}^{n+1}p_{i,k}^{} &&
	%%---------------------------
        \yesnumber\label{eq:proj-1}\\
	%%---------------------------
      \sum_{k=1}^{j}q_{k,j}^{} &=& \sum_{k=j}^{n+1}p_{j,k}^{} &=& 1.
	%%---------------------------
        \yesnumber\label{eq:proj-2}
	%%---------------------------
  \end{IEEEeqnarray*}
\elmma
%%----------------------------------------------------------------
\begin{proof}
We will prove that
  \begin{IEEEeqnarray*}{rCl}
      \sum_{k=1}^{i}q_{k,j}^{}
        &=&
        \begin{cases}
        (z_{1,1}^{}z_{2,2}^{}\cdots z_{i-1,i-1}^{})z_{i,j}^{}
        (z_{i+1,i}^{}z_{i+2,i+1}^{}\dots z_{j,j-1}^{})
        (z_{j+1,j+1}^{}\dots z_{n+1,n+1}^{}) 
        & \\%\text{if }i<j,\\ 
        \hspace{200pt}\text{if }i<j,&\\
        1 \hspace{195pt}\text{if }i=j&
        \end{cases}%\\
	%%---------------------------
        \yesnumber\label{eq:proj-3}\\
	%%---------------------------
      \sum_{k=j}^{n+1}p_{i,k}^{}
        &=&
        \begin{cases}
        (z_{1,1}^{}z_{2,2}^{}\cdots z_{i-1,i-1}^{})
        z_{i,j}^{}
        (z_{i+1,i}^{}z_{i+2,i+1}^{}\dots z_{j,j-1}^{})
        (z_{j+1,j+1}^{}\dots z_{n+1,n+1}^{}) 
        &\\%\text{if }i<j,\\
        \hspace{200pt}\text{if }i<j,&\\
        1 \hspace{195pt}\text{if }i=j.&
        \end{cases}%\\
	%%---------------------------
       \yesnumber\label{eq:proj-4}
	%%---------------------------
  \end{IEEEeqnarray*}
Note that for $i=1$, one has
\begin{IEEEeqnarray*}{rCl}
        q_{1,j}^{} 
    &=& z_{1,j}^{}z_{1,j}^{*}\\
    &=& z_{1,j}^{}(z_{2,1}^{}z_{3,2}^{}\cdots z_{j,j-1}^{})
            (z_{j+1,j+1}^{}z_{j+2,j+2}^{}\cdots z_{n+1,n+1}^{}).
\end{IEEEeqnarray*}
Thus the first equality in (\ref{eq:proj-3}) holds for $i=1$. Next assume that 
it holds for $i-1$ and $1\leq i<j\leq n+1$. Then using (\ref{eq:an0-defrel-9}), 
we get
\begin{IEEEeqnarray*}{rCl}
    \IEEEeqnarraymulticol{3}{l}{
      q_{i,j}^{} + \sum_{k=1}^{i-1}q_{k,j}^{}}\\
        &=& z_{i,j}^{}
        (z_{1,1}^{}z_{2,2}^{}\cdots z_{i-2,i-2}^{}z_{i-1,i-1}^{})
        (z_{i+1,i}^{}z_{i+2,i+1}^{}\cdots z_{j,j-1}^{})
        (z_{j+1,j+1}^{}\cdots z_{n+1,n+1}^{})\\
    && %{}\quad 
       {} + (z_{1,1}^{}z_{2,2}^{}\cdots z_{i-2,i-2}^{})
        z_{i-1,j}^{}z_{i,i-1}^{}
        (z_{i+1,i}^{}z_{i+2,i+1}^{}\cdots z_{j,j-1}^{})
        (z_{j+1,j+1}^{}\cdots z_{n+1,n+1}^{})\\
    &=& (z_{1,1}^{}z_{2,2}^{}\cdots z_{i-2,i-2}^{})
        z_{i,j}^{}z_{i-1,i-1}^{}
        (z_{i+1,i}^{}z_{i+2,i+1}^{}\cdots z_{j,j-1}^{})
        (z_{j+1,j+1}^{}\cdots z_{n+1,n+1}^{})\\
    && %{}\quad 
        {} + (z_{1,1}^{}z_{2,2}^{}\cdots z_{i-2,i-2}^{})
        z_{i-1,j}^{}z_{i,i-1}^{}
        (z_{i+1,i}^{}z_{i+2,i+1}^{}\cdots z_{j,j-1}^{})
        (z_{j+1,j+1}^{}z_{j+2,j+2}^{}\cdots z_{n+1,n+1}^{})\\
    &=& (z_{1,1}^{}\cdots z_{i-2,i-2}^{}z_{i-1,i-1}^{})
        z_{i,j}^{}
        (z_{i+1,i}^{}z_{i+2,i+1}^{}\cdots z_{j,j-1}^{})
        (z_{j+1,j+1}^{}\cdots z_{n+1,n+1}^{}).
\end{IEEEeqnarray*}
Hence the first equality in (\ref{eq:proj-3}) follows. Using this, one now gets
\begin{IEEEeqnarray*}{rCl}
 \sum_{k=1}^{j}q_{k,j}^{}
    &=& \sum_{k=1}^{j-1}q_{k,j}^{}+q_{j,j}^{}\\
    &=& (z_{1,1}^{}z_{2,2}^{}\cdots z_{j-2,j-2}^{})z_{j-1,j}^{}z_{j,j-1}^{}
          (z_{j+1,j+1}^{}z_{j+2,j+2}^{}\cdots z_{n+1,n+1}^{})\\
    && {}\quad {}+ z_{j,j}^{} (z_{1,1}^{}z_{2,2}^{}
          \cdots z_{j-2,j-2}^{}z_{j-1,j-1}^{})
         (z_{j+1,j+1}^{}z_{j+2,j+2}^{}\cdots z_{n+1,n+1}^{})\\
    &=& (z_{1,1}^{}z_{2,2}^{}\cdots z_{j-2,j-2}^{})z_{j-1,j}^{}z_{j,j-1}^{}
           (z_{j+1,j+1}^{}z_{j+2,j+2}^{}\cdots z_{n+1,n+1}^{})\\
    && {}\quad {}+ (z_{1,1}^{}z_{2,2}^{}\cdots z_{j-2,j-2}^{})
          z_{j,j}^{}z_{j-1,j-1}^{} (z_{j+1,j+1}^{}z_{j+2,j+2}^{}
            \cdots z_{n+1,n+1}^{})\\
    &=& z_{1,1}^{}z_{2,2}^{}\cdots z_{j-2,j-2}^{}z_{j-1,j-1}^{}
            z_{j,j}^{}z_{j+1,j+1}^{}z_{j+2,j+2}^{}\cdots z_{n+1,n+1}^{}\\
    &=& 1
\end{IEEEeqnarray*}
Thus we have (\ref{eq:proj-3}).
Proof of (\ref{eq:proj-4}) is similar. First observe that for $j=n+1$, one has
\begin{IEEEeqnarray*}{rCl}
        p_{i,n+1}^{}
        &=& z_{i,n+1}^{*}z_{i,n+1}^{}\\
        &=& (z_{1,1}^{}z_{2,2}^{}\cdots z_{i-1,i-1}^{})
             (z_{i+1,i}^{}z_{i+2,i+1}^{}\dots z_{n+1,n}^{})
              z_{i,n+1}^{}\\
        &=& (z_{1,1}^{}z_{2,2}^{}\cdots z_{i-1,i-1}^{})
               z_{i,n+1}^{}
             (z_{i+1,i}^{}z_{i+2,i+1}^{}\dots z_{n+1,n}^{}).
\end{IEEEeqnarray*}
Thus the first equality in (\ref{eq:proj-4}) holds for $j=n+1$.
Now, assuming that it holds for $j+1$, we have
\begin{IEEEeqnarray*}{rCl}
    \IEEEeqnarraymulticol{3}{l}{
        p_{i,j}^{} + \sum_{k=j+1}^{n+1}p_{i,k}^{}}\\
 &=&
    (z_{1,1}^{}z_{2,2}^{}\cdots z_{i-2,i-2}^{}z_{i-1,i-1}^{})
    (z_{i+1,i}^{}z_{i+2,i+1}^{}\cdots z_{j,j-1}^{})
    (z_{j+1,j+1}^{}z_{j+2,j+2}^{}\cdots z_{n+1,n+1}^{})
    z_{i,j}^{}\\
 && {}\quad + (z_{1,1}^{}z_{2,2}^{}\cdots z_{i-1,i-1}^{})
    z_{i,j+1}^{}
    (z_{i+1,i}^{}z_{i+2,i+1}^{}\cdots z_{j,j-1}^{})z_{j+1,j}^{}
    (z_{j+2,j+2}^{}\cdots z_{n+1,n+1}^{})\\
    &=&
    (z_{1,1}^{}z_{2,2}^{}\cdots z_{i-2,i-2}^{}z_{i-1,i-1}^{})
    (z_{i+1,i}^{}z_{i+2,i+1}^{}\cdots z_{j,j-1}^{})
    z_{j+1,j+1}^{}z_{i,j}^{}
    (z_{j+2,j+2}^{}\cdots z_{n+1,n+1}^{})\\
    &&{}\quad  + (z_{1,1}^{}z_{2,2}^{}\cdots z_{i-1,i-1}^{})
      (z_{i+1,i}^{}z_{i+2,i+1}^{}\cdots z_{j,j-1}^{})
        z_{i,j+1}^{}z_{j+1,j}^{}
    (z_{j+2,j+2}^{}\cdots z_{n+1,n+1}^{})    \\
    &=& (z_{1,1}^{}z_{2,2}^{}\cdots z_{i-1,i-1}^{})
        (z_{i+1,i}^{}z_{i+2,i+1}^{}\cdots z_{j,j-1}^{})
       z_{i,j}^{}
    (z_{j+1,j+1}^{}z_{j+2,j+2}^{}\cdots z_{n+1,n+1}^{})\\
    &=& (z_{1,1}^{}z_{2,2}^{}\cdots z_{i-1,i-1}^{})
    z_{i,j}^{}
    (z_{i+1,i}^{}z_{i+2,i+1}^{}\cdots z_{j,j-1}^{})
    (z_{j+1,j+1}^{}z_{j+2,j+2}^{}\cdots z_{n+1,n+1}^{}).
\end{IEEEeqnarray*}
This proves the first equality in (\ref{eq:proj-4}).
Using  this we now get
\begin{IEEEeqnarray*}{rCl}
    \sum_{k=i}^{n+1}p_{i,k}^{}
    &=&
    \sum_{k=i+1}^{n+1}p_{i,k}^{}+p_{i,i}^{}\\
    &=&
    (z_{1,1}^{}z_{2,2}^{}\cdots z_{i-1,i-1}^{})
    z_{i,i+1}^{}z_{i+1,i}^{}
    (z_{i+2,i+2}^{}z_{i+3,i+3}^{}\cdots z_{n+1,n+1}^{})\\
    &&
     {}\quad {} + (z_{1,1}^{}z_{2,2}^{}\cdots z_{i-1,i-1}^{})
      (z_{i+1,i+1}^{}z_{i+2,i+2}^{}z_{i+3,i+3}^{}\cdots z_{n+1,n+1}^{})
     z_{i,i}^{}\\
    &=&
    (z_{1,1}^{}z_{2,2}^{}\cdots z_{i-1,i-1}^{})
    z_{i,i+1}^{}z_{i+1,i}^{}
    (z_{i+2,i+2}^{}z_{i+3,i+3}^{}\cdots z_{n+1,n+1}^{})\\
    &&
    {}\quad {}+ (z_{1,1}^{}z_{2,2}^{}\cdots z_{i-1,i-1}^{})
    z_{i+1,i+1}^{}z_{i,i}^{}
    (z_{i+2,i+2}^{}z_{i+3,i+3}^{}\cdots z_{n+1,n+1}^{})\\
    &=&
    z_{1,1}^{}z_{2,2}^{}\cdots z_{i-1,i-1}^{}z_{i,i}^{}z_{i+1,i+1}^{}
       z_{i+2,i+2}^{}z_{i+3,i+3}^{}\cdots z_{n+1,n+1}^{}\\
    &=& 1,
\end{IEEEeqnarray*}
which completes the proof.
\end{proof}
%%----------------------------------------------------------------

%%----------------------------------------------------------------
\bppsn %[\texttt{pr:proj-2}]
    %%------------------------------
    \label{pr:proj-2}
    %%------------------------------
  The elements $p_{i,j}^{}$ and $q_{i,j}^{}$ are projections 
  for $1\leq i\leq j\leq n$.  
\eppsn
%%----------------------------------------------------------------
\begin{proof}
We will show that for any $1\leq j\leq n$, if 
  $\left\{\sum_{k=1}^{i} q_{k,j}^{}: 1\leq i\leq j\right\}$
  is a family of projections then 
  $\left\{\sum_{k=1}^{i} q_{k,j+1}^{}: 1\leq i\leq j+1 \right\}$ 
  is also a family of projections. 
  Then by an application of equation~(\ref{eq:proj-2}) in Lemma~\ref{lm:proj-1}, 
  we will successively conclude for $i=1, 2, \ldots, n+1$ that the 
  sums $\sum_{k=1}^{i} q_{k,j}^{}$ are all projections for $i\leq j\leq n+1$. 
  This will imply that all $q_{i,j}^{}$'s (and hence all $p_{i,j}^{}$'s) 
  for $1\leq i\leq j\leq n+1$ are projections.

Now fix $1\leq j\leq n$ and assume that 
$\left\{\sum_{k=1}^{i} q_{k,j}^{}: 1\leq i\leq j\right\}$
  is a family of projections.
Note that by this assumption, it follows that 
$q_{i,j}^{}=\sum_{k=1}^{i} q_{k,j}^{}-\sum_{k=1}^{i-1} q_{k,j}^{}$,
being positive and a difference of two projections, is also a projection.
Hence $p_{i,j}^{}$ is also a projection. 
Now from Lemma~\ref{lm:proj-1}, for $i\leq j$, we get
\[
\sum_{k=1}^{i}q_{k,j+1}^{} 
  = \sum_{k=j+1}^{n+1}p_{i,k}^{}
    = \sum_{k=j}^{n+1}p_{i,k}^{} - p_{i,j}^{}
      = \sum_{k=1}^{i}q_{k,j}^{} - p_{i,j}^{}.
\]
Thus  the element $\sum_{k=1}^{i}q_{k,j+1}^{}$ 
is positive  and is a difference of two projections. 
Therefore it is a projection. For $i=j+1$, 
 $\sum_{k=1}^{j+1}q_{k,j+1}^{}=1$ and hence is a projection.
\end{proof}

%%----------------------------------------------------------------
\blmma %[\texttt{lm:proj-3}]
    %%------------------------------
    \label{lm:proj-3}
    %%------------------------------
For $1\leq j\leq i \leq n+1$, one has
\begin{IEEEeqnarray*}{rClCl}
    \sum_{k=1}^{j}q_{i,k}^{} &=& \sum_{k=i}^{n+1}p_{k,j}^{},&&
	%%---------------------------
       \yesnumber\label{eq:proj-5}\\
	%%---------------------------
    \sum_{k=1}^{i}q_{i,k}^{} &=& \sum_{k=i}^{n+1}p_{k,i}^{} &=& 1.
	%%---------------------------
       \yesnumber\label{eq:proj-6}
	%%---------------------------
\end{IEEEeqnarray*}
\elmma
%%----------------------------------------------------------------
\begin{proof}
We will show that
\begin{IEEEeqnarray*}{rCl}
    \sum_{k=1}^{j}q_{i,k}^{}
    &=&
    \begin{cases}
      (z_{1,1}^{}z_{2,2}^{}\cdots z_{j-1,j-1}^{})
        z_{i,j}^{}
        (z_{j,j+1}^{}z_{j+1,j+2}^{}\dots z_{i-1,i}^{})
           (z_{i+1,i+1}^{}\dots z_{n+1,n+1}^{}), &\\
      \hspace{200pt} \text{if }j<i,&\\
      1 \hspace{195pt} \text{if }j=i,&
      \end{cases}
	%%---------------------------
       \yesnumber\label{eq:proj-7}\\
	%%---------------------------
      \sum_{k=i}^{n+1}p_{k,j}^{}
        &=&
    \begin{cases}
        (z_{1,1}^{}z_{2,2}^{}\cdots z_{j-1,j-1}^{})
        z_{i,j}^{}
        (z_{j,j+1}^{}z_{j+1,j+2}^{}\dots z_{i-1,i}^{})
        (z_{i+1,i+1}^{}\dots z_{n+1,n+1}^{}), &\\
      \hspace{200pt} \text{if }i>j,&\\
        1  \hspace{195pt} \text{if }i=j.&
      \end{cases}
	%%---------------------------
       \yesnumber\label{eq:proj-8}
	%%---------------------------
\end{IEEEeqnarray*}

For $j=1$, the first equality in (\ref{eq:proj-7}) holds as one has
\begin{IEEEeqnarray*}{rCl}
        q_{i,1}^{} 
    &=& z_{i,1}^{}z_{i,1}^{*}\\
    &=& z_{i,1}^{}(z_{1,2}^{}z_{2,3}^{}\cdots z_{i-1,i}^{})
            (z_{i+1,i+1}^{}z_{i+2,i+2}^{}\cdots z_{n+1,n+1}^{}).
\end{IEEEeqnarray*}
Next, assuming that it holds for $j-1$, we have
\begin{IEEEeqnarray*}{rCl}
    \IEEEeqnarraymulticol{3}{l}{
      q_{i,j}^{}+\sum_{k=1}^{j-1}q_{i,k}^{}}\\
    &=& z_{i,j}^{}
    (z_{1,1}^{}z_{2,2}^{}\cdots z_{j-1,j-1}^{})
    (z_{j,j+1}^{}z_{j+1,j+2}^{}\cdots z_{i-1,i}^{})
    (z_{i+1,i+1}^{}z_{i+2,i+2}^{}\cdots z_{n+1,n+1}^{})\\
    && %{}\quad 
    {}+ (z_{1,1}^{}z_{2,2}^{}\cdots z_{j-2,j-2}^{})
    z_{i,j-1}^{}z_{j-1,j}^{}
    (z_{j,j+1}^{}z_{j+1,j+2}^{}\cdots z_{i-1,i}^{})
    (z_{i+1,i+1}^{}z_{i+2,i+2}^{}\cdots z_{n+1,n+1}^{})\\
    &=& (z_{1,1}^{}z_{2,2}^{}\cdots z_{j-2,j-2}^{})
    z_{i,j}^{}z_{j-1,j-1}^{}
    (z_{j,j+1}^{}z_{j+1,j+2}^{}\cdots z_{i-1,i}^{})
    (z_{i+1,i+1}^{}z_{i+2,i+2}^{}\cdots z_{n+1,n+1}^{})\\
    && %{}\quad 
    {}+ (z_{1,1}^{}z_{2,2}^{}\cdots z_{j-2,j-2}^{})
    z_{i,j-1}^{}z_{j-1,j}^{}
    (z_{j,j+1}^{}z_{j+1,j+2}^{}\cdots z_{i-1,i}^{})
    (z_{i+1,i+1}^{}z_{i+2,i+2}^{}\cdots z_{n+1,n+1}^{})\\
    &=& (z_{1,1}^{}z_{2,2}^{}\cdots z_{j-1,j-1}^{})
    z_{i,j}^{}
    (z_{j,j+1}^{}z_{j+1,j+2}^{}\cdots z_{i-1,i}^{})
    (z_{i+1,i+1}^{}z_{i+2,i+2}^{}\cdots z_{n+1,n+1}^{}).
\end{IEEEeqnarray*}
Thus the equality holds.   For $i=j$, one has
\begin{IEEEeqnarray*}{rCl}
    \sum_{k=1}^{i}q_{i,k}^{}
    &=& \sum_{k=1}^{i-1}q_{i,k}^{}+q_{i,i}^{}\\
    &=& (z_{1,1}^{}z_{2,2}^{}\cdots z_{i-2,i-2}^{})
        z_{i,i-1}^{}z_{i-1,i}^{}
        (z_{i+1,i+1}^{}z_{i+2,i+2}^{}\cdots z_{n+1,n+1}^{})\\
    && {}\quad {} + z_{i,i}^{}
     (z_{1,1}^{}z_{2,2}^{}\cdots z_{i-2,i-2}^{})z_{i-1,i-1}^{}
     (z_{i+1,i+1}^{}z_{i+2,i+2}^{}\cdots z_{n+1,n+1}^{})\\
    &=& (z_{1,1}^{}z_{2,2}^{}\cdots z_{i-2,i-2}^{})z_{i,i-1}^{}z_{i-1,i}^{}
      (z_{i+1,i+1}^{}z_{i+2,i+2}^{}\cdots z_{n+1,n+1}^{})\\
    && {}\quad {}+ (z_{1,1}^{}z_{2,2}^{}\cdots z_{i-2,i-2}^{})
        z_{i,i}^{}z_{i-1,i-1}^{}
        (z_{i+1,i+1}^{}z_{i+2,i+2}^{}\cdots z_{n+1,n+1}^{})\\
    &=& z_{1,1}^{}z_{2,2}^{}\cdots z_{i-2,i-2}^{}z_{i-1,i-1}^{}z_{i,i}^{}
          z_{i+1,i+1}^{}z_{i+2,i+2}^{}\cdots z_{n+1,n+1}^{}\\
    &=& 1.
\end{IEEEeqnarray*}
Thus we have (\ref{eq:proj-7}).

Next, note that
\begin{IEEEeqnarray*}{rCl}
        p_{n+1,j}^{}
        &=& z_{n+1,j}^{*}z_{n+1,j}^{}\\
        &=& (z_{1,1}^{}z_{2,2}^{}\cdots z_{j-1,j-1}^{})
             (z_{j,j+1}^{}z_{j+1,j+2}^{}\dots z_{n,n+1}^{})
              z_{n+1,j}^{}\\
        &=& (z_{1,1}^{}z_{2,2}^{}\cdots z_{j-1,j-1}^{})
               z_{n+1,j}^{}
             (z_{j,j+1}^{}z_{j+1,j+2}^{}\dots z_{n,n+1}^{}).
\end{IEEEeqnarray*}
Thus the first equality in (\ref{eq:proj-8}) holds for $i=n+1$. 
Assume it holds for $i+1$.
Then
\begin{IEEEeqnarray*}{rCl}
    \IEEEeqnarraymulticol{3}{l}{
    p_{i,j}^{}+\sum_{t=i+1}^{n+1}p_{t,j}^{} }\\
&=& (z_{1,1}^{}z_{2,2}^{}\cdots z_{j-1,j-1}^{})
    (z_{j,j+1}^{}z_{j+1,j+2}^{}\cdots z_{i-1,i}^{})
    (z_{i+1,i+1}^{}z_{i+2,i+2}^{}\cdots z_{n+1,n+1}^{})
    z_{i,j}^{}\\
&& {}\quad {}+ (z_{1,1}^{}z_{2,2}^{}\cdots z_{j-1,j-1}^{})
    z_{i+1,j}^{}
    (z_{j,j+1}^{}z_{j+1,j+2}^{}\cdots z_{i-1,i}^{})
    z_{i,i+1}^{}
    (z_{i+2,i+2}^{}\cdots z_{n+1,n+1}^{}) \\
&=& (z_{1,1}^{}z_{2,2}^{}\cdots z_{j-2,j-2}^{}z_{j-1,j-1}^{})
    (z_{j,j+1}^{}z_{j+1,j+2}^{}\cdots z_{i-1,i}^{})
    z_{i+1,i+1}^{}z_{i,j}^{}
    (z_{i+2,i+2}^{}\cdots z_{n+1,n+1}^{})\\
 && {}\quad {}+ (z_{1,1}^{}z_{2,2}^{}\cdots z_{j-1,j-1}^{})
    (z_{j,j+1}^{}z_{j+1,j+2}^{}\cdots z_{i-1,i}^{})
    z_{i,i+1}^{}z_{i+1,j}^{}
    (z_{i+2,i+2}^{}\cdots z_{n+1,n+1}^{})     \\
&=& (z_{1,1}^{}z_{2,2}^{}\cdots z_{j-1,j-1}^{})
    (z_{j,j+1}^{}z_{j+1,j+2}^{}\cdots z_{i-1,i}^{})
    z_{i,j}^{}
    (z_{i+1,i+1}^{}z_{i+2,i+2}^{}\cdots z_{n+1,n+1}^{})\\
&=& (z_{1,1}^{}z_{2,2}^{}\cdots z_{j-1,j-1}^{})
    z_{i,j}^{}
    (z_{j,j+1}^{}z_{j+1,j+2}^{}\cdots z_{i-1,i}^{})
    (z_{i+1,i+1}^{}z_{i+2,i+2}^{}\cdots z_{n+1,n+1}^{}).
\end{IEEEeqnarray*}
Thus we have the first equality in (\ref{eq:proj-8}). For $i=j$,
\begin{IEEEeqnarray*}{rCl}
    \sum_{k=j}^{n+1}p_{k,j}^{}
&=&\sum_{k=j+1}^{n+1}p_{k,j}^{}+p_{j,j}^{}\\
&=& (z_{1,1}^{}z_{2,2}^{}\cdots z_{j-1,j-1}^{})
   z_{j+1,j}^{}z_{j,j+1}^{}
   (z_{j+2,j+2}^{}z_{j+3,j+3}^{}\cdots z_{n+1,n+1}^{})\\
&& {}\quad {}+ (z_{1,1}^{}z_{2,2}^{}\cdots z_{j-1,j-1}^{})
       (z_{j+1,j+1}^{}z_{j+2,j+2}^{}\cdots z_{n+1,n+1}^{})
       z_{j,j}^{}\\
&=& (z_{1,1}^{}z_{2,2}^{}\cdots z_{j-1,j-1}^{})
     z_{j+1,j}^{}z_{j,j+1}^{}
     (z_{j+2,j+2}^{}z_{j+3,j+3}^{}\cdots z_{n+1,n+1}^{})\\
&& {}\quad {}+ (z_{1,1}^{}z_{2,2}^{}\cdots z_{j-1,j-1}^{})
      z_{j+1,j+1}^{}z_{j,j}^{}
      (z_{j+2,j+2}^{}z_{j+3,j+3}^{}\cdots z_{n+1,n+1}^{})\\
    &=& z_{1,1}^{}z_{2,2}^{}\cdots z_{j-1,j-1}^{}z_{j,j}^{}z_{j+1,j+1}^{}z_{j+2,j+2}^{}\cdots z_{n+1,n+1}^{}\\
    &=&1.
\end{IEEEeqnarray*}
Thus the proof is complete.
\end{proof}

%%----------------------------------------------------------------
\bppsn %[\texttt{pr:proj-4}]
	%%---------------------------
	\label{pr:proj-4}
	%%---------------------------
  The elements $p_{i,j}^{}$ and $q_{i,j}^{}$ are projections 
  for $1\leq j\leq i\leq n$. 
\eppsn
%%----------------------------------------------------------------
\begin{proof}
By using Lemma~\ref{lm:proj-3}, one can prove exactly as in the 
proof of Proposition~\ref{pr:proj-2} that if 
$\left\{\sum_{k=1}^{j} q_{i,k}^{}: 1\leq j\leq i\right\}$ 
is a family of projections then 
$\left\{\sum_{k=1}^{j} q_{i+1,k}^{}: 1\leq j\leq i+1\right\}$ 
is also a family of projections.
This leads to the required conclusion.
\end{proof}
%%--------------------------------------------------------------------

%%----------------------------------------------------------------
\bcrlre %[\texttt{cr:proj-5}]
	%%---------------------------
	\label{cr:proj-5}
	%%---------------------------
The elements $z_{i,j}^{}$ are partial isometries for all $1\leq i,j\leq n+1$.
\ecrlre
%%----------------------------------------------------------------
%%----------------------------------------------------------------

%%---------------------------------------------------------
\blmma %[\texttt{lm:proj-6}]
	%%---------------------------
	\label{lm:proj-6}
	%%---------------------------
Let $i<k$ and $j<\ell$. Assume that either $i\geq\ell$ or $k\leq j$. Then
\begin{IEEEeqnarray}{rClCl}
  z_{i,\ell}^{}z_{k,j}^{*} &=& 0 &=& z_{i,\ell}^{*}z_{k,j}^{}.
	%%---------------------------
	\label{eq:proj-10}
	%%---------------------------
\end{IEEEeqnarray}
\elmma
%%---------------------------------------------------------
\begin{proof}
Assume that $i\geq\ell$. Then using Lemma~\ref{lm:proj-3}, we have
\begin{IEEEeqnarray*}{rCl}
  \sum_{k=i+1}^{n+1}
    \bigl(z_{i,\ell}^{}z_{k,j}^{*}\bigr)\bigl(z_{i,\ell}^{}z_{k,j}^{*}\bigr)^{*} 
  &=& 
      z_{i,\ell}^{}\left(\sum_{k=i+1}^{n+1} p_{k,j}^{}\right)z_{i,\ell}^{*}\\
  &=& z_{i,\ell}^{} \Bigl(\sum_{s=1}^{j}q_{i+1,s}^{}\Bigr)z_{i,\ell}^{*}\\
  &=& \sum_{s=1}^{j}z_{i,\ell}^{}z_{i+1,s}^{}z_{i+1,s}^{*}z_{i,\ell}^{*}\\
  &=& 0,
\end{IEEEeqnarray*}
since $z_{i,\ell}^{}z_{i+1,s}^{}=0$ for all $1\leq s\leq \ell-1$.
Therefore $z_{i,\ell}^{}z_{k,j}^{*}=0$ for all $j<\ell$ and $i<k$.
Similarly if $k\leq j$ then using Lemma~\ref{lm:proj-1}, we get
\begin{IEEEeqnarray*}{rCl}
  \sum_{\ell=j+1}^{n+1}
    \bigl(z_{i,\ell}^{}z_{k,j}^{*}\bigr)^{*} \bigl(z_{i,\ell}^{}z_{k,j}^{*}\bigr)
  &=& z_{k,j}^{}
             \sum_{\ell=j+1}^{n+1}\bigl(p_{i,\ell}^{}\bigr)z_{k,j}^{*}\\
  &=& z_{k,j}^{}
        \Bigl(\sum_{s=1}^{i}q_{s,j+1}^{}\Bigr)z_{k,j}^{*}\\
  &=& 0.
\end{IEEEeqnarray*}
Therefore as before, $z_{i,\ell}^{}z_{k,j}^{*}=0$.
The second equality can be proved in a similar way using Lemmas \ref{lm:proj-1}
and \ref{lm:proj-3}.
\end{proof}
%%---------------------------------------------------------

%%--------------------------------------------------
\blmma %[\texttt{lm:proj-7}]
	%%---------------------------
	\label{lm:proj-7}
	%%---------------------------
Let $i,j\leq s$ and $k\geq s+1$. Then
\[
q_{i,j}^{}z_{n+1,k}^{}= z_{n+1,k}^{}q_{i,j}^{},\qquad\qquad 
   z_{i,j}^{}p_{n+1,k}^{} = p_{n+1,k}^{}z_{i,j}^{} .
\]
\elmma
%%--------------------------------------------------
\begin{proof}
If $\max\{i,j\}+1<k$, then one has $z_{i,j}^{}z_{n+1,k}^{}= z_{n+1,k}^{}z_{i,j}^{}$
and hence the equality follows. 
Assume $\max\{i,j\}+1=k$ (i.e.\ $k=s+1$ and $\max\{i,j\}=s$). Then
\[
z_{i,j}^{}z_{n+1,k}^{}= z_{n+1,k}^{}z_{i,j}^{}+z_{i,k}^{}z_{n+1,j}^{}.
\]
It follows from Lemmas \ref{lm:proj-1} and \ref{lm:proj-3} that if $i\geq j$, 
then $p_{n+1,j}^{}p_{i,j}^{}=0$ and if $i<j$, then
$p_{i,k}^{}p_{i,j}^{}=0$. Therefore in either case we get from 
the above equation that
\begin{IEEEeqnarray*}{rCl}
q_{i,j}^{}z_{n+1,k}^{} &=& z_{i,j}^{}z_{n+1,k}^{}z_{i,j}^{*}\\
  &=& z_{n+1,k}^{}z_{i,j}^{}z_{i,j}^{*}+z_{i,k}^{}z_{n+1,j}^{}z_{i,j}^{*}\\
  &=& z_{n+1,k}^{}q_{i,j}^{}.
\end{IEEEeqnarray*}
Similarly, We have 
\begin{IEEEeqnarray*}{rCl}
z_{i,j}^{}p_{n+1,k}^{} &=& z_{i,j}^{}z_{n+1,k}^{*}z_{n+1,k}^{}\\
  &=& z_{n+1,k}^{*}z_{i,j}^{}z_{n+1,k}^{}\\
  &=& p_{n+1,k}^{}z_{i,j}^{}+ z_{n+1,k}^{*}z_{i,k}^{}z_{n+1,j}^{}\\
  &=& p_{n+1,k}^{}z_{i,j}^{}+ z_{n+1,k}^{*}z_{n+1,j}^{}z_{i,k}^{}\\
  &=& p_{n+1,k}^{}z_{i,j}^{}.
\end{IEEEeqnarray*}
\end{proof}
%--------------------------------------------------

%%--------------------------------------------------
\blmma %[\texttt{lm:proj-8}]
    %%---------------------------
    \label{lm:proj-8}
    %%---------------------------
Let $1\leq i\leq j\leq n+1$. Then
\[
z_{i,j+1}^{}z_{n+1,j+1}^{*} = z_{n+1,j}^{*}z_{i,j}^{}.
\]
\elmma
%%--------------------------------------------------
\begin{proof}
We have
\begin{IEEEeqnarray*}{rCl}
 \IEEEeqnarraymulticol{3}{l}{
    \left(z_{i,j+1}^{}z_{n+1,j+1}^{*} - z_{n+1,j}^{*}z_{i,j}^{}\right)
      \left(z_{i,j+1}^{}z_{n+1,j+1}^{*} - z_{n+1,j}^{*}z_{i,j}^{}\right)^{*}\qquad}\\
   \qquad &=& z_{i,j+1}^{}p_{n+1,j+1}^{}z_{i,j+1}^{*}
              - z_{i,j+1}^{}z_{n+1,j+1}^{*}z_{i,j}^{*}z_{n+1,j}^{}\\
     &&\quad  {}    + z_{n+1,j}^{*}q_{i,j}^{}z_{n+1,j}^{}  
                  - z_{n+1,j}^{*}z_{i,j}^{}z_{n+1,j+1}^{}z_{i,j+1}^{*}\\
    &=& z_{i,j+1}^{}(p_{n+1,j}^{}+q_{n+1,j+1}^{})z_{i,j+1}^{*}
     - z_{i,j+1}^{}(z_{n+1,j+1}^{}z_{i,j}^{}+z_{n+1,j}^{}z_{i,j+1}^{})^{*}
              z_{n+1,j}^{}\\
    && \quad    
            {}+ z_{n+1,j}^{*}q_{i,j}^{}z_{n+1,j}^{}
               - z_{n+1,j}^{*}(z_{n+1,j+1}^{}z_{i,j}^{}+z_{n+1,j}^{}z_{i,j+1}^{})
               z_{i,j+1}^{*}\\
    &=& z_{i,j+1}^{}p_{n+1,j}^{}z_{i,j+1}^{*}
     - z_{i,j+1}^{}z_{i,j+1}^{*}z_{n+1,j}^{*}z_{n+1,j}^{}\\
     && \quad  
          {}  + z_{n+1,j}^{*}q_{i,j}^{}z_{n+1,j}^{}
     - z_{n+1,j}^{*}z_{n+1,j}^{}z_{i,j+1}^{}z_{i,j+1}^{*}\\
    &=& p_{n+1,j}^{}q_{i,j+1}^{} + z_{n+1,j}^{*}q_{i,j}^{}z_{n+1,j}^{}
           - p_{n+1,j}^{}q_{i,j+1}^{} 
               - z_{n+1,j}^{*}q_{i,j+1}^{}z_{n+1,j}^{}\\
    &=& z_{n+1,j}^{*}(q_{i,j}^{}-q_{i,j+1}^{})z_{n+1,j}^{}.
\end{IEEEeqnarray*}
Summing over $i$ we get
\begin{IEEEeqnarray*}{rCl}
 \IEEEeqnarraymulticol{3}{l}{
    \sum_{i=1}^{j}\left(z_{i,j+1}^{}z_{n+1,j+1}^{*} - z_{n+1,j}^{*}z_{i,j}^{}\right)
      \left(z_{i,j+1}^{}z_{n+1,j+1}^{*} - z_{n+1,j}^{*}z_{i,j}^{}\right)^{*}}\\
    \qquad\qquad\qquad\qquad&=&  z_{n+1,j}^{*}
      \left(\sum_{i=1}^{j}q_{i,j}^{}-\sum_{i=1}^{j}q_{i,j+1}^{}\right)z_{n+1,j}^{}\\
    &=& z_{n+1,j}^{*}\left(I-\sum_{k=j+1}^{n+1}p_{j,k}^{}\right)z_{n+1,j}^{}\\
    &=& z_{n+1,j}^{*}p_{j,j}^{}z_{n+1,j}^{}\\
    &=& 0.
\end{IEEEeqnarray*}
Therefore we have the result.
\end{proof}
%%--------------------------------------------------

%%--------------------------------------------------------------------
\section{Operators in $\pi(C(SU_{0}(n+1)))^{\prime\prime}$}
%%--------------------------------------------------------------------
In this section, we will study the actions of the generating elements
$z_{i,j}^{}$ on a Hilbert space and construct certain operators using the 
strong operator topology available there.  These will turn out to be 
very useful tool in analyzing  the irreducible representations.
Let us start with  a representation $\pi$ of $C(SU_{0}(n+1))$ 
on a Hilbert space $\mathcal{H}$. For all $i$ and $j$, we will denote the 
operators $\pi(z_{i,j}^{})$, $\pi(p_{i,j}^{})$ and $\pi(q_{i,j}^{})$ by 
$Z_{i,j}^{}$, $P_{i,j}^{}$ and $Q_{i,j}^{}$ respectively.

%%------------------------------------------------------
\bppsn %[\texttt{pr:rep-rel-1}]
       %%-----------------------------
       \label{pr:rep-rel-1}
       %%-----------------------------
Let $1\leq s\leq m\leq n$. Assume $Z_{m+1, s}^{}=0$ where $s<m+1$.  
Then $Z_{i,j}^{}=0$ for $i\geq m+1$ and $j \leq s$.
\eppsn
%%------------------------------------------------------
\begin{proof}
Let us assume $s>1$ and take $i=m+1$, $j=s-1$.
From relation~(\ref{eq:an0-defrel-4}), we have $Z_{m+1,s-1}^{}Z_{s,s}^{}=0$, so that
\[
Z_{m+1,s-1}^{}Q_{s,s}^{}=0.
\]
For $1\leq k\leq s-1$, the commutation relation~(\ref{eq:an0-defrel-5}) gives us
\[
Z_{k,s-1}^{}Z_{m+1,s}^{}-Z_{m+1,s}^{}Z_{k,s-1}^{}=Z_{m+1,s-1}^{}Z_{k,s}^{}.
\]
Since $Z_{m+1,s}^{}=0$, we get $Z_{m+1,s-1}^{}Z_{k,s}^{}=0$ which implies
\[
Z_{m+1,s-1}^{}Q_{k,s}^{}=0,\qquad 1\leq k\leq s-1.
\]
Therefore
\begin{equation}
Z_{m+1,s-1}^{}=Z_{m+1,s-1}^{}\left(\sum_{k=1}^{s}Q_{k,s}^{}\right)=0.
       %%-----------------------------
       \label{eq:rep-rel-1}
       %%-----------------------------
\end{equation}
%%------------------------------------------------------
Next, let us assume $m<n$ and take $i=m+2$, $j=s$.
Using relation~(\ref{eq:an0-defrel-4}) we obtain \linebreak
$Z_{m+1,m+1}^{}Z_{m+2,s}^{}=0$, so that
\[
P_{m+1,m+1}^{}Z_{m+2,s}^{}=0.
\]
For $m+2\leq k\leq n+1$, we get from~(\ref{eq:an0-defrel-5})
\[
Z_{m+1,s}^{}Z_{m+2,k}^{}-Z_{m+2,k}^{}Z_{m+1,s}^{}=Z_{m+1,k}^{}Z_{m+2,s}^{}.
\]
Since $Z_{m+1,s}^{}=0$, we get $Z_{m+1,k}^{}Z_{m+2,s}^{}=0$ which implies
\[
P_{m+1,k}^{}Z_{m+2,s}^{}=0,\qquad m+2\leq k\leq n+1.
\]
Therefore
%%------------------------------------------------------
\begin{equation}
Z_{m+2,s}^{}=\left(\sum_{k=m+1}^{n+1}P_{m+1,k}^{}\right)Z_{m+2,s}^{}=0.
       %%-----------------------------
       \label{eq:rep-rel-2}
       %%-----------------------------
\end{equation}
%%------------------------------------------------------
By repeated application of (\ref{eq:rep-rel-1}) and (\ref{eq:rep-rel-2}), 
we get $Z_{i,j}^{}=0$ for $i\geq m+1$ and $j \leq s$.
\end{proof}

%%------------------------------------------------------
\blmma %[\texttt{lm:rep-rel-2}]
       %%-----------------------------
       \label{lm:rep-rel-2}
       %%-----------------------------
Assume $i\neq m$ and $j\neq k$. If 
$Z_{i,j}^{}Q_{m,k}^{}=Q_{m,k}^{}Z_{i,j}^{}$ and 
$Z_{m,k}^{}Q_{i,j}^{}=Q_{i,j}^{}Z_{m,k}^{}$, then one has
$Z_{i,j}^{}Z_{m,k}^{}=Z_{m,k}^{}Z_{i,j}^{}$.
\elmma
%%------------------------------------------------------
\begin{proof}
Computing $\left(Z_{i,j}^{}Z_{m,k}^{}-Z_{m,k}^{}Z_{i,j}^{}\right)
   \left(Z_{i,j}^{}Z_{m,k}^{}-Z_{m,k}^{}Z_{i,j}^{}\right)^{*}$ 
   using the given relations we get the required equality.
\end{proof}
%%------------------------------------------------------

%%------------------------------------------------------
\blmma %[\texttt{lm:rep-rel-3}]
       %%-----------------------------
       \label{lm:rep-rel-3}
       %%-----------------------------
Let $F$ be a subset of $\{1,2,\ldots,n+1\}$. Let $j\not\in F$ and $i\neq m$.
Assume $Z_{m,k}^{}Q_{i,j}^{}=Q_{i,j}^{}Z_{m,k}^{}$ for all $k\in F$ and 
$Z_{i,j}^{}\left(\sum_{k\in F}Q_{m,k}^{}\right)=
  \left(\sum_{k\in F}Q_{m,k}^{}\right)Z_{i,j}^{}$. 
Then
\[
Z_{i,j}^{}Z_{m,k}^{}=Z_{m,k}^{}Z_{i,j}^{} \qquad\text{for all }k\in F.
\] 
\elmma
%%------------------------------------------------------
\begin{proof}
Note that
\begin{IEEEeqnarray*}{rCl}
   \IEEEeqnarraymulticol{3}{l}{
   \sum_{k\in F}\left(Z_{i,j}^{}Z_{m,k}^{}-Z_{m,k}^{}Z_{i,j}^{}\right)
      \left(Z_{i,j}^{}Z_{m,k}^{}-Z_{m,k}^{}Z_{i,j}^{}\right)^{*} }\\
   &=&
     Z_{i,j}^{}\left(\sum_{k\in F}Q_{m,k}^{}\right)Z_{i,j}^{*} 
       + \sum_{k\in F}\left(Z_{m,k}^{}Q_{i,j}^{}Z_{m,k}^{*}\right) \\
    &&\qquad \qquad{} 
         - Z_{i,j}^{}\sum_{k\in F}\left(Z_{m,k}^{}Z_{i,j}^{*}Z_{m,k}^{*}\right)
         - \sum_{k\in F}\left(Z_{m,k}^{}Z_{i,j}^{}Z_{m,k}^{*}\right)Z_{i,j}^{*}\\
   &=&
     Q_{i,j}^{} \left(\sum_{k\in F}Q_{m,k}^{}\right)
       + Q_{i,j}^{}\sum_{k\in F}\left(Z_{m,k}^{}Z_{m,k}^{*}\right)   \\
    &&\qquad \qquad{} 
         - Z_{i,j}^{}Z_{i,j}^{*}\sum_{k\in F}\left(Z_{m,k}^{}Z_{m,k}^{*}\right)
         - Z_{i,j}^{}Z_{i,j}^{*}\sum_{k\in F}\left(Z_{m,k}^{}Z_{m,k}^{*}\right)\\
    &=& 0.
\end{IEEEeqnarray*}
Therefore we have the result.
\end{proof}
%%------------------------------------------------------

Given a representation $\pi$ of $C(SU_{0}(n+1))$, we will 
denote by $r(\pi)\equiv r$ the following
\begin{equation}
r(\pi)\equiv r := \min\{i: 1\leq i\leq n+1, Z_{n+1,i}^{}\neq 0\}.  
       %%-----------------------------
       \label{eq:rep-rel-3}
       %%----------------------------- 
\end{equation}

%%------------------------------------------------------
\bppsn %[\texttt{pr:rep-rel-4}]
       %%-----------------------------
       \label{pr:rep-rel-4}
       %%-----------------------------
Assume $2\leq r\leq n$.  Then one has
\begin{IEEEeqnarray}{rCl}
Z_{i,j}^{}Z_{n+1,k}^{} &=& Z_{n+1,k}^{}Z_{i,j}^{} 
 \quad\text{for }1\leq i\leq n,\; 1\leq j\leq r-1\text{ and }r\leq k\leq n+1.
       %%-----------------------------
       \label{eq:rep-rel-4}
       %%-----------------------------
\end{IEEEeqnarray}
\eppsn
%%------------------------------------------------------
\begin{proof}
Observe that for $1\leq i\leq r-1$, $1\leq j\leq r-1$ and $r\leq k\leq n+1$,
one has $\max\{i,j\}+1\leq r \leq \min\{n+1,k\}$. 
If $\max\{i,j\}+1<\min\{n+1,k\}$, then 
$Z_{i,j}^{}Z_{n+1,k}^{}=Z_{n+1,k}^{}Z_{i,j}^{}$. Since by 
Proposition~\ref{pr:rep-rel-1}, we have $Z_{n+1,j}^{}=0$ for 
$1\leq j\leq r-1$, if $\max\{i,j\}+1=\min\{n+1,k\}$, then one has
%%------------------------------------------------------
\begin{IEEEeqnarray*}{rClCl}
  Z_{i,j}^{}Z_{n+1,k}^{} &=& Z_{n+1,k}^{}Z_{i,j}^{}+Z_{n+1,j}^{}Z_{i,k}^{}
     &=& Z_{n+1,k}^{}Z_{i,j}^{}.
\end{IEEEeqnarray*}
%%------------------------------------------------------
Let us next assume that 
$Z_{i,j}^{}Z_{n+1,k}^{}=Z_{n+1,k}^{}Z_{i,j}^{}$
for all $1\leq j\leq r-1$, $r\leq k\leq n+1$ and for all $1\leq i\leq t<n$.
Then 
$P_{i,j}^{}Z_{n+1,k}^{}=Z_{n+1,k}^{}P_{i,j}^{}$ for all such $j$, $k$ and $i$.
Therefore
\[
Z_{n+1,k}^{}=Z_{n+1,k}^{}\left(\sum_{i=j}^{n}P_{i,j}^{}\right)=
\left(\sum_{i=j}^{n}P_{i,j}^{}\right)Z_{n+1,k}^{}.
\]
Since $P_{i,j}^{}Z_{n+1,k}^{}=Z_{n+1,k}^{}P_{i,j}^{}$ for $1\leq i\leq t$ and
$\sum_{i=t+1}^{n}P_{i,j}^{}=\sum_{m=1}^{j}Q_{t+1,m}^{}$, it follows that
\[
Z_{n+1,k}^{}\left(\sum_{m=1}^{j}Q_{t+1,m}^{}\right)=
\left(\sum_{m=1}^{j}Q_{t+1,m}^{}\right)Z_{n+1,k}^{}.
\]
Therefore
\[
Z_{n+1,k}^{}Q_{t+1,j}^{}=Q_{t+1,j}^{}Z_{n+1,k}^{}
\]
for $1\leq j\leq r-1$.
On the other hand, we also have
\[
Z_{t+1,j}^{}=Z_{t+1,j}^{}
\left(\sum_{k=r}^{n+1}Q_{n+1,k}^{}\right)=
\left(\sum_{k=r}^{n+1}Q_{n+1,k}^{}\right)Z_{t+1,j}^{}.
\]
Therefore by Lemma~\ref{lm:rep-rel-3}, we have
\[
Z_{t+1,j}^{}Z_{n+1,k}^{}=Z_{n+1,k}^{}Z_{t+1,j}^{}
\]
for $r\leq k\leq n+1$ and $1\leq j\leq r-1$.

By repeated application, we now get the result.
\end{proof}
%%------------------------------------------------------

%%--------------------------------------------------------------------
% \subsection{Existence of the SOT sums}  %STARTS
%%--------------------------------------------------------------------
We next look at certain infinite sums of elements from $\pi(C(SU_{0}(n+1)))$
that converge in the strong operator topology (SOT) so that they belong
to the SOT-closure $\pi(C(SU_{0}(n+1)))^{\prime\prime}$ of $\pi(C(SU_{0}(n+1)))$. 
%%------------------------------------------------------
\bppsn %[\texttt{pr:rep-rel-5}]
       %%-----------------------------
       \label{pr:rep-rel-5}
       %%-----------------------------
    Let $r$ be as in (\ref{eq:rep-rel-3}) and let $r\leq k\leq n+1$. 
    If $T_{k}^{}$ is a partial isometry with $T_{k}^{*}T_{k}^{}\leq P_{n+1,k}^{}$ 
    and $T_{k}^{}T_{k}^{*}\leq P_{n+1,k}^{}$, then the operator 
\[
T_{k+1}:= \sum_{m\in\mathbb{N}}Z_{n+1,k+1}^{m}T_{k}^{}(Z_{n+1,k+1}^{m})^{*}
\]
    exists as a strong operator convergent sum and is a partial isometry with 
    $T_{k+1}^{*}T_{k+1}^{}\leq P_{n+1,k+1}^{}$ and 
    $T_{k+1}^{}T_{k+1}^{*}\leq P_{n+1,k+1}^{}$.
    % initial and range spaces contained in $P_{n+1,k+1}^{}\mathcal{H}$.
\eppsn
%%------------------------------------------------------
\begin{proof}
Note that the operators  
\[
 R_{k,m}:= Z_{n+1,k+1}^{m}P_{n+1,k}^{}(Z_{n+1,k+1}^{m})^{*},\qquad m\in\mathbb{N}
\] 
 form a family of orthogonal projections with 
 $\sum_{m\in\mathbb{N}}R_{k,m}^{}\leq P_{n+1,k+1}^{}$. Observe also that
\[
    P_{n+1,k+1}^{}-\sum_{j=0}^{m}R_{k,j}^{}
    =
      Z_{n+1,k+1}^{m+1}(Z_{n+1,k+1}^{m+1})^{*}.
\]
and the sum $\sum_{j\in\mathbb{N}}R_{k,j}^{}$ converges in SOT. Hence
\begin{IEEEeqnarray}{rCl}
 P_{n+1,k+1}^{}-\sum_{j\in\mathbb{N}}R_{k,j}^{}
     &=& \text{SOT-}\lim_{m\to\infty}
       Z_{n+1,k+1}^{m}(Z_{n+1,k+1}^{m})^{*}.
       %%-----------------------------
       \label{eq:rep-rel-5a}
       %%-----------------------------
\end{IEEEeqnarray}
Next, write
\[
T_{k,m}^{}=Z_{n+1,k+1}^{m}T_{k}^{}(Z_{n+1,k+1}^{m})^{*},\qquad m\in\mathbb{N}.
\]
Then one has $T_{k,m}^{} = R_{k,m}^{}T_{k,m}^{}R_{k,m}^{}$ for all $m\in\mathbb{N}$.
It is easy to check that $T_{k,m}^{*}T_{k,m'}^{}=0=T_{k,m}^{}T_{k,m'}^{*}$ 
for $m\neq m'$ and one has
\begin{IEEEeqnarray*}{rCl}
T_{k,m}^{*}T_{k,m}^{} &\leq& P_{n+1,k+1},\\
T_{k,m}^{}T_{k,m}^{*} &\leq& P_{n+1,k+1}.
\end{IEEEeqnarray*}
Therefore the sum $T:=\sum_{m\in\mathbb{N}}T_{k,m}^{}$ converges in SOT and it is a 
partial isometry with both initial and range spaces contained in 
$P_{n+1,k+1}^{}\mathcal{H}$.
\end{proof}
%%------------------------------------------------------
Using the above proposition, one can now define recursively the 
 operators $W_{i,j}$ for $1\leq j\leq n+1$ and
 $j\leq i\leq n+1$ and the operators $U_{i}$ for $r+1\leq i\leq n+1$ as follows:
%%------------------------------------------------------  
  \begin{IEEEeqnarray}{rCl}
 W_{i,j}^{}&:=& 
 \begin{cases}
 Z_{n+1,j}^{},&\text{if }i=j,\\
 % W_{i+1,j}^{} &:=& 
 \sum_{k\in\mathbb{N}}
        Z_{n+1,i}^{k}W_{i-1,j}^{}(Z_{n+1,i}^{k})^{*},&
             \text{if } j+1\leq i\leq n+1.
  \end{cases}
       %%-----------------------------
       \label{eq:rep-rel-5}\\
       %%-----------------------------
  &&\nonumber\\
 U_{i}^{} &:=& 
 \begin{cases}
 Z_{n+1,r}^{},&\text{if }i=r+1,\\
 % U_{i+1}^{} &:=& 
 \sum_{k\in\mathbb{N}}
      Z_{n+1,i}^{k}U_{i-1}^{}(Z_{n+1,i}^{k})^{*},&\text{if } r+2\leq i\leq n+1.
  \end{cases}
       %%-----------------------------
       \label{eq:rep-rel-6}
       %%-----------------------------
  \end{IEEEeqnarray}
%%------------------------------------------------------
Thus we have for $i\geq r+2$,
\begin{IEEEeqnarray*}{rCl}
W_{i,r+1}^{}&=&
  \sum_{\substack{k_{j}\in\mathbb{N}\\(r+2\leq j\leq i)}}^{}
        Z_{n+1,i}^{k_{i}}Z_{n+1,i-1}^{k_{i-1}}\cdots Z_{n+1,r+2}^{k_{r+2}}
        Z_{n+1,r+1}^{}
        (Z_{n+1,r+2}^{k_{r+2}})^{*}\cdots (Z_{n+1,i-1}^{k_{i-1}})^{*}
           (Z_{n+1,i}^{k_{i}})^{*},\\
       %%-----------------------------
       \yesnumber\label{eq:rep-rel-7}\\
       %%-----------------------------
U_{i} &=&
  \sum_{\substack{k_{j}\in\mathbb{N}\\(r+2\leq j\leq i)}}^{}
        Z_{n+1,i}^{k_{i}}Z_{n+1,i-1}^{k_{i-1}}\cdots Z_{n+1,r+2}^{k_{r+2}}
        Z_{n+1,r}^{}
        (Z_{n+1,r+2}^{k_{r+2}})^{*}\cdots (Z_{n+1,i-1}^{k_{i-1}})^{*}
           (Z_{n+1,i}^{k_{i}})^{*},\\
       %%-----------------------------
       \yesnumber\label{eq:rep-rel-8}
       %%-----------------------------
\end{IEEEeqnarray*}
and for $i\geq r+1$,
\begin{IEEEeqnarray*}{rCl}
W_{i,r}^{}&=&
  \sum_{\substack{k_{j}\in\mathbb{N}\\(r+1\leq j\leq i)}}^{}
        Z_{n+1,i}^{k_{i}}Z_{n+1,i-1}^{k_{i-1}}\cdots Z_{n+1,r+1}^{k_{r+1}}
        Z_{n+1,r}^{}
        (Z_{n+1,r+1}^{k_{r+1}})^{*}\cdots (Z_{n+1,i-1}^{k_{i-1}})^{*}
           (Z_{n+1,i}^{k_{i}})^{*}.\\
       %%-----------------------------
       \yesnumber\label{eq:rep-rel-7a}
       %%-----------------------------
\end{IEEEeqnarray*}
We will use these three operators repeatedly in the rest of the paper.
We start with a few observations on these operators.

%%--------------------------------------------------------------------
% \subsection{Existence of the SOT sums}  ENDS
%%--------------------------------------------------------------------

%%------------------------------------------------------
\bppsn %[\texttt{pr:rep-rel-7}]
       %%-----------------------------
       \label{pr:rep-rel-7}
       %%-----------------------------
Let $r$ be as in (\ref{eq:rep-rel-3}). Then
\begin{enumerate}
    \item
    $U_{i}^{}$ is a normal partial isometry for $r+2\leq i\leq n+1$,
    \item
    $W_{i,r}^{}$ is a normal partial isometry for all $r\leq i\leq n+1$,
	\item
	$W_{i,r+1}^{}$ is a partial isometry with source projection 
	and range projection contained in $P_{n+1,i}^{}$ for $r+2\leq i\leq n+1$,
\end{enumerate}
\eppsn
%%------------------------------------------------------
\begin{proof}
%%------------------------------------------------------
By Proposition~\ref{pr:rep-rel-5}, the operators $U_{i}^{}$, 
$W_{i,r}^{}$ and $W_{i,r+1}^{}$ are all partial isometries.
Normality of $U_{i}^{}$ and $W_{i,r}^{}$ follow from the fact that
$P_{n+1,r}=Q_{n+1,r}$. It is easy to see from (\ref{eq:rep-rel-7})
and (\ref{eq:rep-rel-8}) that
%%------------------------------------------------------
\begin{IEEEeqnarray*}{rCl}
 W_{i,r+1}^{*}W_{i,r+1}^{} &\leq & P_{n+1,i}^{},
       %%-----------------------------
       \yesnumber\label{eq:rep-rel-12}\\
       %%-----------------------------
 W_{i,r+1}^{*}W_{i,r+1}^{}-W_{i,r+1}^{}W_{i,r+1}^{*} &=& U_{i}^{*}U_{i}^{}.
       %%-----------------------------
       \yesnumber\label{eq:rep-rel-13}
       %%-----------------------------
\end{IEEEeqnarray*}
%%------------------------------------------------------
Therefore part~3 follows. 
%%------------------------------------------------------
\end{proof}
%%------------------------------------------------------

%%------------------------------------------------------
\bppsn %[\texttt{pr:rep-rel-6}]
       %%-----------------------------
       \label{pr:rep-rel-6}
       %%-----------------------------
Let $r$ be as in (\ref{eq:rep-rel-3}).
Then we have the following relations:
%%------------------------------------------------------  
\begin{IEEEeqnarray}{rCll}
 Z_{i,r}^{}W_{n+1,r+1}^{}-W_{n+1,r+1}^{}Z_{i,r}^{} &=&U_{n+1}^{}Z_{i,r+1}^{},\quad &i\leq r,
       %%-----------------------------
       \label{eq:rep-rel-10}\\
       %%-----------------------------
 Z_{i,r}^{}W_{n+1,r+1}^{*} - W_{n+1,r+1}^{*}Z_{i,r}^{} &=& 0, & i\leq r.
       %%-----------------------------
       \label{eq:rep-rel-11}
       %%-----------------------------
\end{IEEEeqnarray}
%%------------------------------------------------------
\eppsn
%%------------------------------------------------------
\begin{proof}
Since $i\leq r$, the operator $Z_{i,r}^{}$ commutes with 
$Z_{n+1,j}^{}$ for $r+2\leq j\leq n+1$. 
Using this and the relation (\ref{eq:an0-defrel-5}), 
we get (\ref{eq:rep-rel-10}). Equation (\ref{eq:rep-rel-11})
follows by using (\ref{eq:an0-defrel-5}) and (\ref{eq:an0-defrel-8}).
\end{proof}

%%------------------------------------------------------
\bthm %[\texttt{th:rep-rel-9}]
       %%-----------------------------
       \label{th:rep-rel-9}
       %%-----------------------------
Let $r$ be as in (\ref{eq:rep-rel-3}).
	Then the operator $W_{n+1,r}^{}$ commutes with all the $Z_{i,j}^{}$'s.
\ethm
%%------------------------------------------------------
\begin{proof}
%%------------------------------------------------------
% \textred{$1\leq j\leq r-1$:} \\
%%------------------------------------------------------
For $1\leq j\leq r-1$, the commutaions follow from Proposition~\ref{pr:rep-rel-4}.
%%------------------------------------------------------
% \textred{$1\leq÷ i\leq j=r$:} \\
%%------------------------------------------------------
Assume next that $1\leq i\leq j=r$. Then
\begin{IEEEeqnarray*}{rCl}
 \IEEEeqnarraymulticol{3}{l}{
 Z_{i,r}^{}W_{r+1,r}^{} }\\
   &=& Z_{i,r}^{}
       \sum_{k\in\mathbb{N}}
         Z_{n+1,r+1}^{k}Z_{n+1,r}^{}(Z_{n+1,r+1}^{k})^{*}\\
   &=& \sum_{k\in\mathbb{N}}Z_{n+1,r+1}^{k}Z_{i,r}^{}
                   Z_{n+1,r}^{}(Z_{n+1,r+1}^{k})^{*}
        + \sum_{k\in\mathbb{N}}Z_{n+1,r+1}^{k}
                Z_{n+1,r}^{}Z_{i,r+1}^{}Z_{n+1,r}^{}(Z_{n+1,r+1}^{k+1})^{*}\\
   &=& \sum_{k\in\mathbb{N}}Z_{n+1,r+1}^{k}
                Z_{n+1,r}^{2}Z_{i,r+1}^{}(Z_{n+1,r+1}^{k+1})^{*},
\end{IEEEeqnarray*}
and
\begin{IEEEeqnarray*}{rCl}
 W_{r+1,r}^{} Z_{i,r}^{}
   &=& \sum_{k\in\mathbb{N}}Z_{n+1,r+1}^{k}
                   Z_{n+1,r}^{}(Z_{n+1,r+1}^{*})^{k}Z_{i,r}^{}\\
   &=& \sum_{k\in\mathbb{N}}Z_{n+1,r+1}^{k}
                   Z_{n+1,r}^{}Z_{i,r}^{}(Z_{n+1,r+1}^{k})^{*}.
\end{IEEEeqnarray*}
Thus we have
\begin{IEEEeqnarray*}{rCl}
 \IEEEeqnarraymulticol{3}{l}{Z_{i,r}^{}W_{r+1,r}^{} - W_{r+1,r}^{} Z_{i,r}^{}}\\
   \qquad &=& \sum_{k\in\mathbb{N}}
         Z_{n+1,r+1}^{k}Z_{n+1,r}^{}
      \left( Z_{n+1,r}^{}Z_{i,r+1}^{}Z_{n+1,r+1}^{*} - Z_{i,r}^{} \right)
       (Z_{n+1,r+1}^{k})^{*}.
\end{IEEEeqnarray*}
Now
\begin{IEEEeqnarray*}{rCl}
 \IEEEeqnarraymulticol{3}{l}{
    \left(Z_{n+1,r}^{2}Z_{i,r+1}^{}Z_{n+1,r+1}^{*} - Z_{n+1,r}^{}Z_{i,r}^{}\right)
\left(Z_{n+1,r}^{2}Z_{i,r+1}^{}Z_{n+1,r+1}^{*} - Z_{n+1,r}^{}Z_{i,r}^{}\right)^{*}}\\
   \qquad\qquad &=& Z_{n+1,r}^{}Q_{i,r}^{}Z_{n+1,r}^{*}-Q_{n+1,r}^{}Q_{i,r+1}^{}.
\end{IEEEeqnarray*}
Therefore
\begin{IEEEeqnarray*}{rCl}
 \IEEEeqnarraymulticol{3}{l}{
    \sum_{i=1}^{r}
    \left(Z_{n+1,r}^{2}Z_{i,r+1}^{}Z_{n+1,r+1}^{*} - Z_{n+1,r}^{}Z_{i,r}^{}\right)
\left(Z_{n+1,r}^{2}Z_{i,r+1}^{}Z_{n+1,r+1}^{*} - Z_{n+1,r}^{}Z_{i,r}^{}\right)^{*}}\\
    \qquad\qquad&=&  Z_{n+1,r}^{}\left(\sum_{i=1}^{r}Q_{i,r}^{}\right)Z_{n+1,r}^{*}
      - Q_{n+1,r}^{}\left(\sum_{i=1}^{r}Q_{i,r+1}^{}\right)\\
    &=& Q_{n+1,r}^{}- Q_{n+1,r}^{}\left(\sum_{k=r+1}^{n+1}P_{r,k}^{}\right)\\
    &=& Q_{n+1,r}^{}P_{r,r}^{}\\
    &=& P_{n+1,r}^{}P_{r,r}^{}\\
    &=& 0.
\end{IEEEeqnarray*}
Thus we have $Z_{i,r}^{}W_{r+1,r}^{} = W_{r+1,r}^{} Z_{i,r}^{}$ 
for $1\leq i\leq r$. Since
$Z_{i,r}^{}$ commutes with $Z_{n+1,j}^{}$ for all $j\geq r+2$, we get
\[
Z_{i,r}^{}W_{n+1,r}^{} = W_{n+1,r}^{} Z_{i,r}^{}.
\]
%%------------------------------------------------------
% \textred{$1\leq i\leq j=r+1$:} \\
%%------------------------------------------------------
Next, let $1\leq i\leq r+1$. 
Then one has
{\allowdisplaybreaks
\begin{IEEEeqnarray*}{rCl}
 Z_{i,r+1}^{}W_{r+2,r}^{} 
   &=& Z_{i,r+1}^{}
       \sum_{k\in\mathbb{N}}\sum_{s\in\mathbb{N}}
         Z_{n+1,r+2}^{k}Z_{n+1,r+1}^{s}Z_{n+1,r}^{}
       (Z_{n+1,r+1}^{s})^{*}(Z_{n+1,r+2}^{k})^{*}\\
   &=& \sum_{k\in\mathbb{N}}\sum_{s\in\mathbb{N}}Z_{n+1,r+2}^{k}Z_{i,r+1}^{}
                   Z_{n+1,r+1}^{s}Z_{n+1,r}^{}
       (Z_{n+1,r+1}^{s})^{*}(Z_{n+1,r+2}^{k})^{*}\\
       && + \sum_{k\in\mathbb{N}}\sum_{s\in\mathbb{N}}Z_{n+1,r+2}^{k}
                Z_{n+1,r+1}^{}Z_{i,r+2}^{}
                   Z_{n+1,r+1}^{s}Z_{n+1,r}^{}
       (Z_{n+1,r+1}^{s})^{*}(Z_{n+1,r+2}^{k+1})^{*}\\
   &=& \sum_{k\in\mathbb{N}}Z_{n+1,r+2}^{k}Z_{i,r+1}^{}
                   Z_{n+1,r}^{}(Z_{n+1,r+2}^{k})^{*}\\
       && + \sum_{k\in\mathbb{N}}\sum_{s\in\mathbb{N}}
               Z_{n+1,r+2}^{k} Z_{n+1,r+1}^{s+1}Z_{i,r+2}^{}
                   Z_{n+1,r}^{}(Z_{n+1,r+1}^{s})^{*}(Z_{n+1,r+2}^{k+1})^{*}\\
   &=& \sum_{k\in\mathbb{N}}Z_{n+1,r+2}^{k}
                   Z_{n+1,r}^{}(Z_{n+1,r+2}^{k})^{*}Z_{i,r+1}^{}\\
       && + \sum_{k\in\mathbb{N}}\sum_{s\in\mathbb{N}}Z_{n+1,r+2}^{k}
                Z_{n+1,r+1}^{}
                   W_{r+1,r}^{}Z_{i,r+2}^{}(Z_{n+1,r+2}^{k+1})^{*},
\end{IEEEeqnarray*}}
and
\begin{IEEEeqnarray*}{rCl}
 W_{r+2,r}^{} Z_{i,r+1}^{}
   &=& \sum_{k\in\mathbb{N}}Z_{n+1,r+2}^{k}
                   Z_{n+1,r}^{}(Z_{n+1,r+2}^{k})^{*}Z_{i,r+1}^{}\\
   && + \sum_{k\in\mathbb{N}}\sum_{s\in\mathbb{N}}
         Z_{n+1,r+2}^{k}Z_{n+1,r+1}^{s+1}Z_{n+1,r}^{}
       (Z_{n+1,r+1}^{s+1})^{*}(Z_{n+1,r+2}^{k})^{*}Z_{i,r+1}^{}\\
   &=& \sum_{k\in\mathbb{N}}Z_{n+1,r+2}^{k}
                   Z_{n+1,r}^{}(Z_{n+1,r+2}^{k})^{*}Z_{i,r+1}^{}\\
   && + \sum_{k\in\mathbb{N}}
         Z_{n+1,r+2}^{k}Z_{n+1,r+1}^{}W_{r+1,r}^{}
       Z_{n+1,r+1}^{*}Z_{i,r+1}^{}(Z_{n+1,r+2}^{k})^{*}.
\end{IEEEeqnarray*}
Using Lemma~\ref{lm:proj-8}, we now get
\begin{IEEEeqnarray*}{rCl}
 \IEEEeqnarraymulticol{3}{l}{Z_{i,r+1}^{}W_{r+2,r}^{} - W_{r+2,r}^{} Z_{i,r+1}^{}}\\
   % \qquad 
   &=& \sum_{k\in\mathbb{N}}
         Z_{n+1,r+2}^{k}Z_{n+1,r+1}^{}W_{r+1,r}^{}
      \left( Z_{i,r+2}^{}Z_{n+1,r+2}^{*} - Z_{n+1,r+1}^{*}Z_{i,r+1}^{} \right)
       Z_{n+1,r+1}^{*}Z_{i,r+1}^{}(Z_{n+1,r+2}^{k})^{*}\\
       &=& 0.
\end{IEEEeqnarray*}
%%------------------------------------------------------
%%------------------------------------------------------
Thus we have  
$Z_{i,r+1}^{}W_{r+2,r}^{} = W_{r+2,r}^{} Z_{i,r+1}^{}$. 
Since $Z_{i,r+1}^{}$ commutes with $Z_{n+1,j}^{}$ for every $j\geq r+3$, we get
\[
Z_{i,r+1}^{}W_{n+1,r}^{} = W_{n+1,r}^{} Z_{i,r+1}^{}.
\]
%%------------------------------------------------------
% \textred{$1\leq i\leq j$, $r+2\leq j$:} \\
%%------------------------------------------------------
Let $j\geq r+2$ and $1\leq i\leq j$. Then one has
\begin{IEEEeqnarray*}{rCl}
 \IEEEeqnarraymulticol{3}{l}{
 Z_{i,j}^{}W_{j+1,r}^{} }\\
   &=& Z_{i,j}^{}
       \sum_{k\in\mathbb{N}}Z_{n+1,j+1}^{k}W_{j,r}^{}(Z_{n+1,j+1}^{*})^{k}\\
   &=& \sum_{k\in\mathbb{N}}Z_{n+1,j+1}^{k}Z_{i,j}^{}
                           W_{j,r}^{}(Z_{n+1,j+1}^{*})^{k}
    +\sum_{k\geq 1}Z_{n+1,j+1}^{k-1}Z_{n+1,j}^{}Z_{i,j+1}^{}
                                       W_{j,r}^{}(Z_{n+1,j+1}^{*})^{k}\\
   &=& \sum_{k\in\mathbb{N}}Z_{n+1,j+1}^{k}W_{j-1,r}^{}
                                   (Z_{n+1,j+1}^{*})^{k}Z_{i,j}^{}
    +\sum_{k\in\mathbb{N}}Z_{n+1,j+1}^{k}Z_{n+1,j}^{}W_{j,r}^{}
                                       Z_{i,j+1}^{}(Z_{n+1,j+1}^{*})^{k+1},
\end{IEEEeqnarray*}
and
{\allowdisplaybreaks
\begin{IEEEeqnarray*}{rCl}
 W_{j+1,r}^{} Z_{i,j}^{}
   &=& \sum_{k\in\mathbb{N}}Z_{n+1,j+1}^{k}W_{j,r}^{}(Z_{n+1,j+1}^{*})^{k}
           Z_{i,j}^{}\\
   &=& \sum_{k\in\mathbb{N}}\sum_{s\in\mathbb{N}}
        Z_{n+1,j+1}^{k}Z_{n+1,j}^{s}
             W_{j-1,r}^{}(Z_{n+1,j}^{*})^{s}(Z_{n+1,j+1}^{*})^{k} Z_{i,j}^{}\\
   &=& \sum_{k\in\mathbb{N}}Z_{n+1,j+1}^{k}W_{j-1,r}^{}
                                   (Z_{n+1,j+1}^{*})^{k}Z_{i,j}^{}\\
   &&\qquad 
    +\sum_{k\in\mathbb{N}}\sum_{s\geq 1}
                Z_{n+1,j+1}^{k}Z_{n+1,j}^{s}W_{j-1,r}^{}
                     (Z_{n+1,j}^{*})^{s}(Z_{n+1,j+1}^{*})^{k}Z_{i,j}^{}\\
   &=& \sum_{k\in\mathbb{N}}Z_{n+1,j+1}^{k}W_{j-1,r}^{}
                                   (Z_{n+1,j+1}^{*})^{k}Z_{i,j}^{}\\
   &&\qquad 
    +\sum_{k\in\mathbb{N}}
                Z_{n+1,j+1}^{k}Z_{n+1,j}^{}W_{j,r}^{}
                     Z_{n+1,j}^{*}(Z_{n+1,j+1}^{*})^{k}Z_{i,j}^{}.
\end{IEEEeqnarray*}
}
Use Lemma~\ref{lm:proj-8} to get
\begin{IEEEeqnarray*}{rCl}
 \IEEEeqnarraymulticol{3}{l}{Z_{i,j}^{}W_{j+1,r}^{} -W_{j+1,r}^{} Z_{i,j}^{}}\\
   \qquad &=& \sum_{k\in\mathbb{N}}Z_{n+1,j+1}^{k}Z_{n+1,j}^{}W_{j,r}^{}
                            Z_{i,j+1}^{}(Z_{n+1,j+1}^{*})^{k+1}\\
   &&\qquad 
       - \sum_{k\in\mathbb{N}}
                Z_{n+1,j+1}^{k}Z_{n+1,j}^{}W_{j,r}^{}
                     Z_{n+1,j}^{*}(Z_{n+1,j+1}^{*})^{k}Z_{i,j}^{}\\
   &=& \sum_{k\in\mathbb{N}}Z_{n+1,j+1}^{k}Z_{n+1,j}^{}W_{j,r}^{}
                 \left(Z_{i,j+1}^{}Z_{n+1,j+1}^{*} - Z_{n+1,j}^{*}Z_{i,j}^{}\right)
                  (Z_{n+1,j+1}^{*})^{k}\\
   &=& 0.
\end{IEEEeqnarray*}
Thus we have $Z_{i,j}^{}W_{j+1,r}^{} = W_{j+1,r}^{} Z_{i,j}^{}$. 
Since $Z_{i,j}^{}$ commutes with $Z_{n+1,k}^{}$ for all $k\geq j+2$, we get
\[
Z_{i,j}^{}W_{n+1,r}^{} = W_{n+1,r}^{} Z_{i,j}^{}.
\]
%%--------------------------------------------------------------------
Thus $W_{n+1,r}^{}$ commutes with $Z_{i,j}^{}$ for all $i\leq j$.
Therefore using (\ref{eq:an0-defrel-9}), we conclude that 
$W_{n+1,r}^{}$ commutes with $Z_{i,j}^{*}$ for all $i> j$.
Since $W_{n+1,r}^{}$ is normal, it follows that it commutes with all
$Z_{i,j}^{}$'s.
\end{proof}
%%--------------------------------------------------------------------

%%------------------------------------------------------
\bppsn %[\texttt{pr:rep-rel-10}]
     %%-------------------------------
     \label{pr:rep-rel-10}
     %%-------------------------------
Let $r$ be as in (\ref{eq:rep-rel-3}). Then
\begin{IEEEeqnarray*}{rCl}
W_{n+1,r}^{*}U_{n+1}^{} &=& U_{n+1}^{*}U_{n+1}^{}.
\end{IEEEeqnarray*}
\eppsn
%%------------------------------------------------------
\begin{proof}
We have
{\allowdisplaybreaks
\begin{IEEEeqnarray*}{rCl}
\IEEEeqnarraymulticol{3}{l}{W_{n+1,r}^{*}U_{n+1}^{} }\\
  &=& \sum_{\substack{k_{j}\in\mathbb{N}\\(r+1\leq j\leq n+1)}}^{}
  \Biggl(
        Z_{n+1,n+1}^{k_{n+1}}\cdots Z_{n+1,r+1}^{k_{r+1}}
        Z_{n+1,r}^{*}
        (Z_{n+1,r+1}^{*})^{k_{r+1}}\cdots 
           (Z_{n+1,n+1}^{*})^{k_{n+1}}\Biggr)\\
  && \quad {}\times
  \sum_{\substack{k_{j}\in\mathbb{N}\\(r+2\leq j\leq n+1)}}^{}
  \Biggl(
        Z_{n+1,n+1}^{k_{n+1}}\cdots Z_{n+1,r+2}^{k_{r+2}}
        Z_{n+1,r}^{}
        (Z_{n+1,r+2}^{*})^{k_{r+2}}\cdots 
           (Z_{n+1,n+1}^{*})^{k_{n+1}}\Biggr)\\
  &=& \sum_{\substack{k_{j}\in\mathbb{N}\\(r+2\leq j\leq n+1)}}^{}
  \Biggl(
        Z_{n+1,n+1}^{k_{n+1}}\cdots Z_{n+1,r+2}^{k_{r+2}}
         \left(\sum_{k_{r+1}\in\mathbb{N}}
             Z_{n+1,r+1}^{k_{r+1}}Z_{n+1,r}^{*}(Z_{n+1,r+1}^{*})^{k_{r+1}}\right)
             Z_{n+1,r}^{}
        \\
      && \qquad\qquad{}\times
        (Z_{n+1,r+2}^{*})^{k_{r+2}}\cdots 
           (Z_{n+1,n+1}^{*})^{k_{n+1}}\Biggr)\\
  &=& \sum_{\substack{k_{j}\in\mathbb{N}\\(r+2\leq j\leq n+1)}}^{}
     \Biggl(
        Z_{n+1,n+1}^{k_{n+1}}\cdots Z_{n+1,r+2}^{k_{r+2}}
        P_{n+1,r}^{}
        (Z_{n+1,r+2}^{*})^{k_{r+2}}\cdots 
           (Z_{n+1,n+1}^{*})^{k_{n+1}}\Biggr),
\end{IEEEeqnarray*}
}
and
{\allowdisplaybreaks
\begin{IEEEeqnarray*}{rCl}
\IEEEeqnarraymulticol{3}{l}{U_{n+1}^{*}U_{n+1}^{} }\\
  &=& \sum_{\substack{k_{j}\in\mathbb{N}\\(r+2\leq j\leq n+1)}}^{}
  \Biggl(
        Z_{n+1,n+1}^{k_{n+1}}\cdots Z_{n+1,r+2}^{k_{r+2}}
        Z_{n+1,r}^{*}
        (Z_{n+1,r+2}^{*})^{k_{r+2}}\cdots 
           (Z_{n+1,n+1}^{*})^{k_{n+1}}\Biggr)\\
  &&  {}\times\sum_{\substack{k_{j}\in\mathbb{N}\\(r+2\leq j\leq n+1)}}^{}
  \Biggl(
        Z_{n+1,n+1}^{k_{n+1}}\cdots Z_{n+1,r+2}^{k_{r+2}}
        Z_{n+1,r}^{}
        (Z_{n+1,r+2}^{*})^{k_{r+2}}\cdots 
           (Z_{n+1,n+1}^{*})^{k_{n+1}}\Biggr)\\
  &=& \sum_{\substack{k_{j}\in\mathbb{N}\\(r+2\leq j\leq n+1)}}^{}
  \Biggl(
        Z_{n+1,n+1}^{k_{n+1}}\cdots Z_{n+1,r+2}^{k_{r+2}}
        % \left(Z_{n+1,r}^{*} Z_{n+1,r}^{}\right)
        P_{n+1,r}^{}
        (Z_{n+1,r+2}^{*})^{k_{r+2}}\cdots 
           (Z_{n+1,n+1}^{*})^{k_{n+1}}\Biggr).
\end{IEEEeqnarray*}}
Thus we have the required equality.
\end{proof}
%%------------------------------------------------------

%%--------------------------------------------------------------------
\section{Implications of irreducibility}
%%--------------------------------------------------------------------
We now further restrict our attention to the images of
the $z_{i,j}^{}$'s under an irreducible 
representation $\pi$. Irreducibility will help us extract more 
information on the operators $W_{i,j}$'s, which will be needed in 
analyzing irreducible representations.

%%--------------------------------------------------------------------
\bppsn %[\texttt{pr:irr-cor-1}]
     %%-------------------------------
     \label{pr:irr-cor-1}
     %%-------------------------------
Assume $\pi$ is irreducible and let $r$ be as in (\ref{eq:rep-rel-3}), i.e.
\[
r=\min\{1\leq i\leq n+1: \pi(z_{n+1,i})\neq 0\}.
\] 
Then there is a $\lambda\in S^{1}$ such that
$W_{n+1,r}^{}=\lambda I$.
\eppsn
%%--------------------------------------------------------------------
\begin{proof}
Since $W_{n+1,r}^{}$ is normal, by Theorem~\ref{th:rep-rel-9} and by 
the irreducibility of $\pi$ it is a scalar operator. But it is also
a partial isometry and 
\[
W_{n+1,r}^{*}W_{n+1,r}^{}\geq P_{n+1,r}\neq 0.
\]
Therefore we have the result.
\end{proof}
%%--------------------------------------------------------------------

Our next goal is to show that for an irreducible representation $\pi$,
the operator $W_{n+1,r+1}^{}$ is an isometry. For that, we start with
a few commutation relations.

%%--------------------------------------------------------------------
\bppsn %[\texttt{pr:wcomm-1a}]
     %%-------------------------------
     \label{pr:wcomm-1a}
     %%-------------------------------
Let $1\leq j\leq r-1$ and $1\leq i\leq n$. Then
\begin{IEEEeqnarray*}{rCl}
 W_{n+1,r+1}^{}Z_{i,j}^{} &=& Z_{i,j}^{}W_{n+1,r+1}^{}, \\
 W_{n+1,r+1}^{*}Z_{i,j}^{} &=& Z_{i,j}^{}W_{n+1,r+1}^{*}
\end{IEEEeqnarray*}
\eppsn
%%--------------------------------------------------------------------
\begin{proof}
This is an immediate consequence of Proposition~\ref{pr:rep-rel-4}.
\end{proof}
%%--------------------------------------------------------------------

%%--------------------------------------------------------------------
\bppsn %[\texttt{pr:wcomm-1}]
     %%-------------------------------
     \label{pr:wcomm-1}
     %%-------------------------------
\begin{IEEEeqnarray}{rCll}
 W_{i,r+1}^{}Z_{n+1,j}^{} &=& 0, &\quad  r+2\leq i<j.\\
 Z_{n+1,j}^{}W_{i,r+1}^{} &=& W_{i,r+1}^{}Z_{n+1,j}^{},&\quad j\leq i\leq n+1.\\
 Z_{i,j}^{}W_{j,r+1}^{} &=& Z_{i,j}^{}W_{j-1,r+1}^{} &\nonumber\\
     &=& W_{j-1,r+1}^{}Z_{i,j}^{},&\quad i\leq j,\;r+2\leq j.
\end{IEEEeqnarray}
\eppsn
%%--------------------------------------------------------------------
\begin{proof}
The first equality follows from the fact that $Q_{n+1,i}^{}Q_{n+1,j}^{}=0$ 
for $i\neq j$.

For the second equality, observe that if $j<i$, then
\begin{IEEEeqnarray*}{rCl}
 Z_{n+1,j}^{}W_{i,r+1}^{} 
   &=& Z_{n+1,j}^{}
       \sum_{k\in\mathbb{N}}Z_{n+1,i}^{k}W_{i-1,r+1}^{}(Z_{n+1,i}^{k})^{*}\\
   &=& Z_{n+1,j}^{}W_{i-1,r+1}^{},\\
 %%   &=& W_{i-1,r+1}^{}Z_{n+1,j}^{} % WHY?
 W_{i,r+1}^{}Z_{n+1,j}^{}
   &=& \sum_{k\in\mathbb{N}}Z_{n+1,i}^{k}W_{i-1,r+1}^{}(Z_{n+1,i}^{k})^{*}
          Z_{n+1,j}^{}\\
   &=& W_{i-1,r+1}^{}Z_{n+1,j}^{}.
\end{IEEEeqnarray*}
By repeating the above, we get 
$Z_{n+1,j}^{}W_{i,r+1}^{}=Z_{n+1,j}^{}W_{j,r+1}^{}$
and
$W_{i,r+1}^{}Z_{n+1,j}^{}=W_{j,r+1}^{}Z_{n+1,j}^{}$.
Finally, for $j=i$ one has
\begin{IEEEeqnarray*}{rCl}
 Z_{n+1,j}^{}W_{j,r+1}^{} 
   &=& Z_{n+1,j}^{}
       \sum_{k\in\mathbb{N}}Z_{n+1,j}^{k}W_{j-1,r+1}^{}(Z_{n+1,j}^{k})^{*}\\
   &=& \sum_{k\in\mathbb{N}}Z_{n+1,j}^{k+1}W_{j-1,r+1}^{}(Z_{n+1,j}^{k})^{*}\\
   &=& Z_{n+1,j}^{}W_{j-1,r+1}^{} +
       \sum_{k\geq 1}Z_{n+1,j}^{k+1}W_{j-1,r+1}^{}(Z_{n+1,j}^{k})^{*}\\
   &=& Z_{n+1,j}^{}W_{j-1,r+1}^{} Z_{n+1,j}^{*}Z_{n+1,j}^{}+
       \sum_{k\geq 1}Z_{n+1,j}^{k+1}W_{j-1,r+1}^{}(Z_{n+1,j}^{k+1})^{*}Z_{n+1,j}^{}\\
   &=& \sum_{k\geq 1}Z_{n+1,j}^{k}W_{j-1,r+1}^{}(Z_{n+1,j}^{k})^{*}Z_{n+1,j}^{},\\
   &=& W_{j,r+1}^{} Z_{n+1,j}^{} - W_{j-1,r+1}^{}Z_{n+1,j}^{}\\
   &=& W_{j,r+1}^{} Z_{n+1,j}^{}.
\end{IEEEeqnarray*}
The last equation follows from the observation that $Z_{i,i}^{}Z_{n+1,i}^{}=0$.
\end{proof}

%%------------------------------------------------------
\bppsn %[\texttt{pr:wcomm-2}]
     %%-------------------------------
     \label{pr:wcomm-2}
     %%-------------------------------
Let $j\geq r+2$ and $1\leq i\leq j$. Then
\begin{IEEEeqnarray*}{rCl}
 Z_{i,j}^{}W_{j+1,r+1}^{} 
   &=& W_{j+1,r+1}^{} Z_{i,j}^{},\\
 Z_{i,j}^{}W_{j+1,r+1}^{*} 
   &=& W_{j+1,r+1}^{} Z_{i,j}^{*}.
\end{IEEEeqnarray*}
\eppsn
%%------------------------------------------------------
\begin{proof}
Let $j\geq r+2$ and $1\leq i\leq j$. Then one has
\begin{IEEEeqnarray*}{rCl}
 Z_{i,j}^{}W_{j+1,r+1}^{} 
   &=& Z_{i,j}^{}
       \sum_{k\in\mathbb{N}}Z_{n+1,j+1}^{k}W_{j,r+1}^{}(Z_{n+1,j+1}^{*})^{k}\\
   &=& \sum_{k\in\mathbb{N}}Z_{n+1,j+1}^{k}Z_{i,j}^{}
                           W_{j,r+1}^{}(Z_{n+1,j+1}^{*})^{k}\\
   &&\qquad 
    +\sum_{k\geq 1}Z_{n+1,j+1}^{k-1}Z_{n+1,j}^{}Z_{i,j+1}^{}
                                       W_{j,r+1}^{}(Z_{n+1,j+1}^{*})^{k}\\
   &=& \sum_{k\in\mathbb{N}}Z_{n+1,j+1}^{k}W_{j-1,r+1}^{}
                                   (Z_{n+1,j+1}^{*})^{k}Z_{i,j}^{}\\
   &&\qquad +\sum_{k\in\mathbb{N}}Z_{n+1,j+1}^{k}Z_{n+1,j}^{}W_{j,r+1}^{}
                                       Z_{i,j+1}^{}(Z_{n+1,j+1}^{*})^{k+1},
\end{IEEEeqnarray*}
and
\begin{IEEEeqnarray*}{rCl}
 W_{j+1,r+1}^{} Z_{i,j}^{}
   &=& \sum_{k\in\mathbb{N}}Z_{n+1,j+1}^{k}W_{j,r+1}^{}(Z_{n+1,j+1}^{*})^{k}
           Z_{i,j}^{}\\
   &=& \sum_{k\in\mathbb{N}}\sum_{s\in\mathbb{N}}
        Z_{n+1,j+1}^{k}Z_{n+1,j}^{s}
             W_{j-1,r+1}^{}(Z_{n+1,j}^{*})^{s}(Z_{n+1,j+1}^{*})^{k} Z_{i,j}^{}\\
   &=& \sum_{k\in\mathbb{N}}Z_{n+1,j+1}^{k}W_{j-1,r+1}^{}
                                   (Z_{n+1,j+1}^{*})^{k}Z_{i,j}^{}\\
   &&\qquad 
    +\sum_{k\in\mathbb{N}}\sum_{s\geq 1}
                Z_{n+1,j+1}^{k}Z_{n+1,j}^{s}W_{j-1,r+1}^{}
                     (Z_{n+1,j}^{*})^{s}(Z_{n+1,j+1}^{*})^{k}Z_{i,j}^{}\\
   &=& \sum_{k\in\mathbb{N}}Z_{n+1,j+1}^{k}W_{j-1,r+1}^{}
                                   (Z_{n+1,j+1}^{*})^{k}Z_{i,j}^{}\\
   &&\qquad 
    +\sum_{k\in\mathbb{N}}
                Z_{n+1,j+1}^{k}Z_{n+1,j}^{}W_{j,r+1}^{}
                     (Z_{n+1,j}^{*})^{}(Z_{n+1,j+1}^{*})^{k}Z_{i,j}^{}.
\end{IEEEeqnarray*}
Use Lemma~\ref{lm:proj-8} to get
\begin{IEEEeqnarray*}{rCl}
 \IEEEeqnarraymulticol{3}{l}{Z_{i,j}^{}W_{j+1,r+1}^{} -W_{j+1,r+1}^{} Z_{i,j}^{}}\\
   \qquad &=& \sum_{k\in\mathbb{N}}Z_{n+1,j+1}^{k}Z_{n+1,j}^{}W_{j,r+1}^{}
                            Z_{i,j+1}^{}(Z_{n+1,j+1}^{*})^{k+1}\\
   &&\qquad 
       - \sum_{k\in\mathbb{N}}
                Z_{n+1,j+1}^{k}Z_{n+1,j}^{}W_{j,r+1}^{}
                     (Z_{n+1,j}^{*})^{}(Z_{n+1,j+1}^{*})^{k}Z_{i,j}^{}\\
   &=& \sum_{k\in\mathbb{N}}Z_{n+1,j+1}^{k}Z_{n+1,j}^{}W_{j,r+1}^{}
                 \left(Z_{i,j+1}^{}Z_{n+1,j+1}^{*} - Z_{n+1,j}^{*}Z_{i,j}^{}\right)
                  (Z_{n+1,j+1}^{*})^{k}\\
    &=& 0.
\end{IEEEeqnarray*}
Thus we have the required relation.
\end{proof}
%%--------------------------------------------------------------------

%%------------------------------------------------------
\bppsn %[\texttt{pr:wcomm-3a}]
     %%-------------------------------
     \label{pr:wcomm-3a}
     %%-------------------------------
    Let $1\leq i\leq r$. Then $Z_{i,r}^{}$ commutes with 
    $W_{n+1,r+1}^{*}W_{n+1,r+1}^{}$.
\eppsn
%%------------------------------------------------------
\begin{proof}
From (\ref{eq:rep-rel-7}), we have for $r+2\leq j\leq n+1$,
\begin{IEEEeqnarray*}{rCl}
W_{j,r+1}^{*}W_{j,r+1}^{} 
    &=& \!\!\!\sum_{\substack{k_{j}\in\mathbb{N}\\(r+2\leq j\leq n+1)}}^{}
        \hspace{-1em}
        Z_{n+1,j}^{k_{j}}\ldots Z_{n+1,r+2}^{k_{r+2}}
        \left(Z_{n+1,r+1}^{*}Z_{n+1,r+1}^{}\right)
        (Z_{n+1,r+2}^{k_{r+2}})^{*}\ldots (Z_{n+1,j}^{k_{j}})^{*}.
\end{IEEEeqnarray*}
Hence
\begin{IEEEeqnarray*}{rCl}
\IEEEeqnarraymulticol{3}{l}{
   Z_{i,r}^{}W_{n+1,r+1}^{*}W_{n+1,r+1}^{} - 
      W_{n+1,r+1}^{*}W_{n+1,r+1}^{} Z_{i,r}^{}}\\
    \qquad &=& \!\!\!\sum_{\substack{k_{j}\in\mathbb{N}\\(r+2\leq j\leq n+1)}}^{}
        \Biggl(Z_{n+1,n+1}^{k_{n+1}}\ldots Z_{n+1,r+2}^{k_{r+2}}\\
   &&\quad {}  \left(Z_{i,r}^{}Z_{n+1,r+1}^{*}Z_{n+1,r+1}^{}
          - Z_{n+1,r+1}^{*}Z_{n+1,r+1}^{}Z_{i,r}^{}\right)
        (Z_{n+1,r+2}^{k_{r+2}})^{*}\ldots (Z_{n+1,n+1}^{k_{n+1}})^{*}\Biggr).
\end{IEEEeqnarray*}
since $Q_{n+1,r+1}^{}Q_{n+1,r}^{}=0$, we have
\begin{IEEEeqnarray*}{rCl}
  Z_{i,r}^{}Z_{n+1,r+1}^{*}Z_{n+1,r+1}^{}
   &=& Z_{n+1,r+1}^{*}Z_{i,r}^{}Z_{n+1,r+1}^{}\\
   &=& Z_{n+1,r+1}^{*}
     \left(Z_{n+1,r+1}^{}Z_{i,r}^{}+Z_{n+1,r}^{}Z_{i,r+1}^{}\right)\\
   &=& Z_{n+1,r+1}^{*}Z_{n+1,r+1}^{}Z_{i,r}^{}.
\end{IEEEeqnarray*}
Therefore we have the equality
$Z_{i,r}^{}W_{n+1,r+1}^{*}W_{n+1,r+1}^{} = W_{n+1,r+1}^{*}W_{n+1,r+1}^{} Z_{i,r}^{}$.
\end{proof}
%%------------------------------------------------------

%%------------------------------------------------------
\bppsn %[\texttt{pr:wcomm-3b}]
     %%-------------------------------
     \label{pr:wcomm-3b}
     %%-------------------------------
    Let $1\leq i\leq r+1$. Then $Z_{i,r+1}^{}$ commutes with $W_{n+1,r+1}^{*}W_{n+1,r+1}^{}$.
\eppsn
%%------------------------------------------------------
\begin{proof}
We have
\begin{IEEEeqnarray*}{rCl}
\IEEEeqnarraymulticol{3}{l}{
   Z_{i,r+1}^{}W_{r+2,r+1}^{*}W_{r+2,r+1}^{} }\\
    \qquad\qquad &=& Z_{i,r+1}^{}\sum_{k\in\mathbb{N}}
        Z_{n+1,r+2}^{k}\left(Z_{n+1,r+1}^{*}Z_{n+1,r+1}^{}\right)
        (Z_{n+1,r+2}^{k})^{*}\\
    &=& \sum_{k\in\mathbb{N}}
        Z_{n+1,r+2}^{k}\left(Z_{i,r+1}^{}Z_{n+1,r+1}^{*}Z_{n+1,r+1}^{}\right)
        (Z_{n+1,r+2}^{k})^{*}\\
    && \qquad {} + \sum_{k\in\mathbb{N}}
        Z_{n+1,r+2}^{k}\left(Z_{n+1,r+1}^{}Z_{i,r+2}^{}Z_{n+1,r+1}^{*}Z_{n+1,r+1}^{}\right)
        (Z_{n+1,r+2}^{k+1})^{*}\\
    &=& \sum_{k\in\mathbb{N}}
        Z_{n+1,r+2}^{k}Z_{i,r+1}^{}\left(Q_{n+1,r+1}^{}+Q_{n+1,r}^{}\right)
        (Z_{n+1,r+2}^{k})^{*}\\
    && \qquad {} + \sum_{k\in\mathbb{N}}
        Z_{n+1,r+2}^{k}\left(Z_{n+1,r+1}^{}Z_{i,r+2}^{}Z_{n+1,r+2}^{*}\right)
        (Z_{n+1,r+2}^{k})^{*}\\
    &=& \sum_{k\in\mathbb{N}}
        Z_{n+1,r+2}^{k}Z_{i,r+1}^{}Q_{n+1,r}^{}
        (Z_{n+1,r+2}^{k})^{*}\\
    && \qquad {} + \sum_{k\in\mathbb{N}}
        Z_{n+1,r+2}^{k}\left(Z_{n+1,r+1}^{}Z_{i,r+2}^{}Z_{n+1,r+2}^{*}\right)
        (Z_{n+1,r+2}^{k})^{*},
\end{IEEEeqnarray*}
and
{\allowdisplaybreaks
\begin{IEEEeqnarray*}{rCl}
\IEEEeqnarraymulticol{3}{l}{
   W_{r+2,r+1}^{*}W_{r+2,r+1}^{} Z_{i,r+1}^{}}\\
    \qquad\qquad &=& \left(\sum_{k\in\mathbb{N}}
        Z_{n+1,r+2}^{k}\left(Z_{n+1,r+1}^{*}Z_{n+1,r+1}^{}\right)
        (Z_{n+1,r+2}^{k})^{*}\right)Z_{i,r+1}^{}\\
    &=& \left(\sum_{k\in\mathbb{N}}
        Z_{n+1,r+2}^{k}\left(Z_{n+1,r}^{}Z_{n+1,r}^{*}+
            Z_{n+1,r+1}^{}Z_{n+1,r+1}^{*}\right)Z_{i,r+1}^{}
        (Z_{n+1,r+2}^{k})^{*}\right)\\
    &=& \sum_{k\in\mathbb{N}}
        Z_{n+1,r+2}^{k}\left(Q_{n+1,r}^{}Z_{i,r+1}^{}\right)
        (Z_{n+1,r+2}^{k})^{*}\\
    && \qquad {} + \sum_{k\in\mathbb{N}}
        Z_{n+1,r+2}^{k}\left(Z_{n+1,r+1}^{}Z_{n+1,r+1}^{*}Z_{i,r+1}^{}\right)
        (Z_{n+1,r+2}^{k+1})^{*}.
\end{IEEEeqnarray*}
}
Therefore 
\begin{IEEEeqnarray*}{rCl}
\IEEEeqnarraymulticol{3}{l}{
   Z_{i,r+1}^{}W_{r+2,r+1}^{*}W_{r+2,r+1}^{} 
     - W_{r+2,r+1}^{*}W_{r+2,r+1}^{} Z_{i,r+1}^{}}\\
      &=& \sum_{k\in\mathbb{N}}
        Z_{n+1,r+2}^{k}\left(Z_{i,r+1}^{}Q_{n+1,r}^{}-Q_{n+1,r}^{}Z_{i,r+1}^{}\right)
        (Z_{n+1,r+2}^{k})^{*}\\
    && \qquad {} + \sum_{k\in\mathbb{N}}
        Z_{n+1,r+2}^{k}
        Z_{n+1,r+1}^{}\left(Z_{i,r+2}^{}Z_{n+1,r+2}^{*}-Z_{n+1,r+1}^{*}Z_{i,r+1}^{}\right)
        (Z_{n+1,r+2}^{k})^{*}\\
    &=& \sum_{k\in\mathbb{N}}
        Z_{n+1,r+2}^{k}
        Z_{n+1,r+1}^{}\left(Z_{i,r+2}^{}Z_{n+1,r+2}^{*}-Z_{n+1,r+1}^{*}Z_{i,r+1}^{}\right)
        (Z_{n+1,r+2}^{k})^{*}.
\end{IEEEeqnarray*}
The right hand side above vanishes by Lemma~\ref{lm:proj-8}. 
Since $Z_{i,r+1}^{}$ commutes with $Z_{n+1,j}^{}$ and 
$Z_{n+1,j}^{*}$ for $j\geq r+3$, it follows that $Z_{i,r+1}^{}$ 
commutes with $W_{n+1,r+1}^{*}W_{n+1,r+1}^{}$.
\end{proof}
%%------------------------------------------------------

%%------------------------------------------------------
\bppsn %[\texttt{pr:wcomm-3c}]
     %%-------------------------------
     \label{pr:wcomm-3c}
     %%-------------------------------
    Let $\pi$ be irreducible. Then  $W_{n+1,r+1}^{*}W_{n+1,r+1}^{}=I$.
\eppsn
%%------------------------------------------------------
\begin{proof}
From Propositions \ref{pr:wcomm-1a}, \ref{pr:wcomm-2}, \ref{pr:wcomm-3a} 
and \ref{pr:wcomm-3b},
it follows that $W_{n+1,r+1}^{*}W_{n+1,r+1}^{}$ commutes with $Z_{i,j}^{}$ 
for all $1\leq i\leq j\leq n+1$. Therefore using (\ref{eq:an0-defrel-9})
we conclude that $W_{n+1,r+1}^{*}W_{n+1,r+1}^{}$ commutes with $Z_{i,j}^{}$
for all $i$ and $j$. Since 
$W_{n+1,r+1}^{*}W_{n+1,r+1}^{}$ is a projection, and
\[
W_{n+1,r+1}^{*}W_{n+1,r+1}^{}\geq P_{n+1,r+1}^{}\neq 0,
\]
we get the result.
\end{proof}
%%------------------------------------------------------

The above result and Proposition~\ref{pr:irr-cor-1} tell us that the operator
$W_{n+1,r}^{*}W_{n+1,r+1}^{}$ is also an isometry.
We next  decompose the operators $Z_{i,j}^{}$
using this isometry. This will be the main tool in proving
the theorems in the next section.
%%------------------------------------------------------
\bthm %[\texttt{th:factorization-2}]
     %%-------------------------------
     \label{th:factorization-2}
     %%-------------------------------
Assume $\pi$ is irreducible. Define
\begin{IEEEeqnarray}{rClCl}
  W &:=& W_{n+1,r}^{*}W_{n+1,r+1}^{}
    &=& \bar{\lambda}W_{n+1,r+1}^{},
     %%-------------------------------
     \label{eq:wdef}
     %%-------------------------------
\end{IEEEeqnarray}
%%-------------------------------------
{\allowdisplaybreaks
\begin{IEEEeqnarray}{rCl}
  Y_{i,j}^{(1)} & = &
     \begin{cases}
     Z_{i,j}^{} & \text{if }j\not\in\{r,r+1\},\\
     Z_{i,r}^{} W&  \text{if }j=r,\\
     W^{*}Z_{i,r+1}^{} &  \text{if }j=r+1.
     \end{cases}\\
     &&\nonumber \\
     &&\nonumber \\
  Y_{i,j}^{(2)} &=& 
    \begin{cases}
     W & \text{if }i=j=r+1,\\
     W_{}^{*} & \text{if }i=j=r,\\
     I-WW_{}^{*} & \text{if $i,j\in \{r,r+1\}$ and $i\neq j$},\\
     I & \text{if }i=j\not\in\{r,r+1\},\\
     0 & \text{otherwise.}
    \end{cases}
\end{IEEEeqnarray}}
%%-------------------------------------
Then $Y_{i,j}^{(1)}\in \{W,W_{}^{*}\}^{\prime}$ for all $i,j$ and
{\allowdisplaybreaks
\begin{IEEEeqnarray}{rCl}
Z_{i,j}^{} &=& 
   \sum_{k=\min\{i,j\}}^{\max\{i,j\}}Y_{i,k}^{(1)}Y_{k,j}^{(2)}.
     %%-------------------------------
     \label{eq:decomp}
     %%-------------------------------
\end{IEEEeqnarray}}
\ethm
%%------------------------------------------------------

%%--------------------------------------------------------------------

Before getting into the proof of the above theorem, we first prove two
propositions where we show that the elements $Z_{i,r}^{}W$ and $W^{*}Z_{i,r+1}^{}$  
are in the commutant $\{W,W_{}^{*}\}^{\prime}$.
%%--------------------------------------------------------------------
\bppsn %[\texttt{pr:wcomm-7}]
     %%-------------------------------
     \label{pr:wcomm-7}
     %%-------------------------------
Let $W$ be as defined in Theorem~\ref{th:factorization-2}. Then
$Z_{i,r}^{}W\in \{W,W_{}^{*}\}^{\prime}$
for $1\leq i\leq r$.
\eppsn
%%--------------------------------------------------------------------
\begin{proof}
It is enough to show that
$Z_{i,r}^{}W_{n+1,r+1}^{}\in \{W_{n+1,r+1}^{},W_{n+1,r+1}^{*}\}^{\prime}$
From Proposition~\ref{pr:rep-rel-6}, we have
\begin{IEEEeqnarray*}{rCl}
 (Z_{i,r}^{}W_{n+1,r+1}^{}-W_{n+1,r+1}^{}Z_{i,r}^{})W_{n+1,r+1}^{} 
    &=&U_{n+1}^{}Z_{i,r+1}^{}W_{n+1,r+1}^{}.
\end{IEEEeqnarray*}
Now note that
\begin{IEEEeqnarray*}{rCl}
\IEEEeqnarraymulticol{3}{l}{U_{n+1}Z_{i,r+1}^{}W_{n+1,r+1}^{}}\\
 &=& 
  % \sum_{k_{r+2},\ldots,k_{n+1}\in\mathbb{N}}^{}\Biggl(
  \hspace{-1em}
  \sum_{\substack{k_{j}\in\mathbb{N}\\(r+2\leq j\leq n+1)}}^{}
  \hspace{-1em}\Biggl(
     Z_{n+1,n+1}^{k_{n+1}}%Z_{n+1,n}^{k_{n}}
     \cdots Z_{n+1,r+2}^{k_{r+2}}
        Z_{n+1,r}^{}(Z_{n+1,r+2}^{*})^{k_{r+2}}\cdots %(Z_{n+1,n}^{*})^{k_{n}}
           (Z_{n+1,n+1}^{*})^{k_{n+1}}\Biggr)\\
        &&\hspace{2em} Z_{i,r+1}^{} \hspace{-1em}
           \sum_{\substack{k_{j}\in\mathbb{N}\\(r+2\leq j\leq n+1)}}^{}
           \hspace{-1em}\Biggl(
           Z_{n+1,n+1}^{k_{n+1}}%Z_{n+1,n}^{k_{n}}
           \cdots Z_{n+1,r+2}^{k_{r+2}}
        Z_{n+1,r+1}^{}(Z_{n+1,r+2}^{*})^{k_{r+2}}\cdots %(Z_{n+1,n}^{*})^{k_{n}}
           (Z_{n+1,n+1}^{*})^{k_{n+1}}\Biggr)\\
      &=& \hspace{-1em}
      \sum_{\substack{k_{j}\in\mathbb{N}\\(r+2\leq j\leq n+1)}}^{}\hspace{-1em}
      Z_{n+1,n+1}^{k_{n+1}} \cdots Z_{n+1,r+2}^{k_{r+2}}
      \left(Z_{n+1,r}^{}Z_{i,r+1}^{}Z_{n+1,r+1}^{}\right)
      (Z_{n+1,r+2}^{*})^{k_{r+2}}\cdots %(Z_{n+1,n}^{*})^{k_{n}}
           (Z_{n+1,n+1}^{*})^{k_{n+1}}\\
      &=& \hspace{-1em}
      \sum_{\substack{k_{j}\in\mathbb{N}\\(r+2\leq j\leq n+1)}}^{}\hspace{-1em}
      Z_{n+1,n+1}^{k_{n+1}} \cdots Z_{n+1,r+2}^{k_{r+2}}
           \left(Z_{i,r+1}^{}Z_{n+1,r}^{}Z_{n+1,r+1}^{}\right)
      (Z_{n+1,r+2}^{*})^{k_{r+2}}\cdots 
           (Z_{n+1,n+1}^{*})^{k_{n+1}}\\
      &=&  0.
     %%-------------------------------
     \yesnumber\label{eq:wcomm-7pr}\\
     %%-------------------------------
\end{IEEEeqnarray*}
Thus $Z_{i,r}^{}W_{n+1,r+1}^{}$ commutes with $W_{n+1,r+1}^{}$.
One also has
\begin{IEEEeqnarray*}{rClCl}
 (Z_{i,r}^{}W_{n+1,r+1}^{})W_{n+1,r+1}^{*}
    &=&Z_{i,r}^{}(I-U_{n+1}^{}U_{n+1}^{*})
    &=& Z_{i,r}^{},\\
  W_{n+1,r+1}^{*}(Z_{i,r}^{}W_{n+1,r+1}^{})
    &=& Z_{i,r}^{}W_{n+1,r+1}^{*}W_{n+1,r+1}^{}
    &=& Z_{i,r}^{}.
\end{IEEEeqnarray*}
Thus 
$Z_{i,r}^{}W_{n+1,r+1}^{}\in \{W_{n+1,r+1}^{},W_{n+1,r+1}^{*}\}^{\prime}$
for $1\leq i\leq r$.
\end{proof}
%%--------------------------------------------------------------------

%%------------------------------------------------------
\bppsn %[\texttt{pr:wcomm-9}]
     %%-------------------------------
     \label{pr:wcomm-9}
     %%-------------------------------
Let $i\leq r+1 $. Then $W_{}^{*}Z_{i,r+1}^{}\in \{W_{}^{},W_{}^{*}\}^{\prime}$.
\eppsn
%%------------------------------------------------------
\begin{proof}
Let us first prove that 
$W_{n+1,r+1}^{*}Z_{i,r+1}^{}W_{n+1,r+1}^{}=
   W_{n+1,r+1}^{}W_{n+1,r+1}^{*}Z_{i,r+1}^{}$. \\
Observe that
\begin{IEEEeqnarray*}{rCl}
\IEEEeqnarraymulticol{3}{l}{
  Z_{i,r+1}^{}W_{r+2,r+1}^{}}\\
  \qquad &=& \sum_{k\in\mathbb{N}}Z_{n+1,r+2}^{k}
      Z_{i,r+1}^{}Z_{n+1,r+1}^{}(Z_{n+1,r+2}^{*})^{k}
  + \sum_{k\geq 1}Z_{n+1,r+2}^{k-1}Z_{i,r+2}^{}Z_{n+1,r+1}^{2}(Z_{n+1,r+2}^{*})^{k}\\
  &=& \sum_{k\in\mathbb{N}}
         Z_{n+1,r+2}^{k}Z_{i,r+2}^{}Z_{n+1,r+1}^{2}(Z_{n+1,r+2}^{*})^{k+1}.
\end{IEEEeqnarray*}
Therefore
\begin{IEEEeqnarray*}{rCl}
\IEEEeqnarraymulticol{3}{l}{
  W_{n+1,r+1}^{*}Z_{i,r+1}^{}W_{n+1,r+1}^{}}\\
  &=& \hspace{-1em}
  \sum_{\substack{k_{j}\in\mathbb{N}\\(r+2\leq j\leq n+1)}}^{}\hspace{-1em}
        Z_{n+1,n+1}^{k_{n+1}}\ldots Z_{n+1,r+2}^{k_{r+2}}
      \left(Z_{n+1,r+1}^{*}Z_{i,r+2}^{}Z_{n+1,r+1}^{2}Z_{n+1,r+2}^{*}\right)
        (Z_{n+1,r+2}^{k_{r+2}})^{*}\ldots (Z_{n+1,n+1}^{k_{n+1}})^{*}\\
  &=& \hspace{-1em}
  \sum_{\substack{k_{j}\in\mathbb{N}\\(r+2\leq j\leq n+1)}}^{}\hspace{-1em}
  Z_{n+1,n+1}^{k_{n+1}}\ldots Z_{n+1,r+2}^{k_{r+2}}
     \left( Z_{n+1,r+1}^{*}Z_{n+1,r+1}^{2}Z_{i,r+2}^{}Z_{n+1,r+2}^{*}\right)
        (Z_{n+1,r+2}^{k_{r+2}})^{*}\ldots (Z_{n+1,n+1}^{k_{n+1}})^{*}\\
  &=& \hspace{-1em}
  \sum_{\substack{k_{j}\in\mathbb{N}\\(r+2\leq j\leq n+1)}}^{}\hspace{-1em}
      Z_{n+1,n+1}^{k_{n+1}}\ldots Z_{n+1,r+2}^{k_{r+2}}
      \left(P_{n+1,r+1}^{}Z_{n+1,r+1}^{}Z_{i,r+2}^{}Z_{n+1,r+2}^{*}\right)
        (Z_{n+1,r+2}^{k_{r+2}})^{*}\ldots (Z_{n+1,n+1}^{k_{n+1}})^{*}\\
  &=& \hspace{-1em}
  \sum_{\substack{k_{j}\in\mathbb{N}\\(r+2\leq j\leq n+1)}}^{}\hspace{-1em}
       Z_{n+1,n+1}^{k_{n+1}}\ldots Z_{n+1,r+2}^{k_{r+2}}
           \left(  Z_{n+1,r+1}^{}Z_{i,r+2}^{}Z_{n+1,r+2}^{*}\right)
        (Z_{n+1,r+2}^{k_{r+2}})^{*}\ldots (Z_{n+1,n+1}^{k_{n+1}})^{*}.
\end{IEEEeqnarray*}
We also have
\begin{IEEEeqnarray*}{rCl}
\IEEEeqnarraymulticol{3}{l}{
      W_{n+1,r+1}^{}W_{n+1,r+1}^{*}Z_{i,r+1}^{}}\\
   &=& \hspace{-1em}
   \sum_{\substack{k_{j}\in\mathbb{N}\\(r+2\leq j\leq n+1)}}^{}
   \hspace{-1em}
     Z_{n+1,n+1}^{k_{n+1}}\ldots Z_{n+1,r+2}^{k_{r+2}}
      \left(Z_{n+1,r+1}^{}Z_{n+1,r+1}^{*}Z_{i,r+1}^{}\right)
        (Z_{n+1,r+2}^{k_{r+2}})^{*}\ldots (Z_{n+1,n+1}^{k_{n+1}})^{*}.
\end{IEEEeqnarray*}
Hence using Lemma~\ref{lm:proj-8}, we get
\begin{IEEEeqnarray*}{rCl}
\IEEEeqnarraymulticol{3}{l}{
  W_{n+1,r+1}^{*}Z_{i,r+1}^{}W_{n+1,r+1}^{}-
  W_{n+1,r+1}^{}W_{n+1,r+1}^{*}Z_{i,r+1}^{}}\\
  &=& \sum_{\substack{k_{j}\in\mathbb{N}\\(r+2\leq j\leq n+1)}}^{}%{k_{r+2},\ldots,k_{n+1}\in\mathbb{N}}^{}%\hspace{-1em}
      \hspace{-1em}\Biggl(
      Z_{n+1,n+1}^{k_{n+1}}\ldots Z_{n+1,r+2}^{k_{r+2}}
      Z_{n+1,r+1}^{}\\
   &&\hspace{50pt}
       % {}\times 
       \left(Z_{i,r+2}^{}Z_{n+1,r+2}^{*}
       - Z_{n+1,r+1}^{*}Z_{i,r+1}^{}\right)%\\
   % &&\qquad\qquad
        (Z_{n+1,r+2}^{k_{r+2}})^{*}\ldots (Z_{n+1,n+1}^{k_{n+1}})^{*}\Biggr)\\
    &=& 0.
\end{IEEEeqnarray*}
Since $W$ commutes with all the $Z_{i,j}^{}$'s, we get 
$W^{*}Z_{i,r+1}^{}W = WW_{}^{*}Z_{i,r+1}^{}$.
It remains to prove that
$W^{*}Z_{i,r+1}^{}W^{*}=(W^{*})^{2}Z_{i,r+1}^{}$.
Since $W_{n+1,r+1}^{}$ is an isometry, we have
\begin{IEEEeqnarray*}{rCl}
  W_{n+1,r+1}^{*}Z_{i,r+1}^{}W_{n+1,r+1}^{*}W_{n+1,r+1}^{} &=& W_{n+1,r+1}^{*}Z_{i,r+1}^{}\\
      &=& W_{n+1,r+1}^{*}W_{n+1,r+1}^{}W_{n+1,r+1}^{*}Z_{i,r+1}^{}\\
      &=& (W_{n+1,r+1}^{*})^{2}Z_{i,r+1}^{}W.
\end{IEEEeqnarray*}
Therefore it is enough to prove that
\begin{IEEEeqnarray*}{rCl}
\IEEEeqnarraymulticol{3}{l}{
W_{n+1,r+1}^{*}Z_{i,r+1}^{}W_{n+1,r+1}^{*}(I-W_{n+1,r+1}^{}W_{n+1,r+1}^{*})}\\
 \qquad\qquad\qquad  &=&  
  (W_{n+1,r+1}^{*})^{2}Z_{i,r+1}^{}(I-W_{n+1,r+1}^{}W_{n+1,r+1}^{*}).
\end{IEEEeqnarray*}
Note that the left hand side above is 0.
Since by (\ref{eq:rep-rel-13}) we have 
$I-W_{n+1,r+1}^{}W_{n+1,r+1}^{*}=U_{n+1}^{}U_{n+1}^{*}$, it suffices to 
show that $(W^{*})^{2}Z_{i,r+1}^{}U_{n+1}=0$. Now
\begin{IEEEeqnarray*}{rCl}
\IEEEeqnarraymulticol{3}{l}{
Z_{i,r+1}^{}U_{n+1} }\\
   &=&  
   Z_{i,r+1}^{}\left(\sum_{\substack{k_{j}\in\mathbb{N}\\(r+2\leq j\leq n+1)}}^{}
       \hspace{-1em} 
       Z_{n+1,n+1}^{k_{n+1}}\ldots Z_{n+1,r+2}^{k_{r+2}}
   Z_{n+1,r}^{}(Z_{n+1,r+2}^{k_{r+2}})^{*}\ldots (Z_{n+1,n+1}^{k_{n+1}})^{*}\right)\\
   &=&  \sum_{\substack{k_{j}\in\mathbb{N}\\(r+2\leq j\leq n+1)}}^{}
         \hspace{-1em}
         Z_{n+1,n+1}^{k_{n+1}}\ldots Z_{n+1,r+3}^{k_{r+3}}
          Z_{i,r+1}^{}Z_{n+1,r+2}^{k_{r+2}}Z_{n+1,r}^{}
      (Z_{n+1,r+2}^{k_{r+2}})^{*}\ldots (Z_{n+1,n+1}^{k_{n+1}})^{*}\\
   &=&  \sum_{\substack{k_{j}\in\mathbb{N}\\(r+2\leq j\leq n+1)}}^{}
         \hspace{-1em}
         Z_{n+1,n+1}^{k_{n+1}}\ldots Z_{n+1,r+3}^{k_{r+3}}Z_{n+1,r+2}^{k_{r+2}}
           Z_{i,r+1}^{}Z_{n+1,r}^{}
           (Z_{n+1,r+2}^{k_{r+2}})^{*}\ldots (Z_{n+1,n+1}^{k_{n+1}})^{*}\\
   && \quad {}+ \hspace{-1em}
     \sum_{\substack{k_{j}\in\mathbb{N}\\(r+2\leq j\leq n+1)}}^{}
         \hspace{-1em}
         Z_{n+1,n+1}^{k_{n+1}}\ldots  Z_{n+1,r+2}^{k_{r+2}}
           Z_{n+1,r+1}^{}Z_{i,r+2}^{}
      (Z_{n+1,r+2}^{1+k_{r+2}})^{*}(Z_{n+1,r+3}^{k_{r+3}})^{*}
                \ldots (Z_{n+1,n+1}^{k_{n+1}})^{*}.
\end{IEEEeqnarray*}
Hence
{\allowdisplaybreaks
\begin{IEEEeqnarray*}{rCl}
\IEEEeqnarraymulticol{3}{l}{
(W^{*})^{2}Z_{i,r+1}^{}U_{n+1} }\\
   &=&  \left(\sum_{\substack{k_{j}\in\mathbb{N}\\(r+2\leq j\leq n+1)}}^{}
         \hspace{-1em}
         Z_{n+1,n+1}^{k_{n+1}}\ldots Z_{n+1,r+2}^{k_{r+2}}
   (Z_{n+1,r+1}^{*})^{2}
   (Z_{n+1,r+2}^{k_{r+2}})^{*}\ldots (Z_{n+1,n+1}^{k_{n+1}})^{*}\right)\times\\
   && \left( \sum_{\substack{k_{j}\in\mathbb{N}\\(r+2\leq j\leq n+1)}}^{}
         \hspace{-1em}
         Z_{n+1,n+1}^{k_{n+1}}\ldots Z_{n+1,r+3}^{k_{r+3}}Z_{n+1,r+2}^{k_{r+2}}
           Z_{i,r+1}^{}Z_{n+1,r}^{}
           (Z_{n+1,r+2}^{k_{r+2}})^{*}\ldots (Z_{n+1,n+1}^{k_{n+1}})^{*}\right.\\
   &&{}+ \left.\sum_{\substack{k_{j}\in\mathbb{N}\\(r+2\leq j\leq n+1)}}^{}
         \hspace{-1em}
         Z_{n+1,n+1}^{k_{n+1}}\ldots  Z_{n+1,r+2}^{k_{r+2}}
           Z_{n+1,r+1}^{}Z_{i,r+2}^{}
      (Z_{n+1,r+2}^{1+k_{r+2}})^{*}(Z_{n+1,r+3}^{k_{r+3}})^{*}
                \ldots (Z_{n+1,n+1}^{k_{n+1}})^{*}\right)\\
   &=&  \hspace{-1em}
   \sum_{\substack{k_{j}\in\mathbb{N}\\(r+2\leq j\leq n+1)}}^{}
         \hspace{-1em}
         Z_{n+1,n+1}^{k_{n+1}}\ldots Z_{n+1,r+3}^{k_{r+3}}Z_{n+1,r+2}^{k_{r+2}}
           (Z_{n+1,r+1}^{*})^{2}Z_{i,r+1}^{}Z_{n+1,r}^{}
           (Z_{n+1,r+2}^{k_{r+2}})^{*}\ldots (Z_{n+1,n+1}^{k_{n+1}})^{*}\\
   &&{} + \hspace{-1em}\sum_{\substack{k,k_{j}\in\mathbb{N}\\(r+2\leq j\leq n+1)}}^{}
         \hspace{-1em}
         \Biggl(
         Z_{n+1,n+1}^{k_{n+1}}\ldots  Z_{n+1,r+2}^{k_{r+2}}
      (Z_{n+1,r+1}^{*})^{2} (Z_{n+1,r+2}^{k_{r+2}})^{*}(Z_{n+1,r+2}^{})^{k+k_{r+2}}
          Z_{n+1,r+1}^{}Z_{i,r+2}^{}\\
     &&\hspace{180pt}
      (Z_{n+1,r+2}^{k+1+k_{r+2}})^{*}(Z_{n+1,r+3}^{k_{r+3}})^{*}
                \ldots (Z_{n+1,n+1}^{k_{n+1}})^{*}\Biggr)\\
   &=&  \hspace{-1em}
       \sum_{\substack{k,k_{j}\in\mathbb{N}\\(r+2\leq j\leq n+1)}}^{}
         \hspace{-1em}
         \Biggl(
         Z_{n+1,n+1}^{k_{n+1}}\ldots  Z_{n+1,r+2}^{k_{r+2}}
      (Z_{n+1,r+1}^{*})^{2} (Z_{n+1,r+2}^{})^{k}
          Z_{n+1,r+1}^{}Z_{i,r+2}^{}\\
     &&\hspace{180pt} 
          (Z_{n+1,r+2}^{k+1+k_{r+2}})^{*}(Z_{n+1,r+3}^{k_{r+3}})^{*}
                \ldots (Z_{n+1,n+1}^{k_{n+1}})^{*}\Biggr)\\
    &=& 0.
\end{IEEEeqnarray*}
}
Thus the proof is complete.
\end{proof}
%%------------------------------------------------------

For the remaining operators, i.e.\ $W_{}^{*}Z_{i,r+1}^{}$ for $i\geq r+2$
and $Z_{i,r}^{}W$ for $i\geq r+1$, 
we now invoke the relation
(\ref{eq:an0-defrel-9}). Observe that for $i\geq r+2$, one has
\begin{IEEEeqnarray*}{rCl}
\IEEEeqnarraymulticol{3}{l}{
  \left(W_{n+1,r+1}^{*}Z_{i,r+1}^{}\right)^{*}}\\
  \quad &=& \left(Z_{1,1}^{}Z_{2,2}^{}\ldots Z_{r,r}^{}\right)
         \left(Z_{r+1,r+2}^{}\ldots Z_{i-1,i}^{}\right)
           \left(Z_{i+1,i+1}^{}\ldots Z_{n+1,n+1}^{}\right)W_{n+1,r+1}^{}\\
  &=& \left(Z_{1,1}^{}Z_{2,2}^{}\ldots Z_{r-1,r-1}^{}\right)
     \left(Z_{r,r}^{}W_{n+1,r+1}^{}\right)
         \left(Z_{r+1,r+2}^{}\ldots Z_{i-1,i}^{}\right)
           \left(Z_{i+1,i+1}^{}\ldots Z_{n+1,n+1}^{}\right),
\end{IEEEeqnarray*}
and for $i\geq r+1$, one similarly has
\begin{IEEEeqnarray*}{rCl}
\IEEEeqnarraymulticol{3}{l}{
  \left(Z_{i,r}^{}W_{n+1,r+1}^{}\right)^{*}}\\
  \quad &=& W_{n+1,r+1}^{*}\left(Z_{1,1}^{}Z_{2,2}^{}\ldots Z_{r-1,r-1}^{}\right)
         \left(Z_{r,r+1}^{}\ldots Z_{i-1,i}^{}\right)
           \left(Z_{i+1,i+1}^{}\ldots Z_{n+1,n+1}^{}\right)\\
  &=& \left(Z_{1,1}^{}Z_{2,2}^{}\ldots Z_{r-1,r-1}^{}\right)
         W_{n+1,r+1}^{*}Z_{r,r+1}^{}\left(Z_{r+1,r+2}^{}\ldots Z_{i-1,i}^{}\right)
           \left(Z_{i+1,i+1}^{}\ldots Z_{n+1,n+1}^{}\right).
\end{IEEEeqnarray*}
Therefore using the earlier results (Propositions \ref{pr:wcomm-1a}, 
\ref{pr:wcomm-7} and  \ref{pr:wcomm-2}), we now obtain the following proposition.
%%------------------------------------------------------
\bppsn %[\texttt{pr:wcomm-8}]
     %%-------------------------------
     \label{pr:wcomm-8}
     %%-------------------------------
 Let $W$ be as in Theorem~\ref{th:factorization-2}. Then
\begin{IEEEeqnarray*}{rCll}
    W_{}^{*}Z_{i,r+1}^{} &\in & \{W_{}^{},W_{}^{*}\}^{\prime}& \qquad \text{for } i\geq r+2. \\
    Z_{i,r}^{}W  &\in  &\{W,W^{*}\}^{\prime} & \qquad \text{for } i\geq r+1.
\end{IEEEeqnarray*}
\eppsn
%%------------------------------------------------------

\noindent
\textbf{Proof of Theorem~\ref{th:factorization-2}:}\\
%%------------------------------------------------------
From Propositions 
\ref{pr:wcomm-1a}, \ref{pr:wcomm-2} and  
\ref{pr:wcomm-7}--\ref{pr:wcomm-8}, 
it follows that
$Y_{i,j}^{(1)}\in \{W,W^{*}\}^{\prime}$ for all $i,j$.
So it remains to prove (\ref{eq:decomp}). For that,
we need to prove the following equalities
\begin{IEEEeqnarray}{rCll}
  Z_{i,r}^{}WW^{*} &=& Z_{i,r}^{}& \qquad  i\leq r,
     %%-------------------------------
     \label{eq:decomp-1}\\
     %%-------------------------------
  Z_{i,r}^{}W (I-WW^{*})+ W^{*}Z_{i,r+1}^{}W^{} &=& Z_{i,r+1}^{}
              & \qquad  i\leq r,
     %%-------------------------------
     \label{eq:decomp-2}\\
     %%-------------------------------
  W^{*}Z_{i,r+1}^{}W &=& Z_{i,r+1}^{}&  \qquad i\geq r+1,
     %%-------------------------------
     \label{eq:decomp-3}\\
     %%-------------------------------
  W^{*}Z_{i,r+1}^{} (I-WW^{*})+ Z_{i,r}^{}WW^{*}&=& Z_{i,r}^{}&  \qquad i\geq r+1.
     %%-------------------------------
     \label{eq:decomp-4}
     %%-------------------------------
\end{IEEEeqnarray}
For (\ref{eq:decomp-1}), note that
\begin{IEEEeqnarray*}{rCl}
Z_{i,r}^{}(I-WW^{*}) 
    &=&  Z_{i,r}^{}(I-W_{n+1,r+1}^{}W_{n+1,r+1}^{*})\\
    &=&  Z_{i,r}^{}U_{n+1}^{}U_{n+1}^{*}.
\end{IEEEeqnarray*}
Since $i\leq r$, we have
\begin{IEEEeqnarray*}{rCl}
\IEEEeqnarraymulticol{3}{l}{
Z_{i,r}^{}U_{n+1}^{} }\\
    &=& \!\!\!\sum_{\substack{k_{j}\in\mathbb{N}\\(r+2\leq j\leq n+1)}}^{}
        \hspace{-1em}
        Z_{n+1,n+1}^{k_{n+1}}\ldots Z_{n+1,r+3}^{k_{r+3}}Z_{n+1,r+2}^{k_{r+2}}
           \left(Z_{i,r}^{}Z_{n+1,r}^{}\right)
           (Z_{n+1,r+2}^{k_{r+2}})^{*}\ldots (Z_{n+1,n+1}^{k_{n+1}})^{*}\\
    &=& 0.
\end{IEEEeqnarray*}
Therefore
\begin{IEEEeqnarray*}{rCl}
Z_{i,r}^{}WW^{*}
    &=& Z_{i,r}^{}- Z_{i,r}^{}(I-WW^{*})\\
    &=& Z_{i,r}^{}- Z_{i,r}^{}U_{n+1}^{}U_{n+1}^{*}\\
    &=& 0.
\end{IEEEeqnarray*}
Next let us prove (\ref{eq:decomp-2}). Use Proposition~\ref{pr:rep-rel-6}
to get
\begin{IEEEeqnarray*}{rCl}
Z_{i,r}^{}W(I-WW^{*})
    &=& W_{n+1,r}^{*}Z_{i,r}^{}W_{n+1,r+1}^{}(I-WW^{*})\\
    &=& W_{n+1,r}^{*}
       \left(W_{n+1,r+1}^{}Z_{i,r}^{}+U_{n+1}^{}Z_{i,r+1}^{}\right)(I-WW^{*})\\
    &=& WZ_{i,r}^{}(I-WW^{*}) + W_{n+1,r}^{*} U_{n+1}^{}Z_{i,r+1}^{}(I-WW^{*})\\
    &=& WZ_{i,r}^{}U_{n+1}^{}U_{n+1}^{*} + U_{n+1}^{}U_{n+1}^{*}Z_{i,r+1}^{}(I-WW^{*})\\
    &=&  (I-WW^{*})Z_{i,r+1}^{}(I-WW^{*}).
\end{IEEEeqnarray*}
Since $W^{*}Z_{i,r+1}^{}$ commutes with $W$, we get
\begin{IEEEeqnarray*}{rCl}
Z_{i,r}^{}W (I-WW^{*})+ W^{*}Z_{i,r+1}^{}W^{}
  &=& (I-WW^{*})Z_{i,r+1}^{}(I-WW^{*}) + W^{}W^{*}Z_{i,r+1}^{}\\
  &=& Z_{i,r+1}^{}-(I-WW^{*})Z_{i,r+1}^{}WW^{*}\\
  &=& Z_{i,r+1}^{}-U_{n+1}^{*}U_{n+1}^{}Z_{i,r+1}^{}WW^{*}.
\end{IEEEeqnarray*}
By (\ref{eq:wcomm-7pr}), the second term in the right hand side above vanishes.
Therefore we have (\ref{eq:decomp-2}).

For proving (\ref{eq:decomp-3}) and (\ref{eq:decomp-4}), we will need the following equality.
%%------------------------------------------------------
\begin{IEEEeqnarray}{rCl}
W^{*}Z_{i,r+1}^{}-Z_{i,r+1}^{}W^{*} &=& Z_{i,r}^{}(I-WW^{*})\qquad\text{if }i\geq r+1.
%      %%-------------------------------
     \label{eq:wcomm-11}
%      %%-------------------------------
\end{IEEEeqnarray}
%%------------------------------------------------------
Here
\begin{IEEEeqnarray*}{rCl}
\IEEEeqnarraymulticol{3}{l}{
  Z_{i,r+1}^{*}W_{n+1,r+1}^{}-W_{n+1,r+1}^{}Z_{i,r+1}^{*}}\\
  &=& \left(Z_{1,1}^{}Z_{2,2}^{}\ldots Z_{r,r}^{}\right)
         \left(Z_{r+1,r+2}^{}\ldots Z_{i-1,i}^{}\right)
           \left(Z_{i+1,i+1}^{}\ldots Z_{n+1,n+1}^{}\right)W_{n+1,r+1}^{}\\
   && \quad {} - 
        W_{n+1,r+1}^{}\left(Z_{1,1}^{}Z_{2,2}^{}\ldots Z_{r,r}^{}\right)
         \left(Z_{r+1,r+2}^{}\ldots Z_{i-1,i}^{}\right)
           \left(Z_{i+1,i+1}^{}\ldots Z_{n+1,n+1}^{}\right)\\
  &=& \left(Z_{1,1}^{}Z_{2,2}^{}\ldots Z_{r-1,r-1}^{}\right)
     \left(Z_{r,r}^{}W_{n+1,r+1}^{}-W_{n+1,r+1}^{}Z_{r,r}^{}\right)\\
     &&\qquad\qquad\qquad
         \left(Z_{r+1,r+2}^{}\ldots Z_{i-1,i}^{}\right)
           \left(Z_{i+1,i+1}^{}\ldots Z_{n+1,n+1}^{}\right)\\
  &=& \left(Z_{1,1}^{}Z_{2,2}^{}\ldots Z_{r-1,r-1}^{}\right)
     \left(U_{n+1}Z_{r,r+1}^{}\right)
         \left(Z_{r+1,r+2}^{}\ldots Z_{i-1,i}^{}\right)
           \left(Z_{i+1,i+1}^{}\ldots Z_{n+1,n+1}^{}\right)\\
  &=& U_{n+1}\left(Z_{1,1}^{}Z_{2,2}^{}\ldots Z_{r-1,r-1}^{}\right)
     \left(Z_{r,r+1}^{}Z_{r+1,r+2}^{}\ldots Z_{i-1,i}^{}\right)
           \left(Z_{i+1,i+1}^{}\ldots Z_{n+1,n+1}^{}\right)\\
  &=& U_{n+1}Z_{i,r}^{*}.
\end{IEEEeqnarray*}
Therefore we have
\begin{IEEEeqnarray*}{rCl}
W^{*}Z_{i,r+1}^{}-Z_{i,r+1}^{}W^{*} &=& Z_{i,r}^{}U_{n+1}^{*}W_{n+1,r}^{}\\
    &=& Z_{i,r}^{}(I-WW^{*}).
\end{IEEEeqnarray*}

Using the above equality, we get
\begin{IEEEeqnarray*}{rCl}
W^{*}Z_{i,r+1}^{}W - Z_{i,r+1}^{}
  &=& \left(W^{*}Z_{i,r+1}^{}-Z_{i,r+1}^{}W^{*}\right)W\\
  &=& Z_{i,r}^{}(I-WW^{*})W\\
  &=& 0.
\end{IEEEeqnarray*}
Thus we have (\ref{eq:decomp-3}). Next,
\begin{IEEEeqnarray*}{rCl}
\IEEEeqnarraymulticol{3}{l}{
W^{*}Z_{i,r+1}^{} (I-WW^{*})+ Z_{i,r}^{}WW^{*}}\\
   \qquad\qquad &=& \left(Z_{i,r+1}^{}W^{*} + Z_{i,r}^{}(I-WW^{*})\right)(I-WW^{*})+ Z_{i,r}^{}WW^{*}\\
    &=& Z_{i,r}^{}.
\end{IEEEeqnarray*}
Thus we have proved (\ref{eq:decomp-4}).
This completes the proof of Theorem~\ref{th:factorization-2}.
%%------------------------------------------------------

%%---------------------------------------------------
\section{Classification theorem}
%%---------------------------------------------------
We are now in a position to prove a factorization theorem 
for  irreducible representations of the $C^{*}$-algebra $C(SU_{0}(n+1))$
that will subsequently lead to a recursive proof of the classification theorem
for all irreducible representations of $C(SU_{0}(n+1))$.
%%------------------------------------------------------
\bthm[Factorization theorem] % [\texttt{th:factorization-1}]
     %%-------------------------------
     \label{th:factorization-1}
     %%-------------------------------
Let $\pi$ be an irreducible representation of $C(SU_{0}(n+1))$ on a Hilbert space
$\mathcal{H}$. Assume $r\equiv r(\pi)\leq n$ where $r$ is as in (\ref{eq:rep-rel-3}).
Then there is an irreducible representation $\pi_{1}$ of $C(SU_{0}(n+1))$ acting on a 
Hilbert space $\mathcal{H}_{1}$ such that $r(\pi_{1})= r+1$
and $\pi$ is unitarily equivalent to the representation
$\pi_{1}\ast\psi_{s_{r}}$ acting on $\mathcal{H}_{1}\otimes\ell^{2}(\mathbb{N})$
given by
\begin{IEEEeqnarray*}{rClCl}
 \pi_{1}\ast\psi_{s_{r}}(z_{i,j}^{}) & =&
     \sum_{k=\min\{i,j\}}^{\max\{i,j\}}
         {\pi}_{1}(z_{i,k}^{})\otimes \psi_{s_{r}}(z_{k,j}^{}),
     \quad 1\leq i,j\leq n+1.
\end{IEEEeqnarray*}
%%------------------------------------------------------
\ethm
%%--------------------------------------------------------------------
\begin{proof}
From Proposition~\ref{pr:wcomm-3c}, we know that $W$ is an  
isometry. Since $U_{n+1}^{*}U_{n+1}^{}\geq P_{n+1,r}^{}$, 
it follows from (\ref{eq:rep-rel-13})  that $W$ is not a unitary.
By Wold decomposition, there are Hilbert spaces 
$\mathcal{H}_{0}$ and $\mathcal{H}_{1}$ and a unitary $U$ 
on $\mathcal{H}_{0}$ such that 
$(\mathcal{H},W)$ is unitarily equivalent to 
$((\mathcal{H}_{1}\otimes\ell^{2}(\mathbb{N}))\oplus 
 \mathcal{H}_{0}, (I\otimes S^{*})\oplus U)$.
 By Theorem~\ref{th:factorization-2},
the subspace $\mathcal{H}_{1}\otimes\ell^{2}(\mathbb{N})$
is kept invariant by  $Y_{i,j}^{(1)}$, $Y_{i,j}^{(2)}$, $(Y_{i,j}^{(1)})^{*}$ and
$(Y_{i,j}^{(2)})^{*}$ for all $i,j$.
Therefore the subspace $\mathcal{H}_{1}\otimes\ell^{2}(\mathbb{N})$ 
is an invariant subspace for $\pi$. 
By irreducibility of $\pi$, it follows that $\mathcal{H}_{0}=\{0\}$.

%%--------------------------------------------------------------------

Thus there are operators
$Z_{i,j}^{(1)}\in\mathcal{L}(\mathcal{H}_{1})$ and 
$Z_{i,j}^{(2)}\in\mathcal{L}(\ell^{2}(\mathbb{N}))$ such that
\[
Y_{i,j}^{(1)}=Z_{i,j}^{(1)}\otimes I,\qquad
Y_{i,j}^{(2)}=I\otimes Z_{i,j}^{(2)},\qquad 1\leq i,j\leq n+1,
\]
and we have
\[
Z_{i,j}^{}=\sum_{k=\min\{i,j\}}^{\max\{i,j\}}
         Z_{i,k}^{(1)}\otimes Z_{k,j}^{(2)}
   \qquad 1\leq i,j\leq n+1.
\]
In fact one has
\begin{IEEEeqnarray}{rCl}
  Z_{i,j}^{(2)} &=& 
    \begin{cases}
     S^{*} & \text{if }i=j=r+1,\\
     S & \text{if }i=j=r,\\
     P_{0} & \text{if $i,j\in \{r,r+1\}$ and $i\neq j$},\\
     I & \text{if }i=j\not\in\{r,r+1\},\\
     0 & \text{otherwise.}
    \end{cases}
\end{IEEEeqnarray}
Thus  the map
$z_{i,j}^{}\mapsto Z_{i,j}^{(2)}$
defines a representation of $C(SU_{0}(n+1))$ on $\ell^{2}(\mathbb{N})$
which is in fact the representation $\psi_{s_{r}}$ defined in section~3.
Observe that the $Z_{i,j}^{(2)}$'s  belong to the Toeplitz algebra
$\mathscr{T}\subseteq \mathcal{L}(\ell^{2}(\mathbb{N}))$.
Therefore
\begin{equation}
\pi(z_{i,j}^{})\in \mathcal{L}(\mathcal{H}_{1})\otimes \mathscr{T} 
   \quad\text{for all } i,j.
\end{equation}
Let $\sigma:\mathscr{T}\to\mathbb{C}$ be the $*$-homomorphism given by
$S\mapsto 1$. Then one has
\[
Z_{i,j}^{(1)}=(id\otimes\sigma)Z_{i,j}^{}=(id\otimes\sigma)\pi(z_{i,j}^{}).
\]
Thus the map $\pi_{1}:z_{i,j}^{}\mapsto Z_{i,j}^{(1)}$ gives a representation of 
$C(SU_{0}(n+1))$ on $\mathcal{H}_{0}$ and we have
\[
 \pi(z_{i,j}^{}) =\sum_{k=\min\{i,j\}}^{\max\{i,j\}}
   {\pi}_{1}(z_{i,k}^{})\otimes \psi_{s_{r}}(z_{k,j}^{}),\quad 1\leq i,j\leq n+1.
\]
If $T$ belongs to the commutant of ${\pi}_{1}(C(SU_{0}(n+1)))$, then 
$T\otimes I \in \pi(C(SU_{0}(n+1)))^\prime$. By irreducibility of $\pi$, 
it follows that $T=I$. Thus ${\pi}_{1}$ is irreducible. 
Using (\ref{eq:an0-defrel-1}), we get $\pi_{1}(z_{n+1,r}^{})=Z_{n+1,r}^{}W=0$.
Also, note that
\begin{IEEEeqnarray*}{rCl}
  \left(W^{*}Z_{n+1,r+1}^{}\right)^{*} \left(W^{*}Z_{n+1,r+1}^{}\right)
    &=& Z_{n+1,r+1}^{*}WW^{*}Z_{n+1,r+1}^{}\\
    &=& P_{n+1,r+1}- Z_{n+1,r+1}^{*}(I-WW^{*})Z_{n+1,r+1}^{}\\
    &=& P_{n+1,r+1}- Z_{n+1,r+1}^{*}U_{n+1}^{*}U_{n+1}^{}Z_{n+1,r+1}^{}\\
    &=& P_{n+1,r+1} \neq 0.
\end{IEEEeqnarray*}
Thus we have 
\[
\min\{i: 1\leq i\leq n+1, \pi_{1}(z_{n+1,i}^{})\neq 0\} = r+1.
\]
This completes the proof.
\end{proof}

We next come to the main result that gives a parametrization of all 
the irreducible representations.
%%--------------------------------------------------------------------
\bthm %[\texttt{th:main-1}]
     %%-------------------------------
     \label{th:main-1}
     %%-------------------------------
Let $\pi$ be an irreducible representation of $C(SU_{0}(n+1))$ on a 
Hilbert space $\mathcal{H}$.
Then there is a $\lambda=(\lambda_{1},\ldots,\lambda_{n})\in (S^{1})^{n}$ 
and a reduced word $\omega$ in $\mathfrak{S}_{n+1}$ 
of the form
$s_{[a_{k},b_{k}]}s_{[a_{k-1},b_{k-1}]}\ldots  s_{[a_{1},b_{1}]}$ 
where
$1\leq b_{k}<b_{k-1}<\cdots <b_{1}\leq n$, and 
$1\leq a_{i}\leq b_{i}$ for $1\leq i\leq k$,
such that 
$\pi\cong \psi_{\lambda,\omega}$.
\ethm
%%--------------------------------------------------------------------
%%--------------------------------------------------------------------
\begin{proof}
%%--------------------------------------------------------------------
Note that when $\pi\cong \psi_{\lambda,\omega}$, then
$1\leq r\leq n$ implies $b_{1}=n$ and $a_{1}=r$.

It was proved in \cite{GirPal-2022tv} that the statement holds for $n=2$.
%%--------------------------------------------------------------------
Let us assume that it holds for $n-1$.
Let us first deal with the case $r=n+1$.
%%--------------------------------------------------------------------
% \subsection{Step 1: $r=n+1$}
%%--------------------------------------------------------------------
In this case,  from Proposition~\ref{pr:rep-rel-1} it 
follows that $Z_{n+1,j}^{}=0$ for $1\leq j\leq n$. Therefore we have
\begin{enumerate}
  \item
  $P_{n+1,n+1}^{}=Q_{n+1,n+1}^{}=I$,
  \item
  $Z_{i,n+1}^{}=0$ for $1\leq i\leq n$,
  \item
  $Z_{i,j}^{}Z_{n+1,n+1}=Z_{n+1,n+1}Z_{i,j}^{}$ for all $i,j$.
\end{enumerate}
Therefore $Z_{n+1,n+1}^{}$ is a unitary and any spectral subspace 
of $Z_{n+1,n+1}^{}$ is an invariant subspace for $\pi$. 
By irreducibility of $\pi$, it follows that 
$Z_{n+1,n+1}^{}=\mu_{0} I$ for some $\mu_{0}\in S^{1}$. Now define
\[
\pi_{0}(z_{i,j}^{(n-1)})= 
        \begin{cases}
        Z_{i,j}^{} & \text{if }  1\leq j\leq n\text{ and }1\leq i\leq n-1 ,\\
        \mu_{0}Z_{i,j}^{} & \text{if } 1\leq j\leq n \text{ and } i=n,
        \end{cases}
\]
where $z_{i,j}^{(n-1)}$'s denote the generators for $C(SU_{0}(n))$.
Then $\pi_{0}$ is an irreducible representation of $C(SU_{0}(n))$ 
on $\mathcal{H}$. Therefore
there is a $\mu=(\mu_{1},\ldots,\mu_{n-1})\in (S^{1})^{n-1}$ and 
a reduced word 
$\omega=s_{[a_{k},b_{k}]}\ldots s_{[a_{1},b_{1}]}\in
     \mathfrak{S}_{n}\subseteq\mathfrak{S}_{n+1}$  
where $1\leq b_{k}<b_{k-1}<\cdots<b_{1}\leq n-1$ 
such that $\pi_{0}\cong \psi_{\mu,\omega}$.
We thus have
\[
\pi\cong \psi_{\lambda,\omega},
\]
where $\lambda=(\mu_{1},\mu_{2},\ldots,\mu_{n-1},\overline{\mu_{0}})$.

%%--------------------------------------------------------------------
%%--------------------------------------------------------------------
Next, assume $1\leq r\leq n$. By an application of the second factorization
result proved earlier (Theorem~\ref{th:factorization-1}), it follows that
there is an irreducible representation $\pi_{1}$ acting on a Hilbert space
$\mathcal{H}_{1}$ such that $r(\pi_{1})= r+1$ and 
$\pi={\pi}_{1}\ast \psi_{s_{r}}$ i.e.\
\begin{IEEEeqnarray*}{rCl}
\pi(z_{i,j}^{}) &=&
    \sum_{k=\min\{i,j\}}^{\max\{i,j\}}
         {\pi}_{1}(z_{i,k}^{})\otimes \psi_{s_{r}}(z_{k,j}^{}).
\end{IEEEeqnarray*}
Using Lemma~\ref{lm:irred-q0-0} and by repeated application of Theorem~\ref{th:factorization-1},
we conclude that there is an irreducible representation $\pi_{0}$ acting on a Hilbert space
$\mathcal{H}_{0}$ such that $r(\pi_{0})= n+1$ and $\pi=\pi_{0}\ast \psi_{s_{[r,n]}}$.
By the previous case, the result now follows.
\end{proof}
%%---------------------------------------------
%%---------------------------------------------
\brmrk\label{rmrk:one-dim}
Note that one gets the one dimensional representations $\psi_{\lambda,id}^{}$
given by (\ref{eq:irred-q0-2aa}) when $r=m+1$ in each of the steps $m=n, n-1,\ldots, 2$ in the above recursive procedure.
\ermrk
%%---------------------------------------------

%%---------------------------------------------------
\section{Corollaries of the classification theorem}
%%---------------------------------------------------
In this section, we present a few important consequences of the 
classification theorem.
%%---------------------------------------------
\bppsn
 %%---------------------------
 \label{prop:all-irr}
 %%---------------------------
 The set
\[
\left\{\psi_{\lambda,\omega}:\lambda\in (S^{1})^{n},\,
     \omega \text{ a reduced word in }\mathfrak{S}_{n+1}\right\}
\]
gives all inequivalent irreducible representations of $C(SU_{0}(n+1))$.
\eppsn
%%---------------------------------------------
This follows from Theorems~\ref{th:irred-q0-3} and \ref{th:irred-q0-4} in Section~3 
and Theorem~\ref{th:main-1}.

Recall that $\mathcal{O}(SU_{0}(n+1))$ denotes the *-algebra given 
by the relations (\ref{eq:an0-defrel-1})--\ref{eq:an0-defrel-8}) in 
Theorem~\ref{th:an0-defrel}.
Thus $\mathcal{O}(SU_{0}(n+1))$ can be viewed as a dense *-subalgebra
of $C(SU_{0}(n+1))$, and the family $\psi_{\lambda,\omega}$ gives all the 
irreducible *-representations of this algebra.

%%---------------------------------------------
\bppsn
 %%---------------------------
 \label{prop:type1}
 %%---------------------------
$C(SU_{0}(n+1))$ is a $C^{*}$-algebra of type I.
\eppsn
%%---------------------------------------------
\begin{proof}
For each $1\leq r\leq k$, let
$i_{s}$ and $j_{s}$ be nonnegative integers for $a_{r}\leq s\leq b_{r}$.
Let $|e_{i}\rangle\,\langle e_{j}|$ be the rank one operator
$x\mapsto \langle e_{j},x\rangle e_{i}$. Denote by $T_{r}$ the operator
\[
|e_{i_{b_{r}}}\rangle\,\langle e_{j_{b_{r}}}|\otimes
        |e_{i_{b_{r}-1}}\rangle\,\langle e_{j_{b_{r}-1}}|\otimes \cdots
          \otimes |e_{i_{a_{r}}}\rangle\,\langle e_{j_{a_{r}}}|.
\]
It is enough to prove that 
$T_{k}\otimes T_{k-1}\otimes\cdots \otimes T_{1}$ belongs to the 
image $\psi_{\lambda,\omega}(C(SU_{0}(n+1)))$.

We will use the operators
\[
E_{j,i}^{}:= V_{n'(j,i),i}^{*}V_{j,i}^{} 
    |V_{j-1,a_{j-1}}^{}|\cdot|V_{j-2,a_{j-2}}^{}|\cdots|V_{1,a_{1}}^{}| 
\]
described in Corollary~\ref{cr:irred-q0-1c} for this. Recall that
\[
E_{j,i}^{}=I^{\otimes \sum_{s=j+1}^{k}(b_{s}-a_{s}+1)}\otimes 
    \left(P_{0}^{\otimes {(b_{j}+1-i)}} \otimes S^{*}\otimes I^{\otimes(i-a_{j}-1)} \right)
     \otimes P_{0}^{\otimes \sum_{s=1}^{j-1}(b_{s}-a_{s}+1)}.
\]
Now observe that $|e_{i}\rangle\,\langle e_{j}|=(S^{i})^{*}P_{0}S^{j}$ and hence 
\begin{IEEEeqnarray*}{rCl}
  \IEEEeqnarraymulticol{3}{l}{
  E_{k,1+b_{k}}^{i_{b_{k}}}\ldots E_{k,1+a_{k}}^{i_{a_{k}}}
        E_{k,a_{k}}^{}
        (E_{k,1+a_{k}}^{j_{a_{k}}})^{*}\ldots (E_{k,1+b_{k}}^{j_{b_{k}}})^{*}}\\
     \qquad\qquad  \qquad\qquad  &=& 
       T_{k}\otimes P_{0}^{\otimes(b_{k-1}-a_{k-1}+1)}\otimes\cdots
               \otimes P_{0}^{\otimes(b_{1}-a_{1}+1)}.
\end{IEEEeqnarray*}
Thus $T_{k}\otimes P_{0}^{\otimes(b_{k-1}-a_{k-1}+1)}\otimes\cdots
               \otimes P_{0}^{\otimes(b_{1}-a_{1}+1)}$
belongs to $\psi_{\lambda,\omega}(C(SU_{0}(n+1)))$.
Having proved that
\[
T:=T_{k}\otimes T_{k-1}\otimes\cdots\otimes T_{s+1}\otimes
    P_{0}^{\otimes(b_{s}-a_{s}+1)}\otimes\cdots
     \otimes P_{0}^{\otimes(b_{1}-a_{1}+1)}
\]
belongs to $\psi_{\lambda,\omega}(C(SU_{0}(n+1)))$, note that
\begin{IEEEeqnarray*}{rCl}
  \IEEEeqnarraymulticol{3}{l}{
  E_{s,1+b_{s}}^{i_{b_{s}}}\ldots E_{s,1+a_{s}}^{i_{a_{s}}}
        T
        (E_{s,1+a_{s}}^{j_{a_{s}}})^{*}\ldots (E_{s,1+b_{s}}^{j_{b_{s}}})^{*}}\\
        \qquad\qquad &=& T_{k}\otimes T_{k-1}\otimes\cdots\otimes T_{s}\otimes
    P_{0}^{\otimes(b_{s-1}-a_{s-1}+1)}\otimes\cdots
     \otimes P_{0}^{\otimes(b_{1}-a_{1}+1)}.
\end{IEEEeqnarray*}
Thus by repeating the arguent, the result follows.
\end{proof}
%%---------------------------------------------

%%---------------------------------------------
\bppsn
 %%---------------------------
 % \label{prop:CQS}
 %%---------------------------
 \begin{enumerate}
   \item
The map
\[
\Delta(z_{i,j}^{})=\sum_{k=\min\{i,j\}}^{\max\{i,j\}}z_{i,k}^{}\otimes z_{k,j}^{}
\]
extends to a unital $C^{*}$-homomorphism from $C(SU_{0}(n+1))$
to the tensor product $C(SU_{0}(n+1)) \otimes C(SU_{0}(n+1))$.
   \item
   The $C^{*}$-algebra $C(SU_{0}(n+1))$ together with $\Delta$ and
the one dimensional representation $\epsilon=\psi_{id}$ is
a $C^{*}$-bialgebra, i.e.\ a compact quantum semigroup.
\end{enumerate}
\eppsn
%%---------------------------------------------
Proof is identical to the arguments used for the $n=2$ case 
in \cite{GirPal-2022tv}. Observe also that the restriction of 
$\Delta$ and $\epsilon$ to $\mathcal{O}(SU_{0}(n+1))$ turn it 
into a *-bialgebra.
Using this map $\Delta$, one can write the representations
$\psi_{\lambda,\omega}$ as convolution products of $\mychi_{\lambda}$ 
and $\psi_{s_{i}}$'s, just as in the $q\neq 0$ case.
%%---------------------------------------------

    Let $1\leq m < n$. Denote the generators of $C(SU_{0}(m+1))$ and 
    $C(SU_{0}(n+1))$ by $z^{(m)}_{i,j}$ and $z^{(n)}_{i,j}$ respectively, 
    and their comultiplication maps by $\Delta_{m}$ and $\Delta_{n}$ respectively. 
    Then it follows from the defining
    relations (\ref{eq:an0-defrel-1}--\ref{eq:an0-defrel-9}) that the map
    $\phi$ from  $C(SU_{0}(n+1))$ to $C(SU_{0}(m+1))$ given by
\begin{IEEEeqnarray*}{rCl}
  \phi(z^{(n)}_{i,j}) &=& 
       \begin{cases}
       z^{(m)}_{i,j} & \text{if }1\leq i,j\leq m+1,\\
       \delta_{i,j} & \text{otherwise}
       \end{cases}
\end{IEEEeqnarray*}
    is a surjective $C^*$-homomorphism and one has 
    $(\phi\otimes\phi)\Delta_{n}=\Delta_{m}\phi$. 
    Thus $(C(SU_{0}(m+1)), \Delta_{m})$
    is a compact quantum subsemigroup of $(C(SU_{0}(n+1)), \Delta_{n})$. 
    By the same argument used for $C(SU_{0}(3))$ in \cite{GirPal-2022tv}, 
    it follows that 
    $(C(SU_{0}(n+1)), \Delta_{n})$ is not a compact quantum group.

%%---------------------------------------------

%%---------------------------------------------
\bthm
 %%---------------------------
 \label{th:rel-with-M-Y}
 %%---------------------------
The crystallized algebra $C(SU_{0}(n+1))$ is isomorphic to the
$C^{*}$-subalgebra of the space of bounded operators on
$\ell^{2}(\mathbb{Z})^{\otimes n}\otimes 
\ell^{2}(\mathbb{N})^{\otimes\frac{n(n+1)}{2}}$
generated by the limits
\[
\lim_{q\to 0+} \psi^{(q)}\left((-q)^{\min\{i-j,0\}}u_{i,j}^{}(q)\right),
\]
where $\psi^{(q)}$ is the Soibelman representation of $C(SU_{q}(n+1))$.
\ethm
%%---------------------------------------------
\begin{proof}
Similar to the case of nonzero $q$, one can prove that the direct 
integral $\psi^{(0)}$ of the representations $\psi_{\lambda,\omega_{0}}$ 
over $\lambda\in (S^{1})^{n}$,
    where $\omega_{0}$ is the longest word in $\mathfrak{S}_{n+1}$, 
is a faithful representation of $C(SU_{0}(n+1))$ acting on
    $\ell^{2}(\mathbb{Z})^{\otimes n}\otimes 
\ell^{2}(\mathbb{N})^{\otimes\frac{n(n+1)}{2}}$. 
% (We will call this the Soibelman representation of $C(SU_{0}(n+1))$). 
% Let us denote it by $\psi^{(0)}$.
It is simple to check that the limits
$\lim_{q\to 0+}\psi^{(q)}\left((-q)^{\min\{i-j,0\}} u_{i,j}^{}(q)\right)$
exist and one has
\[
\psi^{(0)}(z_{i,j}^{})=
  \lim_{q\to 0+} \psi^{(q)}\left((-q)^{\min\{i-j,0\}}u_{i,j}^{}(q)\right).
\]
Therefore the result follows.
\end{proof}

Finally, we have the following
important property of our crystallized algebra.
%%---------------------------------------------
\bthm
 %%---------------------------
 % \label{th:crystl-property}
 %%---------------------------
Let $\pi^{(0)}$ be an irreducible representation of 
the crystallization $C(SU_{0}(n+1))$ 
(respectively $\mathcal{O}(SU_{0}(n+1))$)
on a Hilbert space $\mathcal{H}$. Then there exist irreducible
representations $\pi^{(q)}$ of $C(SU_{q}(n+1))$  
(respectively $\mathcal{O}(SU_{q}(n+1))$)
on the same Hilbert space $\mathcal{H}$ such that
%%---------------------------
\begin{IEEEeqnarray*}{rCl}
\pi^{(0)}(z_{i,j}^{})=
 \lim_{q\to 0+}
 \pi^{(q)}\left((-q)^{\min\{i-j,0\}}u_{i,j}^{}(q)\right),\qquad
  i,j\in\{1,2,\ldots,n+1\}.
   %%---------------------------
   % \label{eq:}
   %%---------------------------
\end{IEEEeqnarray*}
%%---------------------------
\ethm
%%---------------------------------------------
\begin{proof}
This follows from Theorem~\ref{th:irred-q0-4}, Proposition~\ref{prop:all-irr} 
and Soibelman's result (Theorem 6.2.7, page 121, \cite{KorSoi-1998ab}).
\end{proof}
%%---------------------------------------------

%%---------------------------------------------
\brmrk
 %%---------------------------
 \label{rm:comparison}
 %%---------------------------
Recall (see the discussion after Remark 2.2.2 in \cite{GirPal-2022tv})
that there is a specialization map 
$\theta_{q}:\mathcal{O}_{t}^{\mathbb{Q}[t,t^{-1}]}(SU(n+1))\to 
      \mathcal{O}_{q}(SU(n+1))$
such that the limits in Theorem~\ref{th:rel-with-M-Y} are
$\lim_{q\to 0+}\psi_{}^{(q)}\circ\theta_{q}
         \bigl((-t)^{\min\{i-j,0\}}u_{i,j}(t)\bigr)$.
Thus Theorem~\ref{th:rel-with-M-Y} realizes the crystallization 
$C(SU_{0}(n+1))$ as a $C^{*}$-algebra of bounded operators 
on the Hilbert space that carries the Soibelman representation of the 
family $C(SU_{q}(n+1))$, and is generated by a collection of operators 
that are $q\to 0+$ limits of images of certain matrix entries from 
$\mathcal{O}_{t}^{\mathbb{Q}[t,t^{-1}]}(SU(n+1))$ . 
This is very similar to the description of 
the crystallized $C^{*}$-algebra given by Mattassa \& Yuncken \cite{MatYun-2023aa}.
However, there are crucial differences. In particular, the generating collection of
matrix entries used are different, and more importantly, the coordinate function
algebras from where they come are over different subrings of $\mathbb{Q}(t)$.
We believe that the two crystallized $C^{*}$-algebras are isomorphic, however we do 
not have a proof of this yet.
\ermrk

%%---------------------------------------------

\end{document}